\documentclass[12pt]{amsart}
\usepackage[colorinlistoftodos,bordercolor=orange,backgroundcolor=orange!20,linecolor=orange,textsize=scriptsize]{todonotes}

\usepackage{setspace}
\usepackage{comment}
\usepackage{enumitem}
\usepackage{multicol}
\usepackage[margin=1in]{geometry} 
\usepackage{amsmath,mathrsfs, amsthm,amssymb}
\usepackage{listings}

\usepackage{latexsym}
\usepackage{color}
\usepackage{graphics,graphicx} 
\usepackage{float}
\usepackage{subfigure}
\usepackage{epsfig} 
\usepackage{epstopdf}
\usepackage{caption,subcaption}
\usepackage{afterpage}
\usepackage{algorithm}  
\usepackage{algorithmicx}  
\usepackage{algpseudocode}  
\usepackage{tikz}
\usepackage{cite}
\usepackage{hyperref}
\usepackage{cleveref}
\usepackage{multirow}
\usepackage{braket}
\usepackage{booktabs}
\usepackage{ulem}
\usepackage{array}

\newtheorem{theorem}{Theorem}[section]

\theoremstyle{definition}
\newtheorem{definition}[theorem]{Definition}

\theoremstyle{remark}

\newtheorem{pref}{Preference function} 

\newcommand{\NO}{NeuralOp}

\newcommand{\R}{\mathbb{R}}
\newcommand{\Pp}{\mathcal{P}}
\newcommand{\F}{\mathcal{F}}
\newcommand{\M}{\mathcal{M}}

\newcommand{\Ll}{\mathcal{L}}
\newcommand{\B}{\mathcal{B}}
\newcommand{\K}{\mathcal{K}}

\newcommand{\Am}{\mathbf{A}}
\newcommand{\Mm}{\mathbf{M}}
\newcommand{\Rm}{\mathbf{R}}
\newcommand{\Pm}{\mathbf{P}}
\newcommand{\fv}{\mathbf{f}}
\newcommand{\uv}{\mathbf{u}}
\newcommand{\rv}{\mathbf{r}}

\newcommand{\xv}{\mathbf{x}}

\newcommand{\kr}{Krylov\_ite}
\newcommand{\rela}{Relaxation\_ite}

\newcommand{\Omegab}{\bar{\Omega}}

\DeclareMathOperator{\tr}{tr}

\begin{document}

\title[Meta-solvers for nonlinear PDEs]{Automatic discovery of optimal meta-solvers for time-dependent nonlinear PDEs$^{\, *}$}

\author{Yougkyu Lee$^{\, 1}$}
\author{Shanqing Liu$^{\, 1}$}
\author{Jerome Darbon}
\author{George Em Karniadakis}
\address[Youngkyu Lee, Shanqing Liu, Jerome Darbon and George Em Karniadakis]
{Division of applied mathematics, Brown University. }
\email{ \{Youngkyu\_lee, Shanqing\_liu, Jerome\_darbon, george\_karniadakis\}@brown.edu}

\thanks{$(*)$ This work is supported by the DARPA-DIAL grant HR00112490484.}
\thanks{$(1)$ Youngkyu Lee and Shanqing Liu contributed equally to this work.}

\date{\today}

\begin{abstract}

We present a general and scalable framework for the automated discovery of optimal meta-solvers for the solution of time-dependent nonlinear partial differential equations after appropriate discretization. 
By integrating classical numerical methods (e.g., Krylov based methods) with modern deep  learning components, such as neural operators, our approach enables flexible, on-demand solver design tailored to specific problem classes and objectives.
The fast solvers tackle the large linear system resulting from the Newton-Raphson iteration or by using an implicit-explicit (IMEX) time integration scheme.
Specifically, we formulate solver discovery as a multi-objective optimization problem, balancing various performance criteria such as accuracy, speed, and memory usage. 
The resulting Pareto optimal set provides a principled foundation for solver selection based on user-defined preference functions. 
When applied to problems in reaction–diffusion, fluid dynamics, and solid mechanics, the discovered meta-solvers consistently outperform conventional iterative methods, demonstrating both practical efficiency and broad applicability.

    \end{abstract}

\maketitle

\section{Introduction}

The numerical solution of partial differential equations (PDEs) is at the heart of many problems in computational science and engineering. 
Discretizing these equations typically results in large linear or nonlinear systems that must be solved efficiently and accurately. 
Common discretization techniques include the finite difference method (FDM)~\cite{strikwerda2004finite}, finite element method (FEM)~\cite{bathe2006finite} and spectral or spectral-element  methods~\cite{karniadakis2005spectral}. 
Once discretized, time-dependent and nonlinear PDEs often yield large, nonlinear algebraic systems, 
frequently solved using the Newton-Raphson method~\cite{kelley2003solving}, an iterative technique that linearizes the nonlinear system around successive approximations. 
When time integration is required, especially in stiff systems, implicit-explicit (IMEX) methods~\cite{ascher1997implicit} are commonly used, treating stiff terms implicitly for stability and non-stiff terms explicitly for efficiency. 
These discretization and integration techniques form the foundation of most modern simulation pipelines.

These computational challenges are amplified in large-scale simulations, where the number of unknowns can exceed one billion, as is often encountered in applications of industrial complexity. 
The demand for fast and scalable solvers is especially critical in multi-query contexts such as uncertainty quantification, design, optimization, and control, where the underlying PDE must be solved repeatedly for varying input parameters, such as material properties, boundary conditions, or external forces. 

Historically, iterative solvers such as Jacobi and Gauss-Seidel~\cite{greenbaum1997iterative,saad2003iterative} were among the first methods developed to solve the resulted algebraic  systems. Since the 1970s, Krylov subspace methods have become the standard for large-scale linear systems, with prominent examples including Conjugate Gradient (CG), Generalized Minimal Residual (GMRES), and Bi-Conjugate Gradient Stabilized (BiCGStab)~\cite{saad1993flexible,van2003iterative,van1992bi}. These methods remain a cornerstone of large-scale scientific computing even today.

In recent years, a paradigm shift has emerged with the rise of scientific machine learning (SciML) methods for solving PDEs~\cite{raissi2019physics}. 
Physics-informed neural networks (PINNs) embed governing equations in the loss functions of deep neural networks, allowing solutions to be learned directly from data and physical knowledge. Another growing class of techniques, known as operator learning~\cite{lu2021learning,li2020fourier}, seeks to approximate solution operators in function space, enabling generalization across input conditions. 
Variants of these approaches~\cite{karniadakis2021physics,cai2021physics,linka2022bayesian,diab2024u} seek to replace classical solvers with models based on modern machine learning (ML) techniques. 
These approaches offer several potential advantages in handling high-dimensional problems, irregular geometries, and data-driven scenarios. 
Nonetheless, these methods face challenges—most notably, spectral bias~\cite{rahaman2019spectral}, which limits their ability to capture high-frequency features. 

While hardware performance continues to advance, algorithmic improvements have often delivered even greater computational gains. A prime example is the multigrid method (see, for instance,~\cite{bramble1990parallel,wesseling2004introduction}), which surpassed Moore’s Law in terms of efficiency improvements during the 1980s and 1990s. Multigrid methods function as meta-solvers: they accelerate convergence by combining classical iterative techniques, such as Jacobi and Gauss–Seidel, with a hierarchy of discretizations operating across multiple spatial scales. Despite its conceptual elegance and reliance on well-established solvers, the multigrid method was not fully developed until nearly a century after its foundational components were introduced. This historical delay underscores a recurring theme in computational science—the significant lag between mathematical innovation and its widespread computational application.

\begin{figure}
    \centering
    \includegraphics[width=0.99\linewidth]{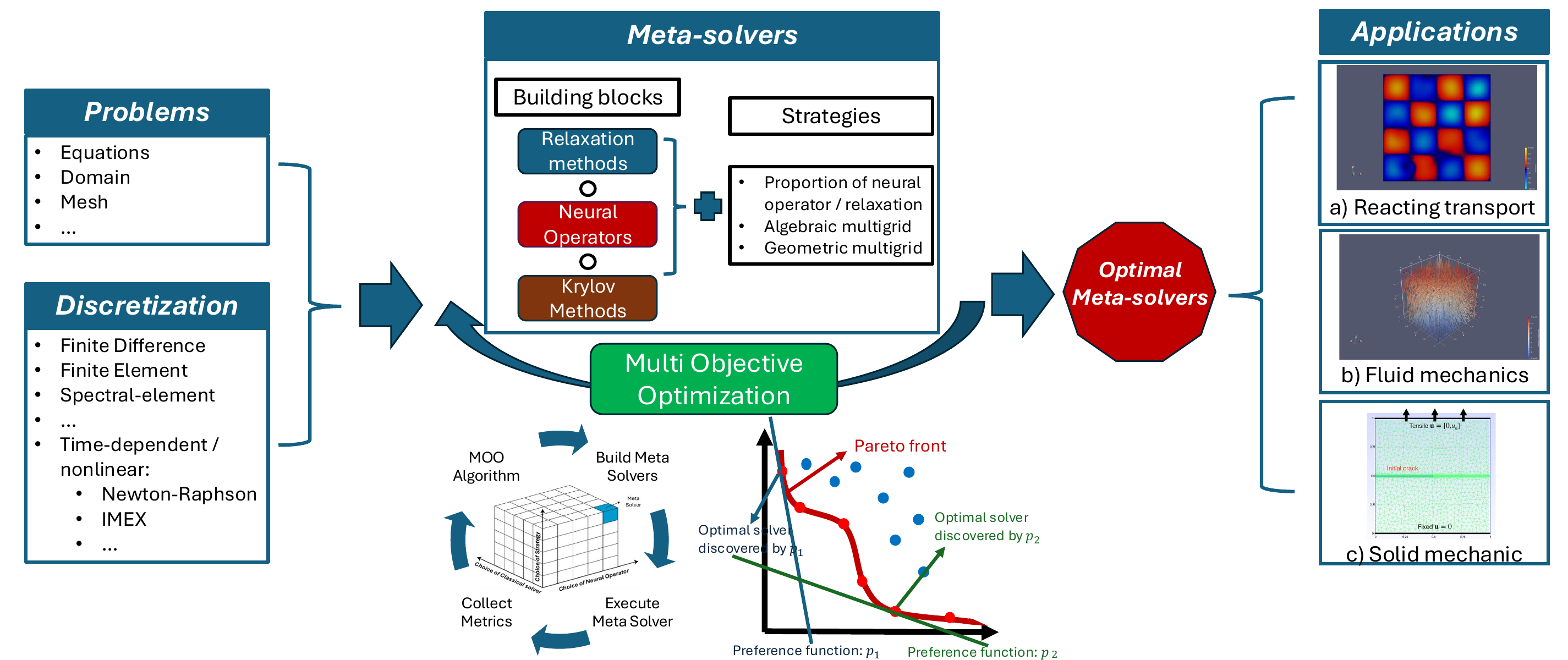}
    \caption{Workflow of the meta-solver methodology. Given a specific problem, the governing equations are first discretized using an appropriate numerical scheme, resulting in large linear or nonlinear systems. 
    To solve these systems, we explore combinations of classical solvers and neural operators, along with various adaptive strategies—yielding broad classes of candidate meta-solvers. 
    These meta-solvers are then optimized using multi-objective optimization (MOO) techniques to identify the high-dimensional Pareto front. 
    The optimal meta-solver is selected based on a user-defined preference function that encodes performance priorities. 
    The bottom figure highlights how MOO and preference functions enable the on-demand discovery of the most suitable meta-solver for a given problem.}
    \label{workflow}
\end{figure}

In this work, we propose a principled and scalable framework to accelerate the discovery of the next generation of fast solvers through the synthesis of meta-solvers. An overview of the proposed framework is illustrated in~\Cref{workflow}. 
This approach stems from the strategic integration of well-established numerical algorithms (e.g., Relaxation and Krylov methods) with emerging computational paradigms, such as neural operators. 
The resulting framework is flexible and generalizable, offering enhanced capabilities for handling non-linearities and time-dependent behavior by incorporating established discretization techniques. 
Furthermore, it is naturally compatible with advanced acceleration strategies, including geometric and algebraic multigrid methods. 
In its fundamental principle, we exploit the spectral properties of neural networks and
basic solvers to search for the optimal meta-solver on-demand for a specific problem across different applications. 

The natural question arises: Which meta-solver is optimal? 
Our findings suggest that no single meta-solver universally outperforms others across all performance metrics, such as accuracy, speed, and memory efficiency. 
To address this, we cast the problem as a multi-objective optimization (MOO) task~\cite{miettinen1999nonlinear,censor1977pareto}. 
Meta-solvers are parametrized in a high-dimensional design space, and their performance is represented as a vector-valued objective function. 
We then identify the Pareto optimal set—those meta-solvers for which no objective can be improved without degrading another. 
To support application-specific solver selection, we introduce the concept of a {\em preference function}, which encodes user-defined priorities among competing performance criteria. 
This enables us to reduce the MOO problem to a classical single-objective formulation, where the objective function is defined as the composition of the preference and performance maps. 
Moreover, our framework allows for preference inference—that is, given a known Pareto-optimal meta-solver, we can infer a corresponding preference function under which it is the optimal choice.

We demonstrate our framework on benchmark problems from diverse application domains, spanning reacting transport (reaction–diffusion equations), fluid mechanics (Navier–Stokes equations), and solid mechanics (brittle fracture). In all cases, the discovered meta-solvers consistently demonstrate significant performance gains over standard Krylov-based methods, showcasing the practical effectiveness and adaptability of our framework.

The organization of this paper is as follows. \Cref{sec-meth} presents our methodology. Specifically, \Cref{subsec-kry} reviews trunk-basis hybrid preconditioners combining Krylov methods and neural operators; \Cref{subsec-timepde} introduces Newton-Raphson and IMEX-based meta-solvers for time-dependent nonlinear PDEs; and \Cref{subsec-moo} discusses multi-objective optimization, Pareto optimality, and preference-based discovery of optimal meta-solvers. Applications and numerical results are provided in \Cref{sec-num}, including reaction-diffusion equations (\Cref{subsec-rd}), Navier–Stokes equations (\Cref{subsec_ns}), and brittle fracture problems (\Cref{subsec-bf}). \Cref{sec-sum} summarizes our contributions and outlines future directions. 
Finally several appendices provides implementation details and further numerical results 
for time-dependent reaction-diffusion equations (\Cref{app_rd}), 
for incompressible Navier-Stokes equation (\Cref{app_ns}) and
for the brittle fracture fracture problem (\Cref{app-bf}).

\section{Methodology}\label{sec-meth}

\subsection{Krylov-based hybrid preconditioners for linear systems}\label{subsec-kry}

\begin{figure}
    \centering
    \includegraphics[width=0.99\linewidth]{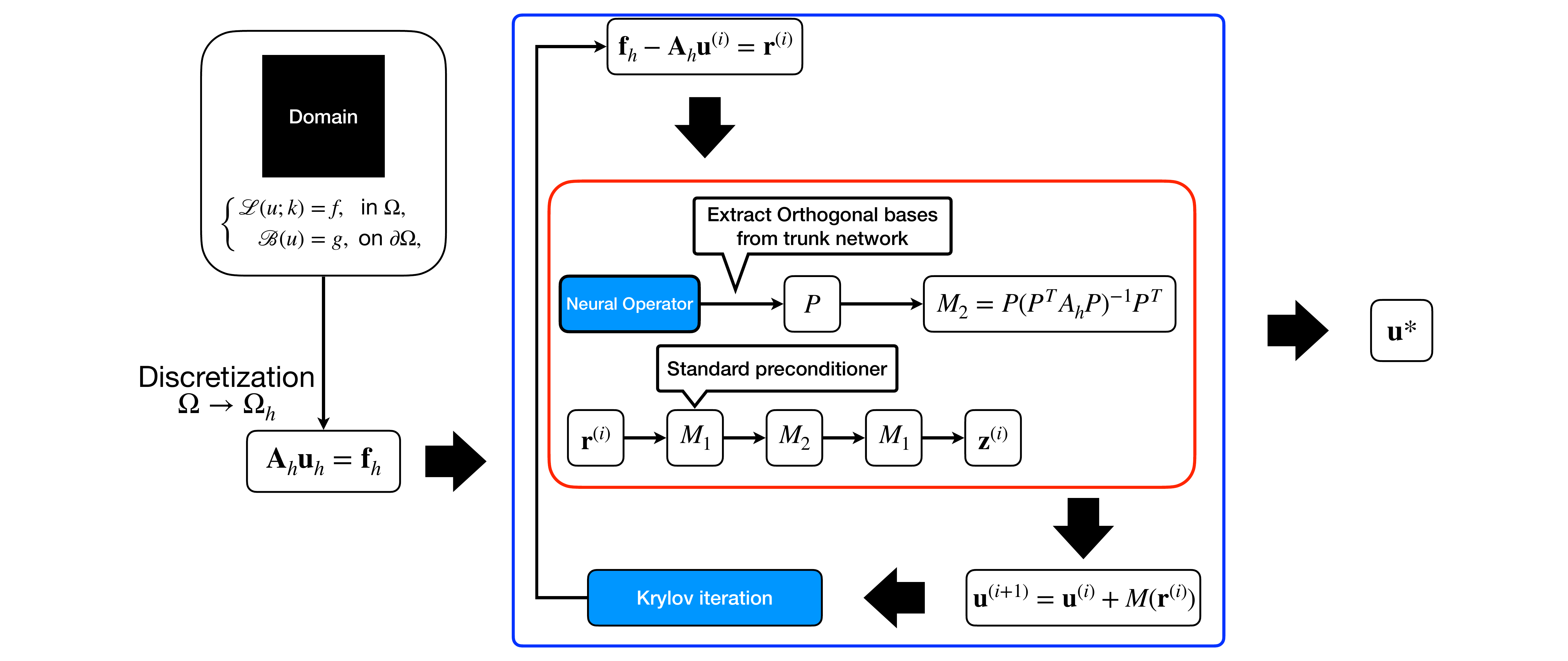}
    \caption{Krylov-based hybrid meta-solvers. These meta-solvers are constructed by combining different neural operators with various Krylov methods. 
In addition, several steps of standard relaxation methods (e.g., Jacobi, Gauss-Seidel) are applied before and after the neural operators, serving as smoothers. Advanced techniques, such as multigrid methods, can also be integrated by replacing classical relaxation steps with combinations of neural operators and relaxation methods.}
    \label{fig_krylov}
\end{figure}
In this section, we recall the Krylov-based hybrid preconditioners for linear systems. 
A schematic overview is provided in~\Cref{fig_krylov}. 
We are interested in numerically approximating the solution $u : \Omegab \to \R$ of the linear differential equation of the form
\begin{equation}\label{eq:linear}
\left\{
\begin{aligned}
& \Ll(u) = f, \ &\text{in } \Omega \ , \\
& \B(u) = g, \ &\text{in } \partial \Omega \ ,
\end{aligned}
\right.
\end{equation}
where $\Omega \subset \R^d$ is an open bounded domain, $\Ll : \R^{\Omega} \to \R^{\Omega}$ is a linear differential operator.
$f:\Omega \to \R$ is a known function that does not explicitly involve $u$. $\B :\R^{\partial \Omega} \to \R^{\partial \Omega}$, $g:\partial \Omega \to \R$ serve as the boundary condition. 
Given a mesh $\Omega^h$ discretizing $\Omega$, a common approach to numerically approximate the solution of system~\eqref{eq:linear} is the 
finite element method, which can be formulated as a linear system on $\Omega^h$, of the form
\begin{equation}
\label{eq:fem}
    \Am \uv = \fv \ ,
\end{equation}
where $\uv$ denotes the nodal coefficients of finite element basis functions. 
The discretized system~\eqref{eq:fem} can be solved using iterative solvers. 
In particular, given an appropriate initialization $\uv^{(0)}$, 
let $\uv^{(i)}$ be the approximate solution after $i$ iterations.  
At the $(i+1)$-th step, the iteration has the form 
\begin{equation}\label{eq::ite}
\left\{ 
\begin{aligned}
& \rv^{(i)} = \fv - \Am \uv^{(i)} \ , \\
& \uv^{(i+1)} = \K (\rv^{(i)}) \ ,
\end{aligned}
\right.
\end{equation}
where $\K$ denotes the update process, which depends on the iterative solver. 

Another approach solving the problem~\eqref{eq:linear} is through the use of neural operators~\cite{kovachki2021neural,lu2021learning}, which approximate the functional operator that maps the function $f$ to the solution $u$ using neural networks. 
However, approximating this operator at the fine level is challenging due to the spectral bias of neural networks. 
To address this difficulty, 
hybrid approaches have been introduced in~\cite{kopanivcakova2024deeponet,lee2024fast, zhang2024blending}, where neural operators are used to define preconditioners for iterative solvers, rather than using neural operators as direct solvers for~\eqref{eq:linear}. 
In this context, 
the preconditioning process is formulated as follows
\begin{equation}
\left\{
    \begin{aligned}
        &\rv^{(i)} = \fv - \Am \uv^{(i)} \ , \\
        &\uv^{\ast} = \uv^{(i)} + \Mm_{1}(\rv^{(i)}) \ , \\
        &\rv^{\ast} = \fv - \Am \uv^{\ast} \ , \\
        &\uv^{\ast\ast} = \uv^{\ast}+ \Mm_{2}(\rv^{\ast}) \ , \\
        &\rv^{\ast\ast} = \fv - \Am \uv^{\ast\ast} \ , \\
        &\uv^{\ast\ast\ast} = \uv^{\ast\ast} + \Mm_{1}(\rv^{\ast\ast}) \ , \\
        &\rv^{\ast\ast\ast} = \fv - \Am \uv^{\ast\ast\ast} \ , \\
        &\uv^{(i+1)}   = \K (\rv^{\ast\ast\ast}) \ ,
    \end{aligned}
\right.
\label{numerical_scheme}
\end{equation}
where $\Mm_{1}$ and $\Mm_{2}$ denote the steps of relaxation method and the inference through the pre-trained neural operator, respectively. 
The update process~\eqref{numerical_scheme} is typically well defined 
for the hybrid preconditioned relaxation methods as in~\cite{zhang2024blending}. 
However, when the preconditioner $\Mm_{2}$ is defined as inference via the neural operator, it cannot be directly used as the preconditioner for Krylov methods.
To address this issue, when DeepONet~\cite{lu2021learning} is utilized as the neural operator, the trunk-basis~(TB) approach, proposed in~\cite{kopanivcakova2024deeponet}, extracts the prolongation operator~$\Pm$ and restriction operator~$\Rm$ from the trunk network of DeepONet. 
This leads to the construction of the second preconditioner
\begin{equation*}
\Mm_{2}:=\Pm(\Rm\Am\Pm)^{-1}\Rm=\Pm\Am_{c}^{-1}\Rm.
\end{equation*}

The $(i,j)$-th component of the matrix of the prolongation operator~$\Pm$ is computed by
\begin{equation}
    [\Pm]_{ij} = T_{j}(\xv_{i}) \ ,
\end{equation}
where $T_{j}(\xv_{i})$ denotes the $j$-th component of the output of the trunk network evaluated at the coordinate $\xv_{i} \in \Omega$.
Note that the restriction operator is defined as the adjoint operator of the prolongation operator, i.e., $\Rm = \Pm^{T}$.
The quality of the prolongation operator $\Pm$ can be further improved utilizing sampling strategies and QR decomposition (see~\cite{kopanivcakova2024deeponet}).

\subsection{Time-dependent nonlinear equations: Newton-Raphson and Implicit-Explicit (IMEX) methods } \label{subsec-timepde}
Consider a general time-dependent nonlinear equation of the form
\begin{equation}
\frac{\partial u}{\partial t} = Lu + G(u) =: F(u) \ ,
\label{eqn:abstract}
\end{equation}
where $Lu$ and $G(u)$ denote the linear and non-linear term, respectively. 
In this section, we propose two meta-solver approaches that integrate neural operators with two prominent classical methods: the \textit{Newton-Raphson} method and the \textit{Implicit-Explicit (IMEX)} method, to enhance the numerical approximation of problem~\eqref{eqn:abstract}. 
Additionally, we introduce a discovery mechanism designed to identify the optimal meta-solver based on the specific problem at hand. 

Before explaining the Newton-Raphson method and the IMEX method, we first perform time discretization.
Note that the Crank–Nicolson scheme~\cite{Crank1996} is used, resulting in the following nonlinear system at each time step:
\begin{equation}
\label{eqn:steady}
\frac{u^{(n+1)} - u^{(n)}}{\Delta t} = \frac{1}{2} F ( u^{(n+1)} ) + \frac{1}{2} F ( u^{(n)} ),
\end{equation}
where $\Delta t$ and $u^{(k)}$ denote the time step size and the solution at $k-th$ time step. 

For the Newton-Raphson method, we rearrange~\eqref{eqn:steady} to obtain the following:
\begin{equation}
\label{eqn:newton}
\F(u^{(n+1)}): =u^{(n+1)} - \frac{1}{2}\Delta t F ( u^{(n+1)} ) - \frac{1}{2}\Delta t F ( u^{(n)} ) - u^{(n)}  = 0.
\end{equation}
The next iterate $u^{(n+1)}$  is then computed by solving $\mathcal{F}(u) = 0$ using the Newton-Raphson method.

For the IMEX method, we apply explicit scheme to the nonlinear term $G(u)$ and implicit scheme to the linear term $Lu$.
We utilize the third-order Adams-Bashforth scheme~\cite{bashforth1883attempt} for explicit scheme and the Crank-Nicolson scheme for implicit scheme, respectively.
The resulting system has the form: 
\begin{equation}
\underbrace{u^{(n+1)} - \frac{1}{2}\Delta t L u^{(n+1)}}_{(M - \frac{1}{2}\Delta t A)U^{(n+1)}}  = 
\underbrace{u^{(n)} + \frac{1}{2}\Delta t L u^{(n)}+\Delta t AB3(u^{(n)}, u^{(n-1)},u^{(n-2)})}_{b^{(n)}},
\label{eqn:time-imex}
\end{equation}
where $U^{(n)}$ and $AB3$ denote the coefficients for finite element approximation at the $n$-th time-step and the 
third-order Adams-Bashforth operator, respectively.
Here, $M$ and $A$ are the mass matrix and the stiffness matrix, respectively.
The third-order Adam-Bashforth operator is given by
\begin{equation}
AB3(u^{(n)}, u^{(n-1)}, u^{(n-2)})=
\left\{
\begin{aligned}
&G(u^{(n)}), &\text{ if } n=0, \\
&\frac{3}{2}G(u^{(n)}) -\frac{1}{2}G(u^{(n-1)}), &\text{ if } n=1, \\
&\frac{23}{12}G(u^{(n)})-\frac{16}{12}G(u^{(n-1})+\frac{5}{12}G(u^{(n)}) &\text{ if } n\ge 2.
\end{aligned}
\right.
\end{equation}
Since the matrix $(M - 0.5\Delta t A)$ is linear, we can compute $U^{(n+1)} = (M - 0.5\Delta t A)^{-1} b^{(n)}$.

Both of the aforementioned methods are classified as linearization techniques, as they transform the original nonlinear problem~\eqref{eqn:abstract} into a sequence of linear systems of the form $A u = f$. 
These systems can be solved using an iterative solver, where the iterations take the form
\begin{equation}\label{eq:ite}
\left\{ 
\begin{aligned}
& r^{(i)} = f - A u^{(i)} \ , \\
& u^{(i+1)} = u^{(i)} +  \K (r^{(i)}) \ ,
\end{aligned}
\right.
\end{equation}
with $\K$ denotes the update process. 
To solve this iterative system, 
we construct the Krylov method based meta-solvers, which are incorporated within the trunk basis hybridization approaches (see a sketch in~\Cref{fig_krylov}).

\subsection{Multi-objective optimization, Pareto optimality and preference function based discovery mechanisms}\label{subsec-moo}

As a result of above construction, the family of trunk-based hybridization preconditioners can be parameterized by a five-dimensional space $\M$, where each dimension corresponds to the choice of neural operator, the choice of Krylov method, the choice of relaxation method, the strategy of applying smoothers, and the level of multi-grid techniques used. 
Thus, by further combining them with the Newton-Raphson method and the IMEX method, we construct two families of meta-solvers for the time-dependent nonlinear system~\eqref{eqn:abstract}. 

Next, we present a preference function-based mechanism for the discovery and rediscovery of optimal meta-solvers. In both families of the aforementioned meta-solvers, performance can be modeled and quantified by a vector-valued function 
 $f:\M \to \R^d$, where for every $x\in \M$, the components $\{f_i(x)\}_{i\in\{1,\dots,d\}}$ represent different performance criteria—such as computational time, relative error, number of iterations, and so on. 
Finding an optimal meta-solver can thus be formulated as a multi-objective optimization problem of the form
\begin{equation}\label{moo}
    \min_{x\in \M} f(x) = [f_1(x),\dots,f_d(x)]^T \ .
\end{equation}

This optimization problem typically does not admit a classical solution that minimizes all criteria simultaneously. Therefore, we instead consider a weak notion of optimality, Pareto optimality, and focus on identifying Pareto optimal solutions (see also~\cite{lee2024automatic,cao2025automatic}).
\begin{definition}
\begin{enumerate}
    \item For every $x^1,x^2\in \M$, we call $x^1$ dominates $x^2$ if $f_i(x^1) \leq f_i(x^2)$ for every $i \in \{1,2,\dots, d\}$ and there exists at least one $j$ such that $f_j(x^1)<f_j(x^2)$. We denote $x^1 \succeq x^2$  $x^1$ dominates $x^2$. 
    \item We say that $x \in \M$ is a Pareto optimal solution of problem~\eqref{moo} if there is no other element $x' \in \M$ such that $x' \succeq x$. We denote $\Pp_f(\M)$ the set of Pareto optimal solutions of problem~\eqref{moo}. 
    \item We call the image of the set of Pareto optimal solutions $\Pp_f(\M)$ by  objective function $f$ the Pareto front $\F_f(\M)$, that is, $\F_f(\M) = \{f(x) \in \R^d \mid x \in \Pp_f(\M)\}$.
\end{enumerate}
\end{definition}

We next present the preference function based methodology to discover the optimal meta-solvers, in various of context. 
First, the performance data for various criteria is evaluated using a re-scaling function. Note that other normalization methods can also be applied to the performance data, that map the original value in each dimension to a common subset of $\R$, allowing a unified evaluation. Here, without loss of generality, we use a re-scaling function. 
In particular, considering the MOO problem~\eqref{moo}, for each dimension $i\in\{1,2,\dots,d\}$, denote
\begin{equation}
\overline{f}_i = \max_{x \in \Pp_f(\M)} f_i(x), \quad \underline{f}_i = \min_{x \in \Pp_f(\M)} f_i(x) \ ,
\end{equation}
where $\Pp_f(\M)$ denotes the set of Pareto optimal solvers of MOO problem~\eqref{moo}, and we assume $\overline{f}_i \neq \underline{f}_i$ for every $i$. 
Then, for a solver $x^k \in \Pp_f(\M)$, its performance is rescaled to
\begin{equation}\label{rescale}
f'_i(x^k) = \frac{f_i(x^k) - \underline{f}_i}{\overline{f}_i - \underline{f}_i} \ .
\end{equation}
User's preferences on different criteria are modeled by a preference function $p:\R^d \to \R$, which maps the image set of performance function to $\R$. 
Then, discovery the optimal meta-solver under this preference is formulated as a classical optimization problem, of the form
\begin{equation}
    \min_{x \in \M} p \circ f(x) \ .
\end{equation}
Here as a 
 prototypical example, we consider the weighted sum function of the rescaled performance. 
 Given a weight $\lambda = (\lambda_1,\lambda_2,\dots,\lambda_d)$ such that $0 \leq \lambda_i \leq 1, \forall \ i \in\{1,2,\dots, d\}$ and $\sum_{i=1}^d\lambda_i = 1$, the preference function $p(\lambda; \cdot):\R^d \to \R$ is defined by
\begin{equation}\label{weight_sum}
f'=(f'_1,f'_2,\dots,f'_d) \mapsto p(\lambda;f') := \sum_{i=1}^d \lambda_i f'_i = \lambda^Tf' \ . 
\end{equation}
Given different weights $\lambda$, the preference function is used to select a solver among the Pareto optimal ones. 
Geometrically, each weighted sum preference function corresponds to a $(d-1)$-dimensional hyperplane, where the weights determine the orientation of the hyperplane. One can incrementally decrease the constant term, starting from $0$, until the hyperplane intersects the Pareto front. The point of intersection corresponds to the performance of the optimal solver being selected, and the constant term represents the inverse of the weighted sum.

The re-discovery of one particular (parameterized) solver $x^j \in \Pp_{f'}(\M)$ is equivalent to finding a weight $\lambda^j$ such that $p(\lambda;f'(x^j))$ is the minimum for every $x \in \Pp_{f'}(\M)$. This can be formulated as a linear programming problem, that has the form:
\begin{equation}\label{lp_discover}
\begin{aligned}
&\min_{\lambda} \ \lambda^T f'(x^j) \\
&\text{ s.t. }
\left\{ 
\begin{aligned}
& \lambda^T f'(x^j) \leq \lambda^T f'(x^{-j}),\ \forall \ x^{-j} \in \Pp_{f'}(\M) \ ,\\
& 0 \leq \lambda_i \leq 1 \text{ and } \sum \lambda_i=1 \ .
\end{aligned}
\right.
\end{aligned}
\end{equation}
It is geometrically equivalent to finding the tangent hyperplane of the Pareto front at the point $f'(x^j)$. Note that the LP problem~\eqref{lp_discover} can be solved efficiently using modern LP solvers (e.g., CPLEX). 

\section{Numerical Results}\label{sec-num}

In the following, we apply the two families of meta-solvers to a series of time-dependent nonlinear equations arising from reacting flows, fluid mechanics, and solid mechanics. For each case, we discover optimal meta-solvers under various preference functions and compare their performance against standard (vanilla) iterative solvers. 

\subsection{Solving time-dependent reaction-diffusion equations}\label{subsec-rd} 
\subsubsection{The model and implementation details}
In this section, we apply our methodology for solving the time-dependent reaction-diffusion equation in a two-dimensional domain. In particular: Given the domain $\Omega = [0, 1]^{2} \subset \R^{2}$ and time $T=[0,1]$, coefficient $k(\mathbf{x})$, and force term $ f(\mathbf{x})$, the 2-d time-dependent reaction-diffusion equation is given by
\begin{equation}
\left\{\begin{split}
\frac{\partial u}{\partial t} &= \nabla \cdot (k \nabla u) + R(u) + f, \,\,\, \text{ in } \Omega \times (0,1], \\
u &= u_{0}, \text{ in } t=0, \\
u &= 0, \text{ in } \partial\Omega \times (0,1],
\end{split}
\right.
\label{eqn:time-reaction-diffusion}
\end{equation}
where $R(u) = u - u^{2}$~(Fisher type)~\cite{Fisher1937THEWO}. 
To show the generalizability of our methodology, we consider two cases of different coefficient $k(\mathbf{x})$ and form term $f(\mathbf{x})$. 
In the first case, we take the diffusion coefficient to satisfy the following:
\begin{subequations}
\begin{equation}
k \sim N(1.0, K(\mathbf{x},\mathbf{x}')=0.3e^{-\frac{\Vert \mathbf{x}-\mathbf{x}'\Vert^{2}}{2 \times 0.1^{2}}}), \text{ and } k \geq 0.1 \ .
\end{equation}
For the external force $f$, we take
\begin{equation}
    f \sim N(0.0, K(\mathbf{x},\mathbf{x}')=e^{-\frac{\Vert \mathbf{x}-\mathbf{x}'\Vert^{2}}{2 \times 0.1^{2}}}).
\end{equation}
\end{subequations}
The second case corresponds to smaller correlation length in the system. 
In particular, we take the new diffusion coefficient to satisfying the following: 
\begin{subequations}
\begin{equation}
k \sim N(10.0, K(\mathbf{x},\mathbf{x})=3.0e^{-\frac{\Vert \mathbf{x}-\mathbf{x}'\Vert^{2}}{2 \times 0.01^{2}}}), \text{ and } k \geq 1.0 \ .
\end{equation}
\end{subequations}
The other implementation parameters remain the same as the previous equation. 
In the following section, we present the discovery and rediscovery of IMEX-based optimal meta-solvers for problems with large correlation lengths, and Newton–Raphson based optimal meta-solvers for small correlation lengths. Comprehensive numerical results—including the composition of all Pareto-optimal meta-solvers, various 3D projections of the Pareto fronts, and the discovery and rediscovery results for both approaches in both approaches—are provided in~\Cref{app_rd}.

\subsubsection{Discovery and re-discovery of optimal meta-solvers}

For the parametrization of the meta-solvers, we consider 
\begin{equation}
    \M = \NO \times \kr \times \rela \times S \times \{1,2,3\} \ ,
\end{equation}
such that $\NO= \{\text{DeepONet, U-DeepONet, KAN, JacobiKAN, CheyKAN}\}$, $\kr  = \{ \text{FGMRES, CG, BiCGStab} \}$, $\rela = \{\text{Jacobi, Gauss-Seidel, SOR, SSOR}\}$, strategies of applying relaxation method $S = \{\text{1-1-1, 3-1-3, 5-1-5, 7-1-7, 9-1-9}\}$ and the levels of multi-grid is chosen from $\{1,2,3\}$. The evaluation is then conducted in a six-dimensional performance function $f = (f_1,\dots, f_6)$, where each dimension represents a particular performance criterion. In particular, $f_1$ is the relative error, $f_2$ is the computational time, $f_3$ is the number of iterations, $f_4$ is the memory allocation, $f_5$ is the MACs and $f_6$ is the training time, respectively.

We implement the preference function based method to discover the optimal solver. In particular, the performance data at each dimension are first rescaled to a real value in $[0,1]$, and the preference functions are applied in the rescaled data. Here as a prototypical example, we applied the weighted sum preference functions.

\begin{figure}
    \centering
    \begin{minipage}[t]{0.45\textwidth}
    \centering
        \vspace{1.5cm}
        \subfigure[3-d projection of Pareto front and optimal meta-solvers discovered by preference functions.]{
            \centering
            \includegraphics[width=\linewidth]{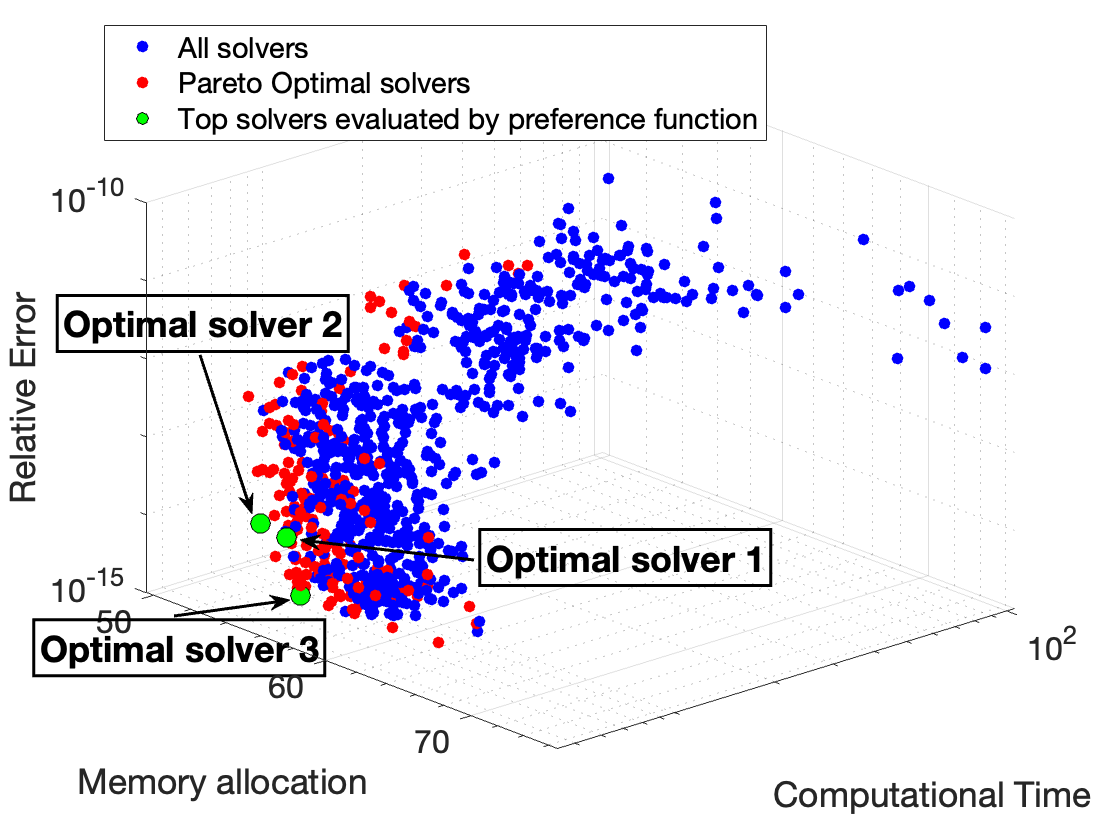}}
    \end{minipage}%
    \hfill
    \begin{minipage}[t]{0.55\textwidth}
        \begin{table}[H]
    \centering
     \begin{minipage}{\linewidth}
        \centering
        \captionof{subtable}{\small The top-3 meta-solvers.}
        \resizebox{\linewidth}{!}{
        \begin{tabular}{c|c|c|c|c|c}
\hline \hline 
 &Neural operator & Classical solver & Multi-grid & Relaxation & Strategies \\ 
\hline \hline
     Top-1 & U-Net & BiCGStab & 2-level & Gauss-Seidel & 3-1-3 \\ \hline
     Top-2 & DeepONet & BiCGStab & 2-level & Gauss-Seidel &  3-1-3 \\ \hline
     Top-3 & U-Net & BiCGStab & 2-level & Gauss-Seidel & 5-1-5 \\ \hline
\end{tabular}
}
    \end{minipage}

    \vspace{1em} 

    \begin{minipage}{\linewidth}
        \centering
        \captionof{subtable}{\small Performance of the top-3 solvers.}
        \resizebox{\linewidth}{!}{ 
          \begin{tabular}{c|c|c|c|c|c|c}
\hline \hline 
 &Relative error & Com. time  & \# of ite & Memory & MACs & Training time \\ 
\hline \hline
Top-1 & $1.64 \times 10^{-14}$ & 5.11  & 384 & 56.39 & $2.55 \times 10^{11}$ & 10973  \\ \hline
Top-2& $1.62 \times 10^{-14}$ & 5.34  & 384 & 54.46  & $2.21 \times 10^{11}$ & 17235.9 \\ \hline
Top-3& $2.17 \times 10^{-15}$ & 5.99  & 384 & 55.87 & $2.56 \times 10^{11}$ & 10973  \\ \hline
\end{tabular}
        }
    \end{minipage}
    \begin{minipage}{\linewidth}
        \centering
        \captionof{subtable}{\small Performance of vanilla method.}
        \scriptsize
          \begin{tabular}{c|c|c}
\hline \hline 
 & Com. time  & \# of ite  \\ 
 \hline \hline
 CG & 127.5 & 49024 \\ \hline
 FGMRES & 189.6 & 70272 \\ \hline
 BiCGStab & 77.6 & 29184 \\ \hline
\end{tabular}
        
    \end{minipage}
    \caption{Composition and performance of optimal meta-solvers}
    \label{table_rd_1}
\end{table}
    \end{minipage}
    \caption{\small IMEX based meta-solvers for solving the 2-d time-dependent nonlinear reaction-diffusion equation. (a) the projection of the 6-dimensional Pareto front into three-dimensional, time-error-memory, space, with Pareto optimal meta-solvers marked in red and the optimal meta-solvers discovered by preference function are highlighted in green. Tables: (A) the parameters of the top-3 meta-solvers discovered by preference function, (B) the 6-dimensional performance of the top-3 meta-solvers, (C) comparing the computational time and number of iteration for the vanilla Newton-Raphson method, using different Krylov methods as iterative solvers. Up to an error of machine precision, our meta-solvers obtain $\approx15$ times speedup in computational time and $>76$ times speedup in number of iterations.}
\label{fig_rd_1}
\end{figure}

In~\Cref{fig_rd_1}, we present the results obtained using IMEX based meta-solvers  for problems with large correlation lengths. 
In particular,~\Cref{fig_rd_1} (a) shows a three-dimensional projection of the Pareto front. For this analysis, we define the preference function to emphasize the computational time and relative error, 
that is \[p^1_{rd}(f) =  \frac{1}{24} \big( 10 (f_1 + f_2) +  (f_3 + f_4 + f_5 + f_6 ) \big) \ . \] 
The performance of the top three meta-solvers is highlighted with green dots. 
In~\Cref{table_rd_1} (A), we present the composition of the top-3 meta-solvers discovered by the preference function $p_{rd}^1$, and their performance is shown in \Cref{table_rd_1} (B).

\begin{figure}
    \centering
    \begin{minipage}[t]{0.45\textwidth}
    \centering
        \vspace{1.5cm}
        \subfigure[3-d projection of Pareto front and optimal meta-solvers discovered by preference functions.]{
            \centering
            \includegraphics[width=\linewidth]{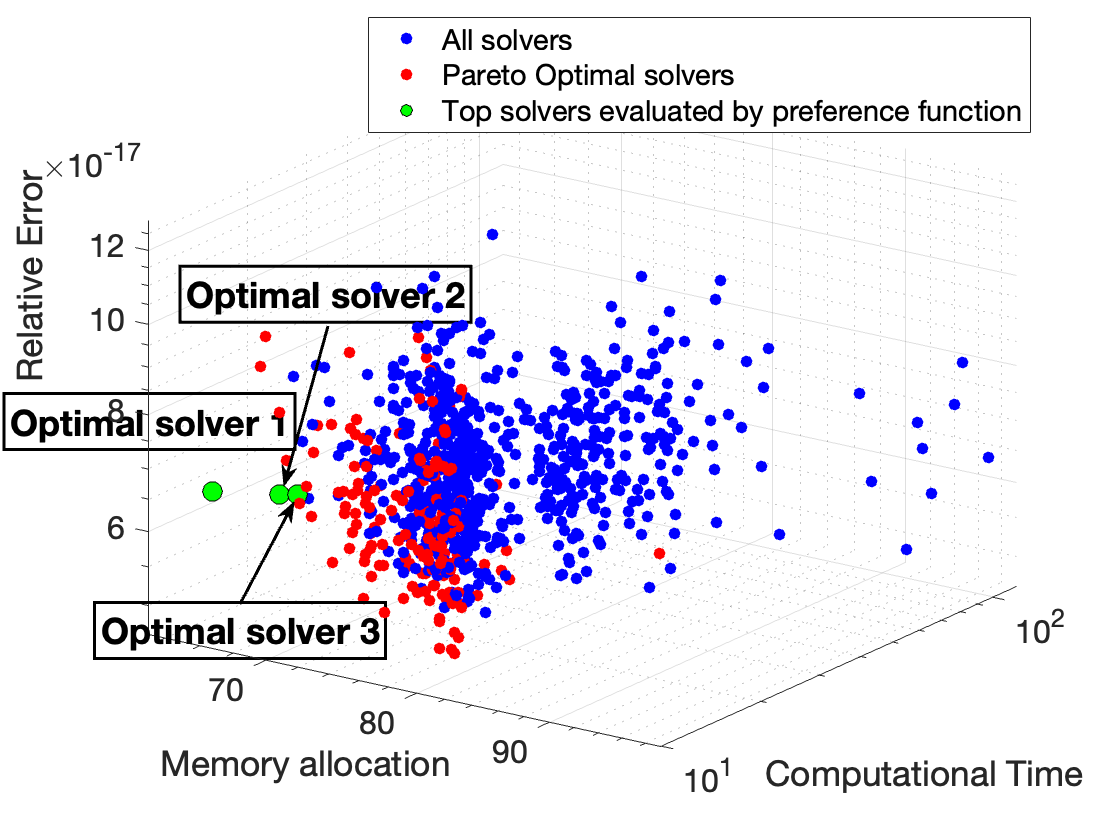}}
    \end{minipage}%
    \hfill
    \begin{minipage}[t]{0.55\textwidth}
        \begin{table}[H]
    \centering
     \begin{minipage}{\linewidth}
        \centering
        \captionof{subtable}{\small The top-3 meta-solvers.}
        \resizebox{\linewidth}{!}{
        \begin{tabular}{c|c|c|c|c|c}
\hline \hline 
 &Neural operator & Classical solver & Multi-grid & Relaxation & Strategies \\ 
\hline \hline
     Top-1 &  U-DeepONet & BiCGStab & 3-level & Jacobi & 3-1-3 \\ \hline
     Top-2 & U-DeepONet & BiCGStab & 2-level & Jacobi &  5-1-5 \\ \hline
     Top-3 & U-DeepONet & BiCGStab & 3-level & Gauss-Seidel & 7-1-7 \\ \hline
\end{tabular}
}
    \end{minipage}

    \vspace{1em} 

    \begin{minipage}{\linewidth}
        \centering
        \captionof{subtable}{\small Performance of the top-3 solvers.}
        \resizebox{\linewidth}{!}{ 
          \begin{tabular}{c|c|c|c|c|c|c}
\hline \hline 
 &Relative error & Com. time  & \# of ite & Memory & MACs & Training time \\ 
\hline \hline
Top-1& $6.16\times 10 ^{-17}$ & 15.71 & 520 & 63.09 & $3.95\times 10^{11} $   & 9010.38 \\ \hline
Top-2&  $6.21 \times 10^{-17}$ & 17.15 & 512 & 66.14 & $3.50 \times 10^{11}$  & 9010.38  \\ \hline
Top-3&  $6.38\times 10^{-17}$ & 15.91 & 384 & 67.87  & $2.96 \times 10^{11} $& 9010.38 \\\hline
\end{tabular}
        }
    \end{minipage}
    \begin{minipage}{\linewidth}
        \centering
        \captionof{subtable}{\small Performance of vanilla method.}
        \scriptsize
          \begin{tabular}{c|c|c}
\hline \hline 
 & Com. time  & \# of ite  \\ 
 \hline \hline
 CG & 260.3 & 71887 \\ \hline
 FGMRES & 426.4 & 128000 (maximum) \\ \hline
 BiCGStab & 183.2 & 50893 \\ \hline
\end{tabular}
        
    \end{minipage}
    \caption{Composition and performance of optimal meta-solvers}
    \label{table_rd}
\end{table}
    \end{minipage}
    \caption{\small Newton-Raphson based meta-solvers for solving the 2-d time-dependent nonlinear reaction-diffusion equation with small correlation length. Figures: (a) the projection of the 6-dimensional Pareto front into three-dimensional, time-error-memory, space, with Pareto optimal meta-solvers marked in red and the optimal meta-solvers discovered by preference function are highlighted in green. Tables: (A) the parameters of the top-3 meta-solvers discovered by preference function, (B) the 6-dimensional performance of the top-3 meta-solvers, (C) comparing the computational time and number of iteration for the vanilla Newton-Raphson method, using different Krylov methods as iterative solvers. Up to an error of machine precision, our meta-solvers obtain $\approx11.7$ times speedup in computational time and $>100$ times speedup in number of iterations.}
\label{fig_rd}
\end{figure}

In~\Cref{fig_rd}, we present the results obtained using IMEX based meta-solvers. 
In particular,~\Cref{fig_rd} (a) shows a three-dimensional projection of the Pareto front. For this analysis, we define the preference function as the average of all performance criteria, 
that is \[p_{rd}^2(f) = \frac{1}{6} \sum_{i=1}^6f_i \ . \] 
The performance of the top three meta-solvers is highlighted with green dots. 
In~\Cref{table_rd} (A), we present the composition of the top-3 meta-solvers discovered by the preference function $p_{rd}^1$, and their performance is shown in \Cref{table_rd} (B). 

It is worth emphasizing that our preference function based approach is both broad and flexible, allowing users to define and apply their own preference functions in order to identify the most suitable (i.e., optimal) meta-solver for their specific objectives. Conversely, this approach also enables the discovery of a preference function under which a particular meta-solver becomes the optimal choice. 
As an illustration,  we present below various weighted-sum preference functions along with the corresponding optimal meta-solvers identified by each, for the IMEX baed meta-solvers solving reaction-diffusion problen with small correlation length. 
The weights are presented in the following order: $\lambda_1$ for relative error, $\lambda_2$ for computational time, $\lambda_3$ for number of iterations, $\lambda_4$ for memory allocation, $\lambda_5$ for MACs, and $\lambda_6$ for training time. 
\begin{subequations}
\begin{equation}
\left\{
\begin{aligned}
&\text{Preference function $p_{rd}^a(f) = (0.264, 0.443, 0, 0.094, 0.072, 0.127)^T f$  } \\
&\text{Optimal solver}: x^a = (\text{U-Net}, \text{FGMRES}, \text{2-level}, \text{SSOR}, 5-1-5) \ .
\end{aligned}
\right.
\end{equation}

\begin{equation}
\left\{
\begin{aligned}
&\text{Preference function $p_{rd}^b(f) = (0.240, 0, 0.003, 0.036, 0.368, 0.354)^T f$ } \\
&\text{Optimal solver}: x^b = (\text{JacobiKAN}, \text{BiCGStab}, \text{2-level}, \text{SSOR}, 7-1-7) \ .
\end{aligned}
\right.
\end{equation}

\begin{equation}\label{rd_rediscovery}
\left\{
\begin{aligned}
&\text{Preference function $p_{rd}^c(f) = (0.043, 0.771, 0, 0.030, 0.156, 0)^T f$ } \\
&\text{Optimal solver}: x^b = (\text{KAN}, \text{FGMRES}, \text{3-level}, \text{SSOR}, 3-1-3) \ .
\end{aligned}
\right.
\end{equation}

\end{subequations}

Take, for example, the preference function and its corresponding optimal meta-solver discovered in~\eqref{rd_rediscovery}. If a user places a strong emphasis on computational time ($\approx77\%$ weight), and MACs ($\approx16\%$ weight), then the optimal meta-solver for this preference profile consists of a Neural Operator KAN combined with the Krylov method FGMRES, using SSOR as the smoother. The smoother is applied with the 3-1-3 strategy and leverages a three-level multigrid technique.

\subsection{Solving incompressible Navier-Stokes equations}\label{subsec_ns} 
\subsubsection{The model and implementation details}

In this section, 
we consider a vector-valued nonlinear equation, the 3-dimensional incompressible Navier-Stokes equation.
Given domain $\Omega = [0,1]^3 \subseteq \R^{3}$, the incompressible Naiver Stokes equation is given by
\begin{equation}
\left\{\begin{split}
\rho \frac{\partial \mathbf{u}}{\partial t} + \rho (\mathbf{u} \cdot \nabla)\mathbf{u}&= -\nabla p + \mu \Delta \mathbf{u}, \,\,\, \text{ in } \Omega \times (0,1], \\
\nabla \cdot \mathbf{u} &= 0, \,\,\, \text{ in } \Omega \times (0,1], \\
\mathbf{u} &= \mathbf{u}_{0}, \text{ at } t=0, \\
\mathbf{u} &= \mathbf{u}_{D}, \text{ on } \partial\Omega \times (0,1],
\end{split}
\right.
\label{eqn:imcompressible-navier-stokes}
\end{equation}
where $\rho = \mu = 1$.
For the reference solution, we used the fully 3D Navier-Stokes solution proposed in~\cite{ethier1994exact}, that is, 
\begin{align*}
    u_{1}(x, y, z, t) &= -a [e^{ax}\sin(ay+dz)+e^{az}\cos(ax+dy)]e^{-d^{2}t}, \\
    u_{2}(x, y, z, t) &= -a [e^{ay}\sin(az+dx)+e^{ax}\cos(ay+dz)]e^{-d^{2}t}, \\
    u_{3}(x, y, z, t) &= -a [e^{az}\sin(ax+dy)+e^{ay}\cos(az+dx)]e^{-d^{2}t}, \\
    p(x, y, z, t) &= -\frac{a^{2}}{2} [e^{2ax} + e^{2ay} + e^{2az} +  2 \sin(ax+dy)\cos(az+dx)e^{a(y+z)} \\
    &+ 2 \sin(ay+dz)\cos(ax+dy)]e^{a(z+x)} + 2 \sin(az+dx)\cos(ay+dz)]e^{a(x+y)}]e^{-2d^{2}t},
\end{align*}
where $a$ and $d$ are given. To show the generalizability of our methodology, we consider two cases of different coefficients. 
In the first case, we take $a = \frac{\pi}{2}$ and $d = \frac{\pi}{4}$. In the second case, we take $a = \pi$ and $d = \frac{\pi}{8}$. 
In the following section, we present the discovery and rediscovery results for IMEX-based meta-solvers in both cases. Full implementation details for training and testing, along with comprehensive numerical results, are provided in~\Cref{app_ns}.

\subsubsection{Discovery and re-discovery of optimal meta-solvers}

We take the same parameterization space as in~\Cref{subsec-rd} for the meta-solvers, that is
\begin{equation}
    \M = \NO \times \kr \times \rela \times S \times \{1,2,3\} \ .
\end{equation}

For the evaluation, we take an 8-dimensional performance function $f = (f_1,\dots, f_8)$, where each dimension represents a particular performance criterion. In particular, $f_1$ is the relative error for velocity $u$, $f_2$ is the relative error for pressure $p$, $f_3$ is the computational time, $f_4$ is the number of iterations, $f_5$ is the memory allocation, $f_6$ is the MACs, $f_7$ is the average MACs and $f_8$ is the training time, respectively.

\begin{figure}
    \centering
    \begin{minipage}[t]{0.45\textwidth}
    \centering
        \vspace{1.5cm}
        \subfigure[3-d projection of Pareto front and optimal meta-solvers discovered by preference functions.]{
            \centering
            \includegraphics[width=\linewidth]{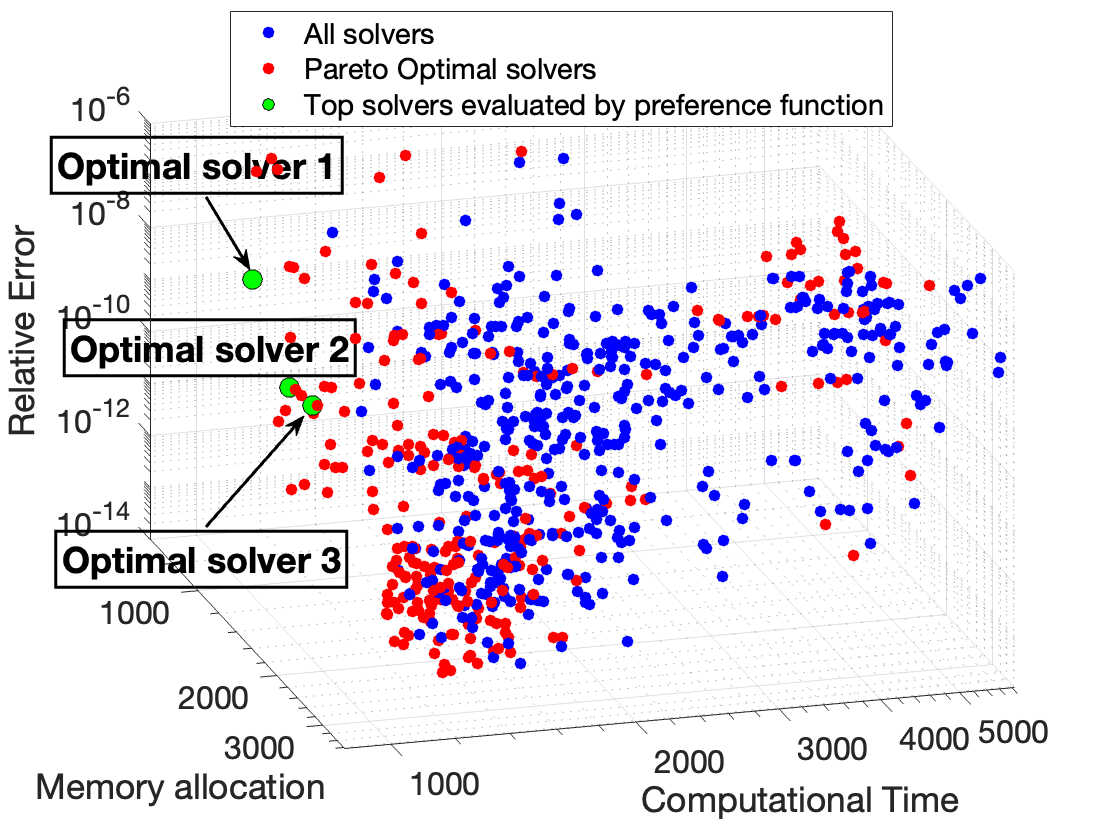}}
    \end{minipage}%
    \hfill
    \begin{minipage}[t]{0.55\textwidth}
    \begin{table}[H]
    \centering
     \begin{minipage}{\linewidth}
        \centering
        \captionof{subtable}{\small The top-3 meta-solvers.}
        \resizebox{\linewidth}{!}{
        \begin{tabular}{c|c|c|c|c|c}
\hline \hline 
 &Neural operator & Classical solver & Multi-grid & Relaxation & Strategies \\ 
\hline \hline
     Top-1 & KAN & CG & 3-level & SSOR & 5-1-5 \\ \hline
     Top-2 & DeepONet & CG & 3-level & Gauss-Seidel &  7-1-7 \\ \hline
     Top-3 & U-DeepONet & CG & 3-level & SSOR & 3-1-3 \\ \hline
\end{tabular}
}
    \end{minipage}

    \vspace{1em} 

    \begin{minipage}{\linewidth}
        \centering
        \captionof{subtable}{\small Performance of the top-3 solvers.}
        \resizebox{\linewidth}{!}{ 
          \begin{tabular}{c|c|c|c|c|c|c|c|c}
\hline \hline 
 &Relative error:u& Relative error:p& Com. time  & \# of ite & Memory & MACs&Ave. MACs & Training time \\ 
\hline \hline
Top-1 &$3.36 \times 10^{-11}$ & $2.28 \times 10^{-9}$ & 1095.7 & 6543 & 809.4 & $3.52 \times 10^{12}$ & $1.11 \times 10^{12}$ & 1086.8 \\
\hline
Top-2& $9.55 \times 10^{-14}$ & $1.34 \times 10^{-11}$ & 1229.4 & 7959 & 784.1 & $3.74 \times 10^{12}$ & $1.10 \times 10^{12}$ & 26.4 \\
\hline
Top-3& $1.60 \times 10^{-13}$ & $5.33 \times 10^{-12}$ & 1317.7 & 8505 & 777.6 & $4.05 \times 10^{12}$ & $1.07 \times 10^{12}$ & 49.6 \\
\hline
\end{tabular}
        }
    \end{minipage}
    \begin{minipage}{\linewidth}
        \centering
        \captionof{subtable}{\small Performance of vanilla method.}
        \scriptsize
          \begin{tabular}{c|c|c}
\hline \hline 
 & Com. time  & \# of ite  \\ 
 \hline \hline
 CG & 13824.8 & 324501 \\ \hline
 BiCGStab & 9946.89 & 220176 \\ \hline
  GMRES & -- & -- \\ \hline
\end{tabular}
        
    \end{minipage}
    \caption{Composition and performance of optimal meta-solvers}
    \label{table_ns_1}
\end{table}
    \end{minipage}
    \caption{\small IMEX based meta-solvers for solving the 3-d incompressible Navier-Stokes equations in case $a=\frac{\pi}{2}$ and $d = \frac{\pi}{4}$. Figures: (a) the projection of the 8-dimensional Pareto front into three-dimensional, time-error-memory, space, with Pareto optimal meta-solvers marked in red and optimal meta-solvers discovered by preference function highlighted in green. Tables: (A) the parameters of the top-3 meta-solvers discovered by preference function, (B) the 8-dimensional performance of the top-3 meta-solvers, (C) comparing the computational time and number of iteration for the vanilla IMEX method, using different Krylov methods as iterative solvers. Up to an error of machine precision, our meta-solvers obtain $\approx9.1$ times speedup in computational time and $>30$ times speedup in number of iterations.}
    \label{fig_ns_1}
\end{figure}

In~\Cref{fig_ns_1}, we present the results obtained using IMEX based meta-solvers for solving the case with $a=\frac{\pi}{2}$ and $d = \frac{\pi}{4}$. In particular,~\Cref{fig_ns_1} (a) shows a three-dimensional projection of the Pareto front. For this analysis, we define the preference function to be the average of all performance criteria.  
The performance of the top three meta-solvers is highlighted with green dots. 
In~\Cref{table_ns_1} (A), we present the composition of the top-3 meta-solvers discovered by the preference function $p_{ns}^1$, and their performance is shown in \Cref{table_ns_1} (B). 

\begin{figure}
    \centering
    \begin{minipage}[t]{0.45\textwidth}
    \centering
        \vspace{1.5cm}
        \subfigure[3-d projection of Pareto front and optimal meta-solvers discovered by preference functions.]{
            \centering
            \includegraphics[width=\linewidth]{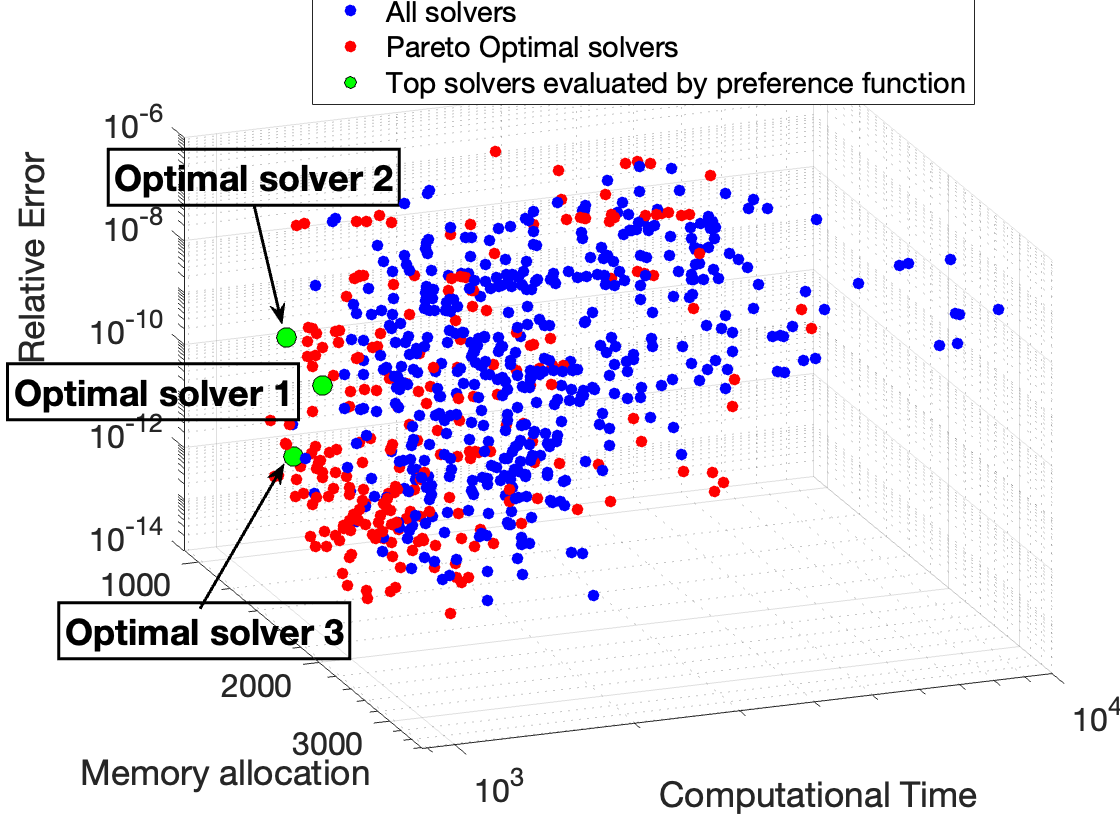}}
    \end{minipage}%
    \hfill
    \begin{minipage}[t]{0.55\textwidth}
    \begin{table}[H]
    \centering
     \begin{minipage}{\linewidth}
        \centering
        \captionof{subtable}{\small The top-3 meta-solvers.}
        \resizebox{\linewidth}{!}{
        \begin{tabular}{c|c|c|c|c|c}
\hline \hline 
 &Neural operator & Classical solver & Multi-grid & Relaxation & Strategies \\ 
\hline \hline
     Top-1 &  DeepONet & CG & 2-level & SSOR & 5-1-5 \\ \hline
     Top-2 &  JacobiKAN & CG & 2-level & Gauss-Seidel &  5-1-5 \\ \hline
     Top-3 &  DeepONet & CG & 2-level & Gauss-Seidel & 7-1-7 \\ \hline
\end{tabular}
}
    \end{minipage}

    \vspace{1em} 

    \begin{minipage}{\linewidth}
        \centering
        \captionof{subtable}{\small Performance of the top-3 solvers.}
        \resizebox{\linewidth}{!}{ 
          \begin{tabular}{c|c|c|c|c|c|c|c|c}
\hline \hline 
 &Relative error:u& Relative error:p& Com. time  & \# of ite & Memory & MACs&Ave. MACs & Training time \\ 
\hline \hline
Top-1 &$3.34 \times 10^{-11}$ & $2.60 \times 10^{-9}$ & 883.9 & 6249 & 2042.4 & $3.44 \times 10^{12}$ & $1.09 \times 10^{12}$ & 26.4 \\
\hline
Top-2&  $9.36 \times 10^{-11}$ & $1.72 \times 10^{-9}$ & 991.7 & 8519 & 1403.8 & $3.61 \times 10^{12}$ & $1.02 \times 10^{12}$ & 115.9 \\
\hline
Top-3&  $6.57 \times 10^{-14}$ & $8.21 \times 10^{-12}$ & 1016.1 & 7730 & 1408.3 & $3.66 \times 10^{12}$ & $1.08 \times 10^{12}$ & 26.4 \\
\hline
\end{tabular}
        }
    \end{minipage}
    \begin{minipage}{\linewidth}
        \centering
        \captionof{subtable}{\small Performance of vanilla method.}
        \scriptsize
          \begin{tabular}{c|c|c}
\hline \hline 
 & Com. time  & \# of ite  \\ 
 \hline \hline
 CG & 12736 & 326122 \\ \hline
 BiCGStab & 9832.6 & 218408 \\ \hline
  GMRES & -- & -- \\ \hline
\end{tabular}
        
    \end{minipage}
    \caption{Composition and performance of optimal meta-solvers}
    \label{table_ns}
\end{table}
    \end{minipage}
    \caption{\small IMEX based meta-solvers for solving the 3-d incompressible Navier-Stokes equations in case $a=\pi$ and $d = \frac{\pi}{8}$. Figures: (a) the projection of the 8-dimensional Pareto front into three-dimensional, time-error-memory, space, with Pareto optimal meta-solvers marked in red and optimal meta-solvers discovered by preference function highlighted in green. Tables: (A) the parameters of the top-3 meta-solvers discovered by preference function, (B) the 8-dimensional performance of the top-3 meta-solvers, (C) comparing the computational time and number of iteration for the vanilla IMEX method, using different Krylov methods as iterative solvers. Up to an error of machine precision, our meta-solvers obtain $\approx11$ times speedup in computational time and $>30$ times speedup in number of iterations.}
    \label{fig_ns}
\end{figure}

In~\Cref{fig_ns}, we present the results obtained using IMEX based meta-solvers for solving the case with $a=\pi$ and $d = \frac{\pi}{8}$. In particular,~\Cref{fig_ns} (a) shows a three-dimensional projection of the Pareto front. For this analysis, we define the preference function to emphasize the relative error for $p$ and $u$, as well as the computational time, 
that is after a rescaling process,  \[p_{ns}(f) = \frac{1}{25}   \Big( 5(f_1+f_2) + 10 f_3 +(\sum_{i=4}^bf_i) \Big) \ . \] 
The performance of the top three meta-solvers is highlighted with green dots. 
In~\Cref{table_ns} (A), we present the composition of the top-3 meta-solvers discovered by the preference function $p_{ns}^1$, and their performance is shown in \Cref{table_ns} (B).

We also present below various weighted-sum preference functions along with the corresponding optimal meta-solvers identified by each, for the case with $a = \pi$ and $d = \frac{\pi}{8}$. 
The weights are presented in the following order: $\lambda_1$ for relative error of $u$, $\lambda_2$ for relative error of $p$, $\lambda_3$ for computational time, $\lambda_4$ for number of iterations, $\lambda_5$ for memory allocation, $\lambda_6$ for MACs, $\lambda_7$ for average MACs, and $\lambda_8$ for training time. 

\begin{subequations}
\begin{equation}
\left\{
\begin{aligned}
&\text{Preference function $p^a_{ns}(f) = (0,0,0,0.551,0.262,0,0.187,0)^Tf$  } \\
&\text{Optimal solver}: x^a = (\text{JacobiKAN}, \text{CG}, \text{3-level}, \text{SSOR}, 5-1-5) \ .
\end{aligned}
\right.
\end{equation}

\begin{equation}
\left\{
\begin{aligned}
&\text{Preference function $p^b_{ns}(f) = (0.512,0,0,0.023,0.007,0.105,0.354,0.001)^Tf$ } \\
&\text{Optimal solver}: x^a = (\text{JacobiKAN}, \text{CG}, \text{2-level}, \text{Gauss-Seidel}, 3-1-3) \ .
\end{aligned}
\right.
\end{equation}

\begin{equation}\label{ns_rediscovery}
\left\{
\begin{aligned}
&\text{Preference function $p^c_{ns}(f) = (0,0.516,0,0.435,0.001,0,0.048, 0)^Tf$ } \\
&\text{Optimal solver}: x^a = (\text{JacobiKAN}, \text{CG}, \text{2-level}, \text{SSOR}, 7-1-7) \ .
\end{aligned}
\right.
\end{equation}

\end{subequations}

Take, for example, the preference function and its corresponding optimal meta-solver discovered in~\eqref{ns_rediscovery}. 
If a user places a strong emphasis on relative error of $p$ ($\approx52\%$ weight), and number of iterations ($\approx44\%$ weight), then the optimal meta-solver for this preference profile consists of a Neural Operator JacobiKAN combined with the Krylov method CG, using SSOR as the smoother. The smoother is applied with the 7-1-7 strategy and leverages a two-level multigrid technique.

\subsection{Solving brittle fracture problems}\label{subsec-bf}
\subsubsection{The model and implementation details}
In this section, we consider the brittle fracture problem. In particular, 
assume a domain $\Omega \subset \R^{2}$ is given. 
The external boundary $\Gamma$ is decomposed into the Dirichlet boundary $\Gamma_{D}$ and Neumann boundary $\Gamma_{N}$, i.e., $\Gamma = \Gamma_{D} \cup \Gamma_{N}$.
Moreover, the Dirichlet boundary $\Gamma_{D}$ consists of a homogeneous boundary $\Gamma_{D,0}$ and non-homogeneous~(loading) boundary $\Gamma_{D,1}$. 
Let $\Omega_{c} \subset \Omega$ be a crack set.
Then, we consider the $n$-th loading step of the phase-field modeling of brittle fracture.
The pair of solution spaces $\{V,Q\}$ is given by
\begin{align}
    V &= \{\mathbf{u} \in [H^{1}(\Omega)]^{d}\colon \mathbf{u}=0 \text{ on } \Gamma_{D,0}, \mathbf{u}=\mathbf{u}_{n} \text{ on } \Gamma_{D,1}\}, \\
    Q &= \{\alpha \in H^{1}(\Omega)\colon 0 \leq \alpha \leq 1, \alpha \geq \alpha_{n-1} \text{ in } \Omega\}. 
\end{align}
Then, the pair of solutions $\{\mathbf{u}, \alpha\}$ is obtained by solving the following minimization problem:
\begin{subequations}
\begin{equation}
\min_{\mathbf{u} \in V, \alpha \in Q} \mathcal{E}_{n}(\mathbf{u}, \alpha) := \mathcal{E}_{el} + \mathcal{E}_{d} - \mathcal{E}_{w} + \mathcal{E}_{irr},
\label{eqn:phase-field}
\end{equation}
where
\begin{equation}
\left\{\begin{split}
\mathcal{E}_{el}&=\int_{\Omega} \left(a(\alpha)\Psi^{+}(\varepsilon(\mathbf{u}))+\Psi^{-}(\varepsilon(\mathbf{u})\right)dx, \\
\mathcal{E}_{d}&=\frac{G_{c}}{c_{w}}\int_{\Omega} (\frac{w(\alpha)}{\ell}+\ell \vert \nabla \alpha \vert^{2})dx, \\
\mathcal{E}_{w} &= \int_{\Omega} \mathbf{f}_{n} \cdot \mathbf{u} dx + \int_{\Gamma_{N}}\mathbf{t}_{n} \cdot \mathbf{u} ds, \\
\mathcal{E}_{irr} &= \frac{\gamma}{2}\int_{\Omega} \langle \alpha - \alpha_{n-1} \rangle_{-}^{2} dx.
\end{split}
\right.
\label{eqn:phase-field-def}
\end{equation}
\end{subequations}
This minimization problem~\eqref{eqn:phase-field} is addressed via the Karush-Kuhn-Tucker (KKT) condition, employing an alternating minimization strategy. The resulting local subproblems are solved using Newton-Raphson based meta-solvers. 

In this implementation, 
we enrich the parameterization space, in particular the set of relaxation methods, by incorporating the Chebyshev semi-iterative method. Additionally, we replace the geometric multigrid method with the algebraic multigrid (AMG) approach. 
For the evaluation, we consider a similar 8-dimensional performance function as in~\Cref{subsec_ns}, with the distinction that $f_2$ now representing the relative error for $\alpha$. 
The implementation details, a brief description of Chebyshev semi-iterative method, together with full numerical results are presented in~\Cref{app-bf}.

\subsubsection{Discovery and re-discovery of optimal meta-solvers}

\begin{figure}
    \centering
    \begin{minipage}[t]{0.45\textwidth}
    \centering
        \vspace{1.5cm}
        \subfigure[3-d projection of Pareto front and optimal meta-solvers discovered by preference functions.]{
            \centering
            \includegraphics[width=\linewidth]{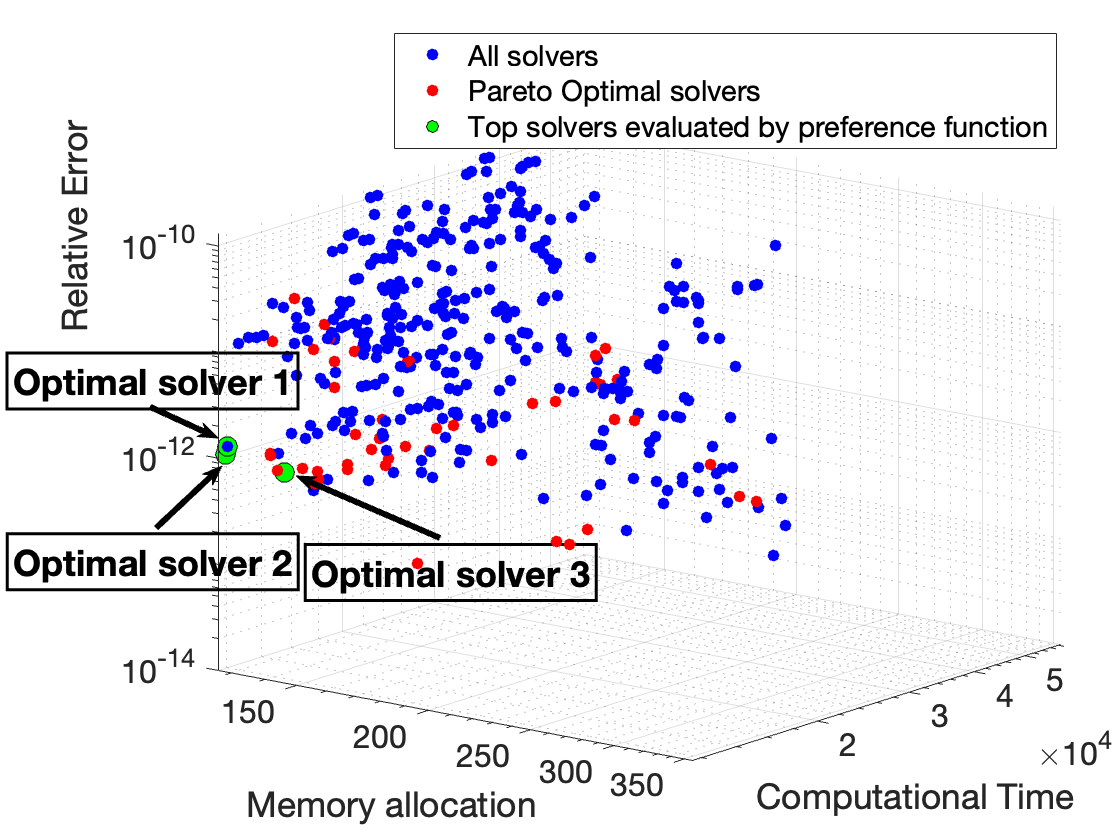}}
    \end{minipage}%
    \hfill
    \begin{minipage}[t]{0.55\textwidth}
        \begin{table}[H]
    \centering
     \begin{minipage}{\linewidth}
        \centering
        \captionof{subtable}{\small The top-3 meta-solvers.}
        \resizebox{\linewidth}{!}{
        \begin{tabular}{c|c|c|c|c|c}
\hline \hline 
 &Neural operator & Classical solver & Relaxation & Strategies & AMG \\ 
\hline \hline
     Top-1 &  U-DeepONet & CG  & Chebyshev & 1-1-1& YES \\ \hline
     Top-2 & DeepONet & CG  & Chebyshev &  1-1-1& YES \\ \hline
     Top-3 & U-DeepONet & BiCGStab  & Chebyshev & 1-1-1& YES \\ \hline
\end{tabular}
}
    \end{minipage}

    \vspace{1em} 

    \begin{minipage}{\linewidth}
        \centering
        \captionof{subtable}{\small Performance of the top-3 solvers.}
        \resizebox{\linewidth}{!}{ 
          \begin{tabular}{c|c|c|c|c|c|c|c|c}
\hline \hline 
 &relative error:u& relative error:$\alpha$& Com. time  & \# of ite & Memory & MACs&Ave. MACs & Training time \\ 
\hline \hline
Top-1& $4.10 \times 10^{-6}$  & $1.23 \times 10^{-12}$ & 11937.4 & 521754 & 126.9 & $2.63 \times 10^{13}$ & $8.96 \times 10^{12}$ & 2110.97 \\
\hline
Top-2&  $4.10 \times 10^{-6}$  & $1.04 \times 10^{-12}$ & 11844.2 & 521426 & 126.7 & $2.63 \times 10^{13}$ & $8.95 \times 10^{12}$ & 2388.54 \\
Top-3&  $4.10 \times 10^{-6}$  & $4.77 \times 10^{-13}$ & 15300.9 & 363977 & 127.1 & $3.55 \times 10^{13}$ & $1.21 \times 10^{13}$ & 2110.97 \\
\hline
\end{tabular}
        }
    \end{minipage}
    \begin{minipage}{\linewidth}
        \centering
        \captionof{subtable}{\small Performance of vanilla method.}
        \scriptsize
          \begin{tabular}{c|c|c}
\hline \hline 
 & Com. time  & \# of ite  \\ 
 \hline \hline
 CG & 43000 & 57771176 \\ \hline
\end{tabular}
        
    \end{minipage}
    \caption{Composition and performance of optimal meta-solvers}
    \label{table_bf}
\end{table}
    \end{minipage}
    \caption{\small Newton-Raphson based meta-solvers for solving the brittle fracture problem. Figures: (a) A sketch of the generated mesh for computation, (b) projection of the eight-dimensional Pareto front into three-dimensional, time-error-memory, space, with Pareto optimal meta-solvers marked in red and optimal meta-solvers discovered by preference function highlighted in green. Tables: (A) the parameters of the top-three meta-solvers discovered by preference function, (B) the eight-dimensional performance of the top-three meta-solvers, (C) comparing the computational time and number of iteration for the vanilla Newton-Raphson method, using CG as iterative solvers. Up to an error of machine precision, our meta-solvers obtain $\approx3.6$ times speedup in computational time and $>150$ times speedup in number of iterations.}
\label{fig_bf}
\end{figure}

In~\Cref{fig_bf}, we present the results obtained using Newton-Raphson based meta-solvers. In particular,~\Cref{fig_bf} (a) provides a schematic illustration of the generated mesh for solving this problem. 
\Cref{fig_bf} (b) shows a three-dimensional projection of the Pareto front. For this analysis, we define the preference function to be the average of the performance. In particular, the performance data in each dimension is first rescaled to a real value in$[0,1]$, then, 
\[p_{bf}(f) = \frac{1}{8}  \sum_{i=1}^8 f_i \ . \] 
The performance of the top three meta-solvers is highlighted with green dots. 
In~\Cref{table_bf} (A), we present the composition of the top-3 meta-solvers discovered by the preference function $p_{ns}^1$, and their performance is shown in \Cref{table_bf} (B). 

For other preference functions and their corresponding optimal meta-solvers, we present the results below. The weights are ordered as follows: $\lambda_1$ for relative error of $u$, $\lambda_2$ for relative error of $a$, $\lambda_3$ for computational time, $\lambda_4$ for number of iterations, $\lambda_5$ for memory allocation, $\lambda_6$ for MACs, $\lambda_7$ for average MACs, and $\lambda_8$ for training time. 

\begin{subequations}
\begin{equation}
\left\{
\begin{aligned}
&\text{Preference function $p^a_{bf}(f) = (0,0,0.03,0.73,0,0.25,0,0)^Tf$  } \\
&\text{Optimal solver}: x^a = (\text{DeepONet}, \text{CG}, \text{with AMG}, \text{Chebyshev}, 1-1-1) \ .
\end{aligned}
\right.
\end{equation}

\begin{equation}
\left\{
\begin{aligned}
&\text{Preference function $p^b_{bf}(f) = (0,0.79,0,0,0,0,0.02,0.19)^Tf$ } \\
&\text{Optimal solver}: x^a = (\text{U-DeepONet}, \text{CG}, \text{no AMG}, \text{Chebyshev}, 7-1-7) \ .
\end{aligned}
\right.
\end{equation}

\begin{equation}\label{bf_rediscovery}
\left\{
\begin{aligned}
&\text{Preference function $p^c_{bf}(f) = (0,0.22,0.02,0.54,0,0,0.23, 0)^Tf$ } \\
&\text{Optimal solver}: x^a = (\text{KAN}, \text{BiCGStab}, \text{with AMG}, \text{Chebyshev}, 1-1-1) \ .
\end{aligned}
\right.
\end{equation}

\end{subequations}

Take, for example, the preference function and its corresponding optimal meta-solver discovered in~\eqref{bf_rediscovery}. 
If a user places  emphasis on relative error of $a$ ($\approx22\%$ weight), number of iterations ($\approx54\%$ weight), and average MACs ($\approx23\%$ weight), then the optimal meta-solver for this preference profile consists of a Neural Operator KAN combined with the Krylov method BiCGStab, using Chebyshev semi-iterative method as the smoother. The smoother is applied with the 1-1-1 strategy and leverages the algebraic multigrid technique.

\section{Summary}\label{sec-sum}

\subsection{Summary of contributions}
In the present work, we introduce two families of meta-solvers for solving time-dependent nonlinear equations. For time evolution, we adopt a Crank–Nicolson framework and apply the meta-solvers in the spatial domain while advancing in time. To handle the nonlinear terms, we enhance Krylov-based meta-solvers, originally designed for linear systems, by hybridizing them with the Newton–Raphson method and the implicit-explicit (IMEX) scheme. This leads to the development of two novel classes of meta-solvers tailored for time-dependent nonlinear problems.

After parameterizing both the meta-solvers and the performance metrics, we implement a Pareto optimality-based evaluation methodology. 
Building on this, we apply a preference function based discovery approach to identify optimal solvers, as well as a linear programming based rediscovery methodology to retrieve specific known meta-solvers from among the Pareto-optimal set. 

We apply these two classes of nonlinear meta-solvers to benchmark problems in reacting transport (reaction–diffusion equations), fluid mechanics (Navier–Stokes equations), and solid mechanics (brittle fracture). In all cases, the discovered optimal meta-solvers demonstrate substantial performance gains. Below, we summarize the speedups in computational time and iteration count achieved by the optimal meta-solvers compared to the best classical iterative solvers.
\begin{itemize}
\item Reaction–diffusion equations (two test cases):
IMEX-based meta-solvers achieve approximately $15$ times speedup in computational time and $> 76  $ times in number of iterations. Newton–Raphson based meta-solvers achieve approximately $11.7 $ times speedup in computational time and $>100 $ times in number of iterations.
\item Navier–Stokes equations (two test cases): IMEX-based meta-solvers achieve approximately $9.1 $ times  speedup in computational time and $>30$ times  in number of iterations for the first case, and approximately $11$ times and $>30$ times, respectively, for the second case. 
\item Brittle fracture problem: The Newton-Rapson based meta-solvers achieves approximately $3.6$ times  speedup in computational time and $>150$ times in number of iterations. 
 \end{itemize}

\subsection{Discussion}

The existing results already demonstrate strong effectiveness and significant improvements that generalize well across diverse application domains. Nevertheless, they also inspire further possibilities for advancement and development. 

First all of, the discovered meta-solvers demonstrate up to 10 times improvement in computational time and 100 times reduction in the number of iterations.
This implies that the proposed meta-solvers significantly reduce the condition number of linearized problem to which Krylov methods are applied, leading to a dramatic decrease in iteration count.
Although each iteration of the meta-solvers incurs a higher computational cost compared to vanilla Krylov methods, the overall runtime can still improve by up to 10 times.
Furthermore, Newton-Rapson based meta-solvers have to assemble both the linear system and an appropriate preconditioner at every Newton iteration, which imposes an additional computational overhead.
As a future research direction, it would be interesting to develop a purely nonlinear meta-solver that can eliminate such costs.

Another important observation is that the Pareto front forms the boundary of a nonconvex set. 
This inherent nonconvexity poses significant challenges from an optimization perspective. 
In particular, we find that up to 50\% of Pareto-optimal meta-solvers cannot be discovered using our current linear programming approach. 
This limitation arises because, geometrically, the weighted-sum preference function can only capture solutions on the convex hull of the Pareto front. 
However, to the best of our knowledge, no existing optimization techniques can fully explore nonconvex Pareto fronts in such high-dimensional spaces. 
As part of future work, we aim to investigate the topological structure of Pareto fronts in different problem domains and explore advanced tools from optimization and optimal control theory to address the nonconvex regions of the Pareto front.

\newcommand{\etalchar}[1]{$^{#1}$}

\newpage
\appendix

\section{Implementation details and further numerical results for solving time-dependent reaction-diffusion equations}\label{app_rd}

Let us recall that the time-dependent reaction-diffusion equation we consider has the following form: 

\begin{equation}
\left\{\begin{split}
\frac{\partial u}{\partial t} &= \nabla \cdot (k \nabla u) + R(u) + f, \,\,\, \text{ in } \Omega \times (0,1], \\
u &= u_{0}, \text{ in } t=0, \\
u &= 0, \text{ in } \partial\Omega \times (0,1],
\end{split}
\right.
\label{app:eqn:time-reaction-diffusion}
\end{equation}
where $R(u) = u - u^{2}$ is a Fisher-type reaction term, the spatial domain is $\Omega = [0, 1]^{2} \subset \R^{2}$, and the time interval is $T=[0,1]$. The coefficient $k(\mathbf{x})$ and force term $ f(\mathbf{x})$ are given. In the first implementation, we take the diffusion coefficient to satisfy the following:
\begin{equation}
k \sim N(1.0, K(\mathbf{x},\mathbf{x}')=0.3e^{-\frac{\Vert \mathbf{x}-\mathbf{x}'\Vert^{2}}{2 \times 0.1^{2}}}), \text{ and } k \geq 0.1 \ .
\end{equation}
For the external force $f$, we take
\begin{equation}
    f \sim N(0.0, K(\mathbf{x},\mathbf{x}')=e^{-\frac{\Vert \mathbf{x}-\mathbf{x}'\Vert^{2}}{2 \times 0.1^{2}}}).
\end{equation}

To discretize the spatial domain $\Omega$, we perform the triangulation and utilize the linear finite element with mesh size $h=1/30$.
For the time discretization scheme, we utilize the backward Euler method, which can allow us to use the time step $\Delta t=h=1/30$.
Moreover, since the neural operator only acts on the steady-state problem induced by time discretization scheme, it is enough to generate the training samples consisting of the steady-state problem.
To train neural operators, we take $30,000$ samples of $k$, $f$, and the reference solution $u$.
The reference solution is simply obtained by using direct solver~\cite{MUMPS:1,MUMPS:2}.
All neural operators are trained using AdamW~\cite{loshchilov2018decoupled} with batch size $1,000$ until the relative $L^{2}$ error is lower than $8\%$ or the number of iteration reaches $20,000$.

\subsection{Identifying all optimal meta-solvers in Pareto sense} 
For solving equation~\eqref{app:eqn:time-reaction-diffusion}, 
we implement both Newton–Raphson-based and IMEX-based meta-solvers. Our parameterization results in a total of 900 meta-solvers for each approach. Among these, 138 are identified as Pareto optimal for the Newton–Raphson method, and 161 for the IMEX method. 
In the following, we first summarize the composition of the Pareto optimal sets by counting the occurrences of different components used in the construction of the meta-solvers.

\begin{table}[H]
    \centering
    \caption{The composition of the set of Pareto optimal solvers by counting the number of elements in each dimension, for solving time-dependent 2-d reaction-diffusion equation, using Newton-Raphson based and IMEX based meta-solvers.}
    \begin{minipage}{\linewidth}
        \centering
        \captionof{subtable}{Different neural operators.}
\resizebox{\linewidth}{!}{
        \begin{tabular}{c| c c c  c c|c } 
\hline
Neural Op & DeepONet & U-DeepONet  & KAN & JacobiKAN& ChebyKAN&Total \\
\hline
\# in Pareto opt for Newton  & 26 & 46 & 15 & 31 & 20 & 138\\
\hline
\# in Pareto opt for IMEX  & 29 & 41 & 36 & 24 & 31 & 161\\
\hline
\end{tabular}}

    \end{minipage}

    \vspace{1em} 

    \begin{minipage}{\linewidth}
        \centering
        \captionof{subtable}{Different Krylov solvers.}
           \begin{tabular}{c| c c c c } 
\hline
Classical solvers & FGMRES & CG & BiCGStab \\
\hline
\# in Pareto opt for Newton  & 47 & 23 & 68  \\
\hline
\# in Pareto opt for IMEX  & 52 & 25 & 84  \\
\hline
\end{tabular}
        
    \end{minipage}
    
    \vspace{1em}

    \begin{minipage}{\linewidth}
        \centering
        \captionof{subtable}{Different smoothers.}
           \begin{tabular}{c| c c c c    } 
\hline
Smoother & GS & Jacobi & SOR & SSOR  \\
\hline
\# in Pareto opt for Newton & 23 & 18 & 38 & 59  \\
\hline
\# in Pareto opt for IMEX & 32 & 13 & 44 & 72  \\
\hline
\end{tabular}
        
    \end{minipage}

\vspace{1em}

\begin{minipage}{\linewidth}
        \centering
        \captionof{subtable}{Different strategies of applying smoothers.}
          \begin{tabular}{c| c c c c c } 
\hline
Strategies for smoother & 1-1-1 & 3-1-3 & 5-1-5 & 7-1-7 & 9-1-9  \\
\hline
\# in Pareto opt for Newton & 14 & 39 & 31 & 26 & 28 \\
\hline
\# in Pareto opt for IMEX & 16 & 42 & 34 & 30 & 39 \\
\hline
\end{tabular}
        
    \end{minipage}
    
\vspace{1em}

\begin{minipage}{\linewidth}
        \centering
        \captionof{subtable}{Different levels of multigrid.}
    \begin{tabular}{c| c c c } 
\hline
Levels in multi-grid & 1-level & 2-level &  3-level   \\
\hline
\# in Pareto opt for Newton & 1 & 76 & 61  \\
\hline
\# in Pareto opt for IMEX & 16 & 94 & 51  \\
\hline
\end{tabular}
       \end{minipage}
    \label{rd_com1}
\end{table}

Moreover, we plot the Pareto fronts for both Newton-Paphson method, in~\Cref{front_2drd_new}, and IMEX method, in~\Cref{front_2drd_im}, for the computational time, relative error and number of iteration, for the relative error, MACs and Memory allocation, and for the  Relative error, convergence rate and MACs.  

\begin{figure}[H]
\centering
\subfigure[Computational time -- Relative error -- \# of iterations]{
\includegraphics[width=0.45\textwidth]{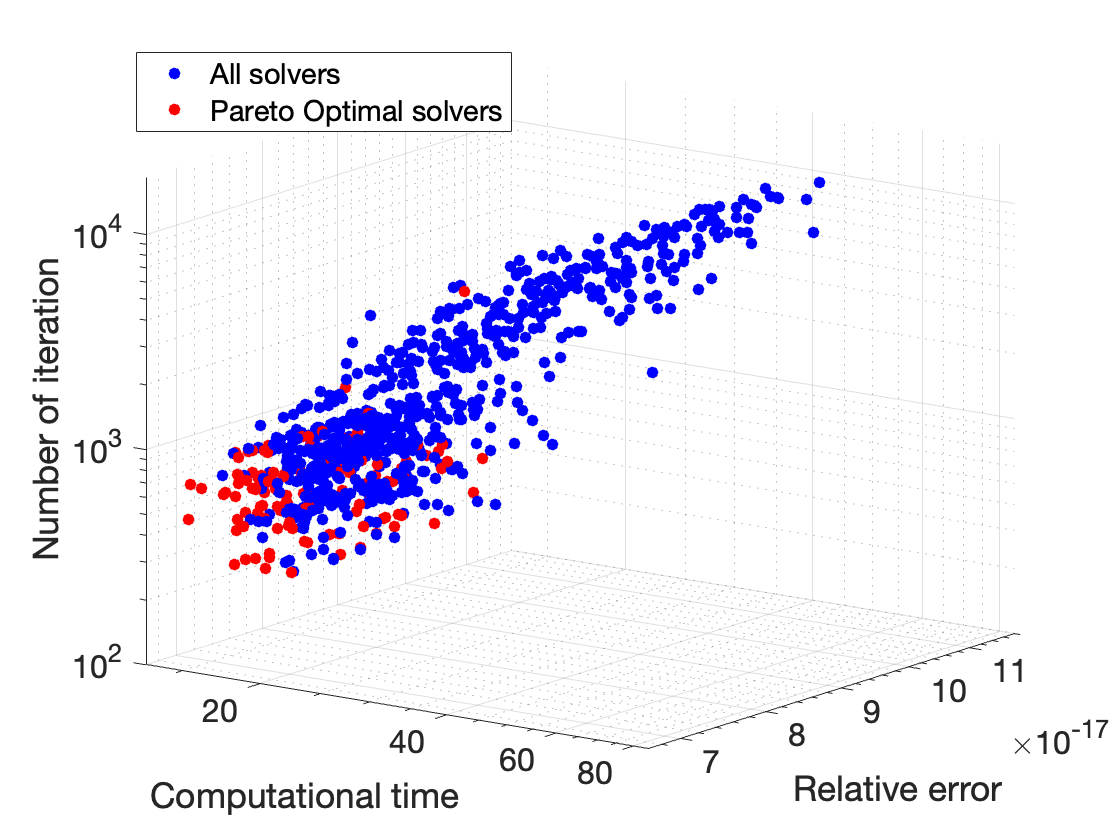}}
\subfigure[Relative error -- MACs -- Memory allocation]{
\includegraphics[width=0.45\textwidth]{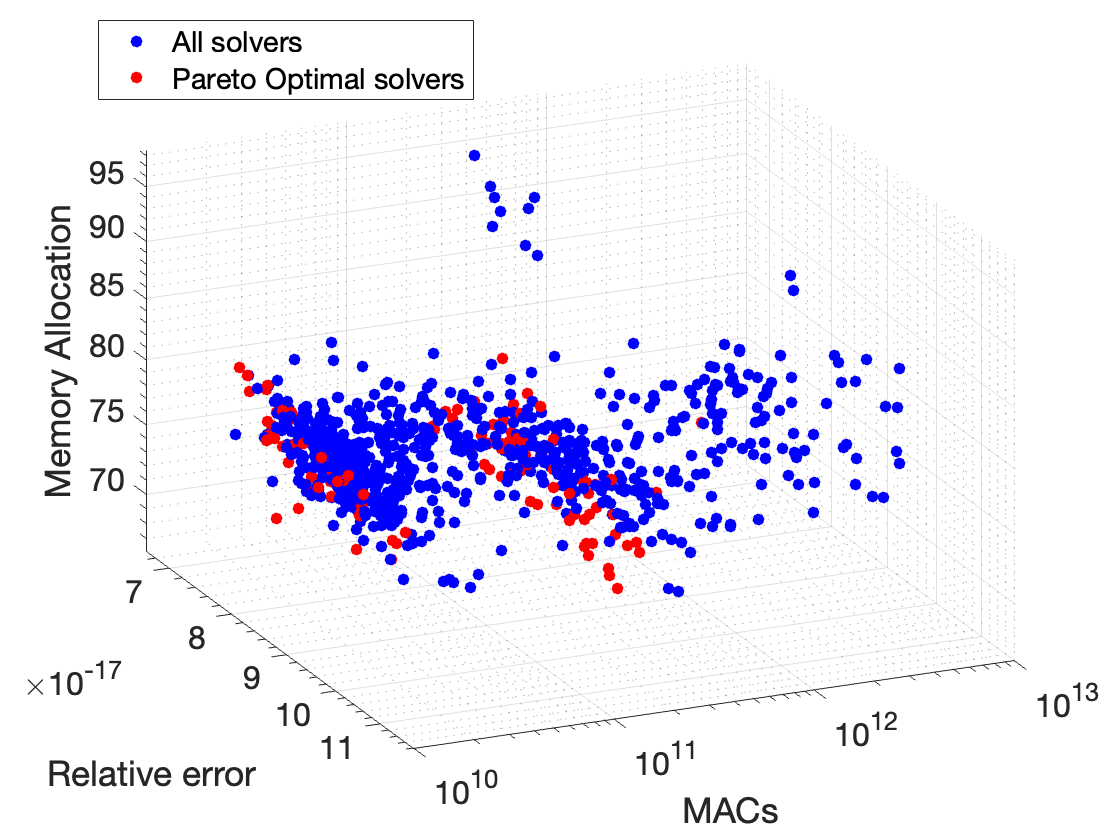}}
\subfigure[Relative error -- MACs -- \# of iterations]{
\includegraphics[width=0.45\textwidth]{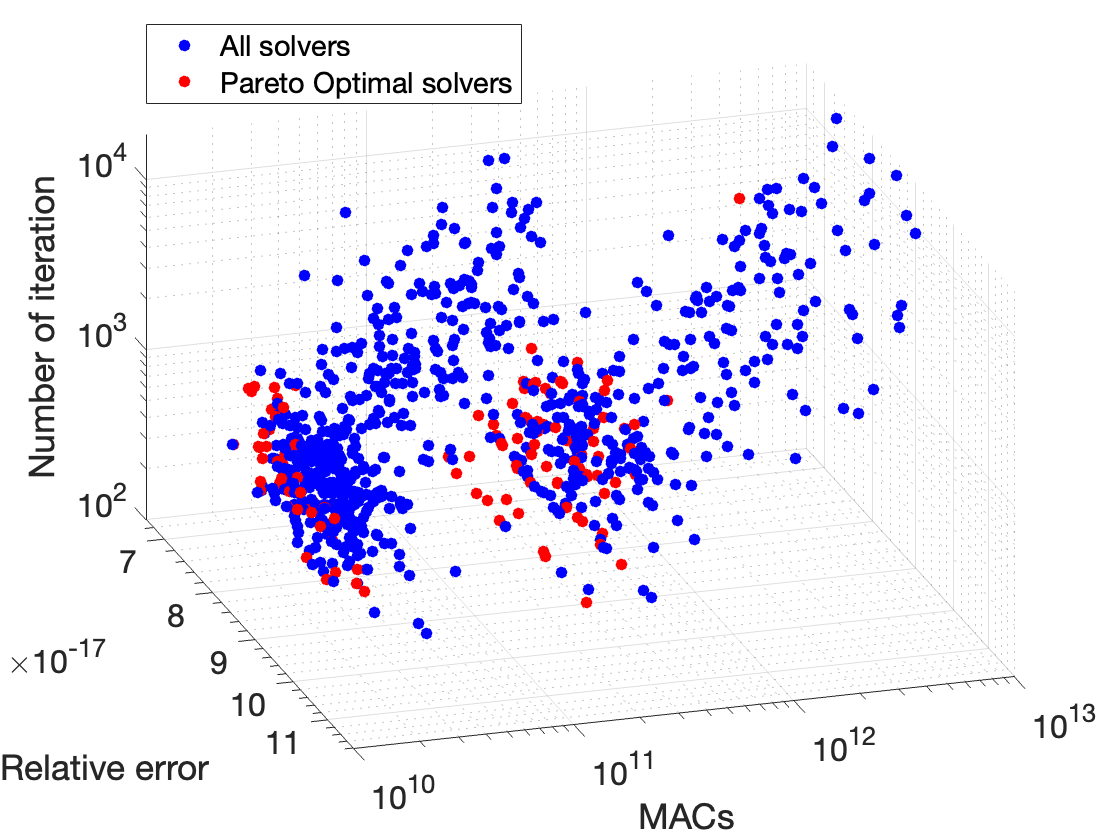}}
\caption{Pareto fronts for solving 2-D time-dependent reaction-diffusion equation using Newton-Raphson method. All solvers are depicted in blue while Pareto optimal solvers are highlighted in red. The``gap'' due to the adaption thus the improvement of performance by multi-grid techniques.}
\label{front_2drd_new}
\end{figure}

\begin{figure}[H]
\centering
\subfigure[Computational time -- Relative error -- \# of iterations]{
\includegraphics[width=0.45\textwidth]{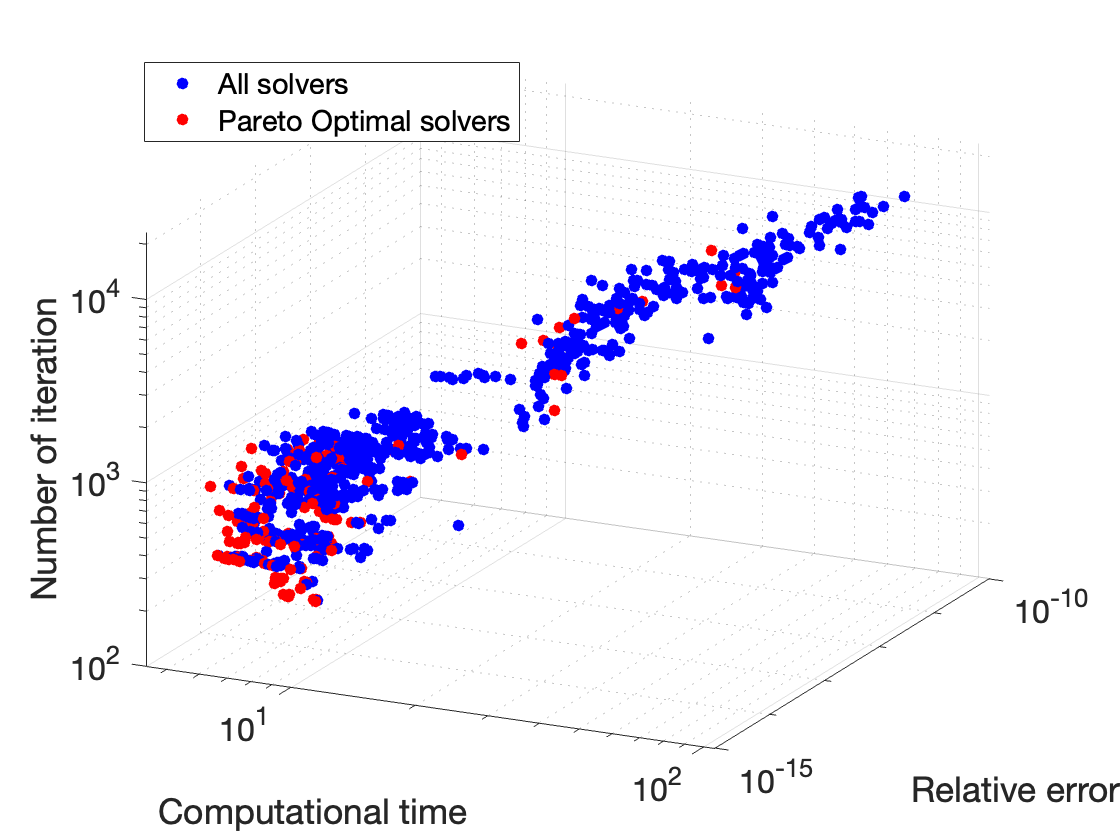}}
\subfigure[Relative error -- MACs -- Memory allocation]{
\includegraphics[width=0.45\textwidth]{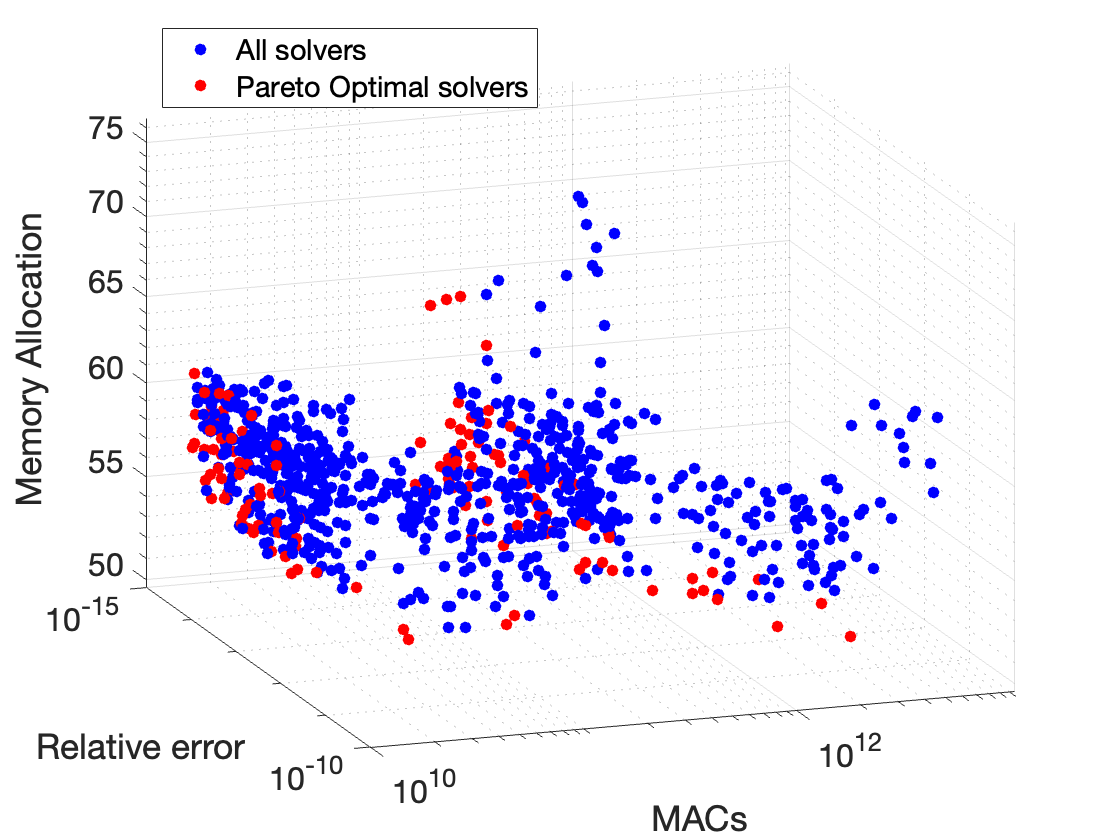}}
\subfigure[Relative error -- MACs -- \# of iterations]{
\includegraphics[width=0.45\textwidth]{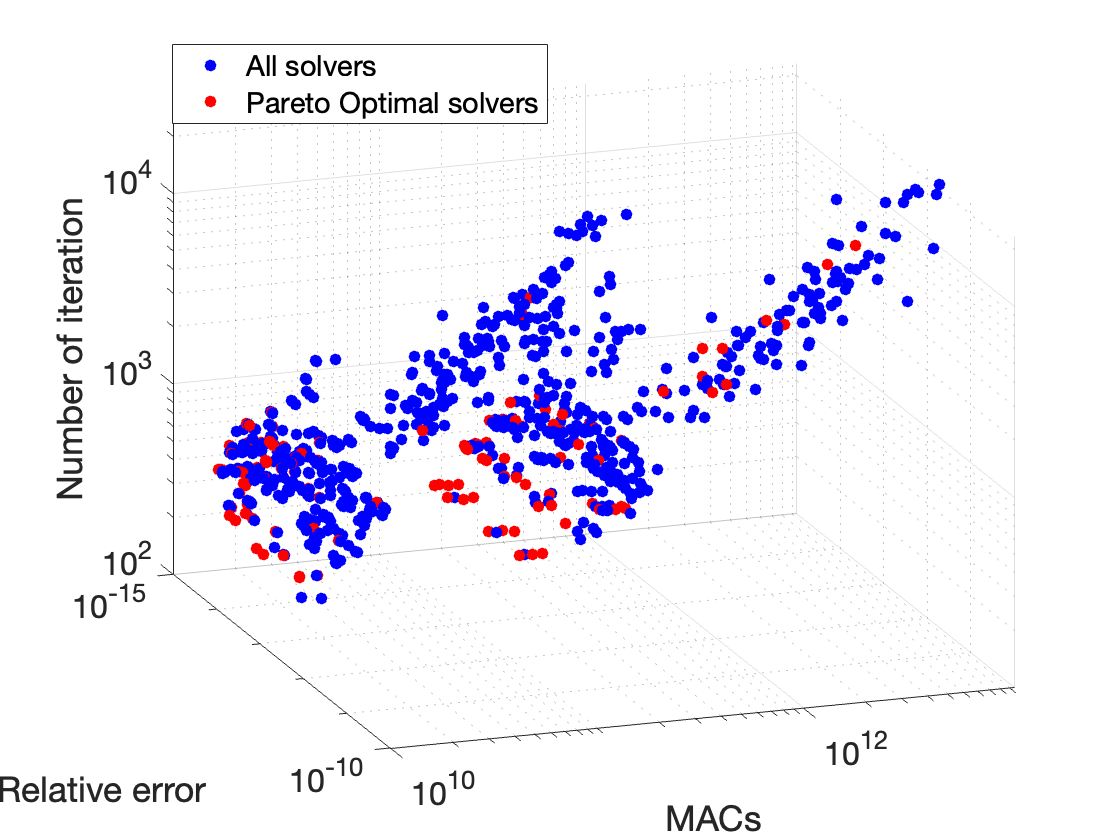}}
\caption{Pareto fronts for solving 2-D time-dependent reaction-diffusion equation using IMEX method. All solvers are depicted in blue while Pareto optimal solvers are highlighted in red. The``gap'' due to the adaption thus the improvement of performance by multi-grid techniques.}
\label{front_2drd_im}
\end{figure}

\subsection{Discovery and re-discovery of optimal meta-solvers by preference functions}
In this section, we present the results of applying the preference function based methodology to discover the optimal meta-solvers, in various of context. 
\begin{pref}
$p^1(f') = \frac{1}{6}(\sum_{i=1}^6 f'_i)$. The average of all the re-scaling performance. The top 3 solvers and their performance using this preference function for Newton-Raphson method is shown in~\Cref{top3_rd_p1new}. 
The top 3 solvers and their performance using this preference function for IMEX method is shown in~\Cref{top3_rd_p1im}.
\end{pref}

\begin{table}[H]
    \centering
    \caption{Top-3 solvers evaluated by preference function $p^1$ for Newton-Raphson method, for solving 2-d time-dependent reaction-diffusion equation.}
     \begin{minipage}{\linewidth}
        \centering
        \captionof{subtable}{The top-3 meta-solvers}
        \resizebox{\linewidth}{!}{
        \begin{tabular}{c|c|c|c|c|c}
\hline \hline 
 &Neural operator & Classical solver & Multi-grid & Relaxation & Strategies \\ 
\hline \hline
    Top 1 solver & U-Net & CG & 3-level & Jacobi & 3-1-3 \\ \hline
Top 2 solver & U-Net & CG & 2-level & SOR &  7-1-7 \\ \hline
Top 3 solver  & U-Net & FGMRES & 2-level & SSOR & 5-1-5 \\ \hline
\end{tabular}
}
    \end{minipage}

    \vspace{1em} 

    \begin{minipage}{\linewidth}
        \centering
        \captionof{subtable}{Performance and rank of performance of the top-3 solvers}
        \resizebox{\linewidth}{!}{ 
          \begin{tabular}{c|c|c|c|c|c|c}
\hline \hline 
 &Error & Com. time  & \# of ite & Memory & MACs & Training time \\ 
\hline \hline
Top 1 solver & $7.87 \times 10^{-17}$ & 16.1725 & 896 & 66.1221 & $6.63 \times 10^{11}$ & 9010.38 \\ \hline
Top 2 solver &  $7.55 \times 10^{-17}$ & 16.2697 & 640 & 69.126  & $4.27 \times 10^{11}$ & 9010.38 \\ \hline
Top 3 solver &  $7.72 \times 10^{-17}$ & 14.9918 & 512 & 72.0439 & $3.42 \times 10^{11}$ & 9010.38 \\ \hline
\end{tabular}
        }
    \end{minipage}
    \label{top3_rd_p1new}
\end{table}

\begin{table}[H]
    \centering
    \caption{Top-3 solvers evaluated by preference function $p^1$ for IMEX method, for solving 2-d time-dependent reaction-diffusion equation.}
     \begin{minipage}{\linewidth}
        \centering
        \captionof{subtable}{The top-3 solvers}
        \resizebox{\linewidth}{!}{
        \begin{tabular}{c|c|c|c|c|c}
\hline \hline 
 &Neural operator & Classical solver & Multi-grid & Relaxation & Strategies \\ 
\hline \hline
Top 1 solver & U-Net & CG & 2-level & SOR & 7-1-7 \\ \hline
Top 2 solver & U-Net & CG & 2-level & SOR &  5-1-5 \\ \hline
Top 3 solver & U-Net & CG & 2-level & Gauss-Seidel & 9-1-9 \\ \hline
\end{tabular}
}
    \end{minipage}

    \vspace{1em} 

    \begin{minipage}{\linewidth}
        \centering
        \captionof{subtable}{Performance and rank of performance of the top-3 solvers}
        \resizebox{\linewidth}{!}{ 
          \begin{tabular}{c|c|c|c|c|c|c}
\hline \hline 
 &Error & Com. time  & \# of ite & Memory & MACs & Training time \\ 
\hline \hline
Top 1 solver  & $1.60 \times 10^{-14}$ & 7.66828 & 640 & 52.6182 & $4.21 \times 10^{11}$ & 10973 \\ \hline
Top 2 solver  &  $1.23 \times 10^{-13}$ & 6.85867 & 640 & 52.582  & $4.19 \times 10^{11}$ & 10973 \\ \hline
Top 3 solver &  $1.11 \times 10^{-14}$ & 8.47799 & 640 & 53.2295 & $4.22 \times 10^{11}$ & 10973 \\ \hline
\end{tabular}
        }
    \end{minipage}
    \label{top3_rd_p1im}
\end{table}

\begin{pref}
 $p^2(f') =  \frac{1}{24} \big( 10 (f'_1 + f'_2) +  (f'_3 + f'_4 + f'_5 + f'_6 ) \big)$. This preference function means that the computational time and relative error are ten times more important then other criteria. 
 The top 3 solvers and their performance using this preference function for Newton-Raphson method is shown in~\Cref{top3_rd_p2new}. 
The top 3 solvers and their performance using this preference function for IMEX method is shown in~\Cref{top3_rd_p2im}.
\end{pref}

\begin{table}[H]
    \centering
    \caption{Top-3 solvers evaluated by preference function $p^2$ for Newton-Raphson method, for solving 2-d time-dependent reaction-diffusion equation.}
     \begin{minipage}{\linewidth}
        \centering
        \captionof{subtable}{The top-3 solvers}
        \resizebox{\linewidth}{!}{
        \begin{tabular}{c|c|c|c|c|c}
\hline \hline 
 &Neural operator & Classical solver & Multi-grid & Relaxation & Strategies \\ 
\hline \hline
Top 1 solver & JacobiKAN & BiCGStab & 3-level & Jacobi & 3-1-3 \\ \hline
Top 2 solver & JacobiKAN & FGMRES & 3-level & SOR &  5-1-5 \\ \hline
Top 3 solver & JacobiKAN & CG & 3-level & SOR & 5-1-5 \\ \hline
\end{tabular}
}
    \end{minipage}

    \vspace{1em} 

    \begin{minipage}{\linewidth}
        \centering
        \captionof{subtable}{Performance and rank of performance of the top-3 solvers}
        \resizebox{\linewidth}{!}{ 
          \begin{tabular}{c|c|c|c|c|c|c|c}
\hline \hline 
 &Error & Com. time &  \# of ite. & Memory & MACs & Training time \\ 
\hline \hline
Top 1 solver & $6.67 \times 10^{-17}$ & 15.3513 & 512 & 78.335  & $2.93\times 10^{10}$ & 46384.5\\ \hline
Top 2 solver & $6.98 \times 10^{-17}$ & 13.9129 & 640 & 79.0752 & $2.53\times 10^{10}$ & 46384.5 \\ \hline
Top 3 solver & $6.92 \times 10^{-17}$ & 14.7723 & 642 & 77.54   & $2.69\times 10^{10}$ & 46384.5 \\ \hline
\end{tabular}
        }
    \end{minipage}
    \label{top3_rd_p2new}
\end{table}

\begin{table}[H]
    \centering
    \caption{Top-3 solvers evaluated by preference function $p^2$ for IMEX method, for solving 2-d time-dependent reaction-diffusion equation.}
     \begin{minipage}{\linewidth}
        \centering
        \captionof{subtable}{The top-3 solvers}
        \resizebox{\linewidth}{!}{
        \begin{tabular}{c|c|c|c|c|c}
\hline \hline 
 &Neural operator & Classical solver & Multi-grid & Relaxation & Strategies \\ 
\hline \hline
Top 1 solver & U-Net & BiCGStab & 2-level & Gauss-Seidel & 3-1-3 \\ \hline
Top 2 solver & DeepONet & BiCGStab & 2-level & Gauss-Seidel &  3-1-3 \\ \hline
Top 3 solver & U-Net & BiCGStab & 2-level & Gauss-Seidel & 5-1-5 \\ \hline
\end{tabular}
}
    \end{minipage}

    \vspace{1em} 

    \begin{minipage}{\linewidth}
        \centering
        \captionof{subtable}{Performance and rank of performance of the top-3 solvers}
        \resizebox{\linewidth}{!}{ 
          \begin{tabular}{c|c|c|c|c|c|c|c}
\hline \hline 
 &Error & Com. time &  \# of ite. & Memory & MACs & Training time \\ 
\hline \hline
Top 1 solver & $1.64 \times 10^{-14}$ & 5.11194  & 384 & 56.3916 & $2.55 \times 10^{11}$ & 10973  \\ \hline
Top 2 solver & $1.62 \times 10^{-14}$ & 5.34499  & 384 & 54.457  & $2.21 \times 10^{11}$ & 17235.9 \\ \hline
Top 3 solver & $2.17 \times 10^{-15}$ & 5.99003  & 384 & 55.8711 & $2.56 \times 10^{11}$ & 10973  \\ \hline
\end{tabular}
        }
    \end{minipage}
    \label{top3_rd_p2im}
\end{table}

\begin{figure}[H]
\centering
\subfigure[The top-3 meta-solvers discovered by preference function $p^1$.]{
\includegraphics[width=0.45\textwidth]{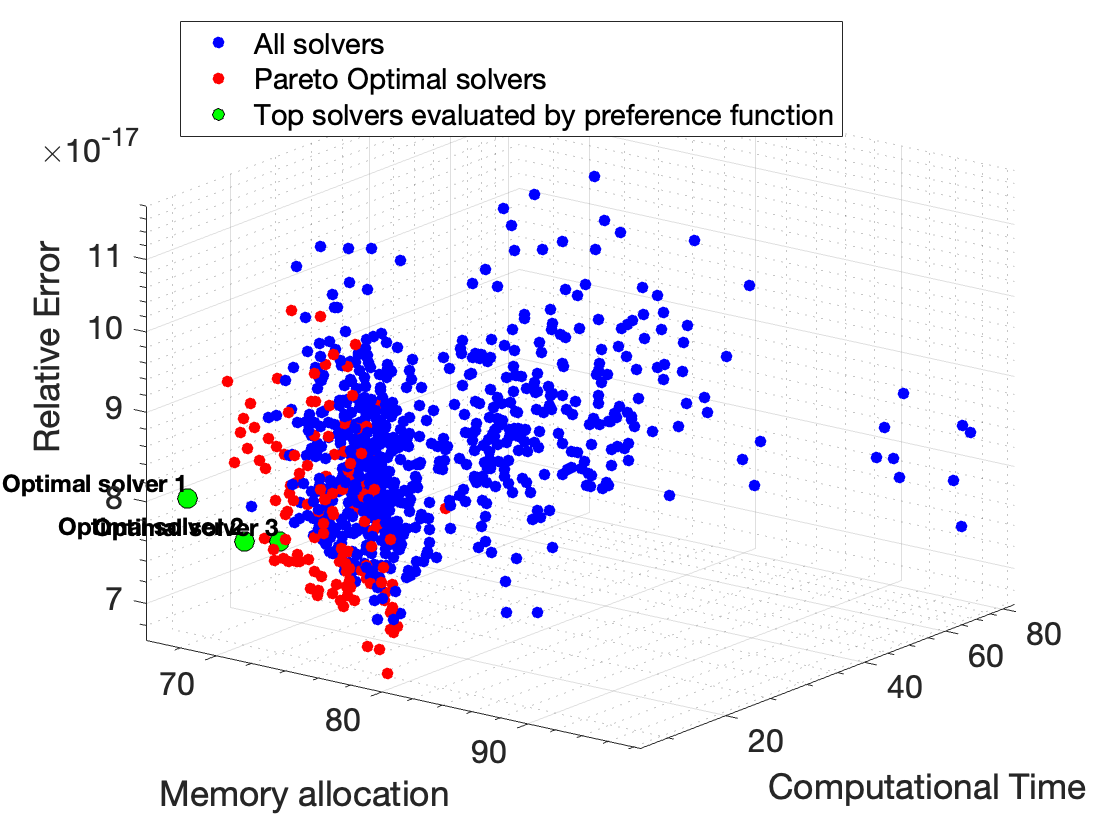}}
\subfigure[The top-3 meta-solvers discovered by preference function $p^2$.]{
\includegraphics[width=0.45\textwidth]{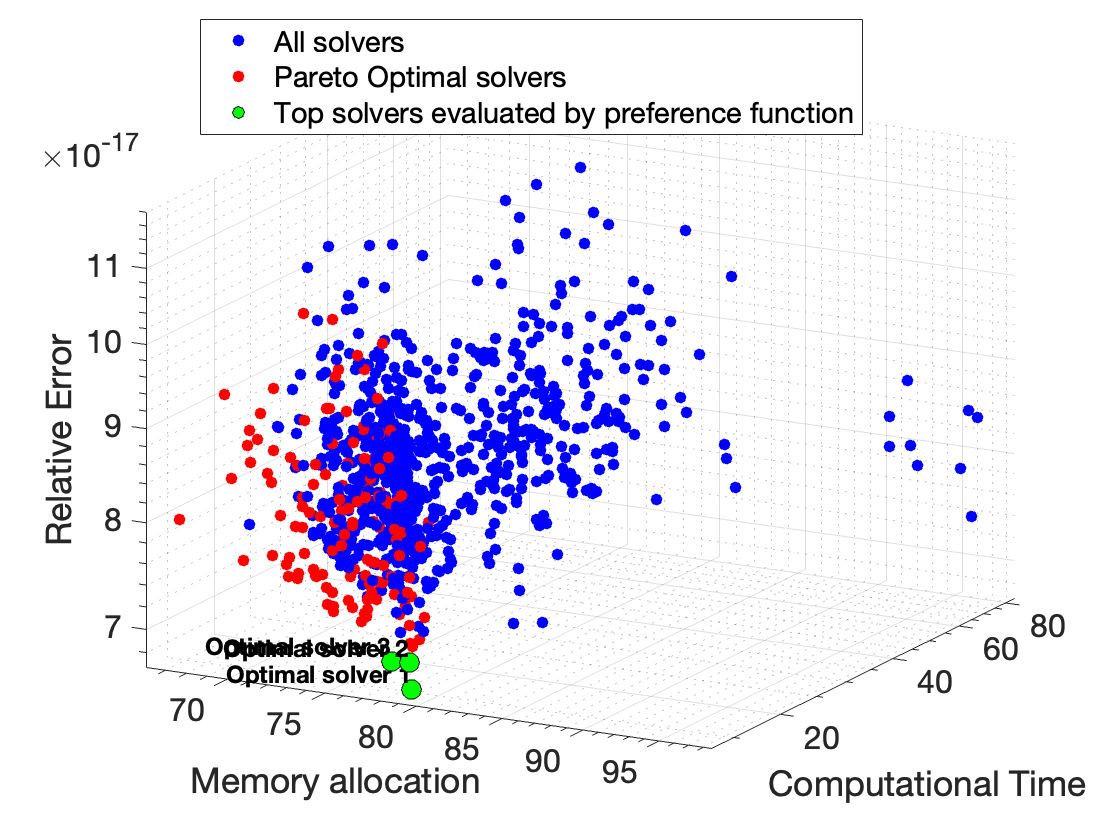}}
\caption{Pareto front: Computational Time -- Relative Error -- Memory allocation, with top 3 solvers evaluated by $p^1$ and $p^2$, for Newton-Raphson method.}
\label{discover_new}
\end{figure}

\begin{figure}[H]
\centering
\subfigure[The top-3 meta-solvers discovered by preference function $p^1$.]{
\includegraphics[width=0.45\textwidth]{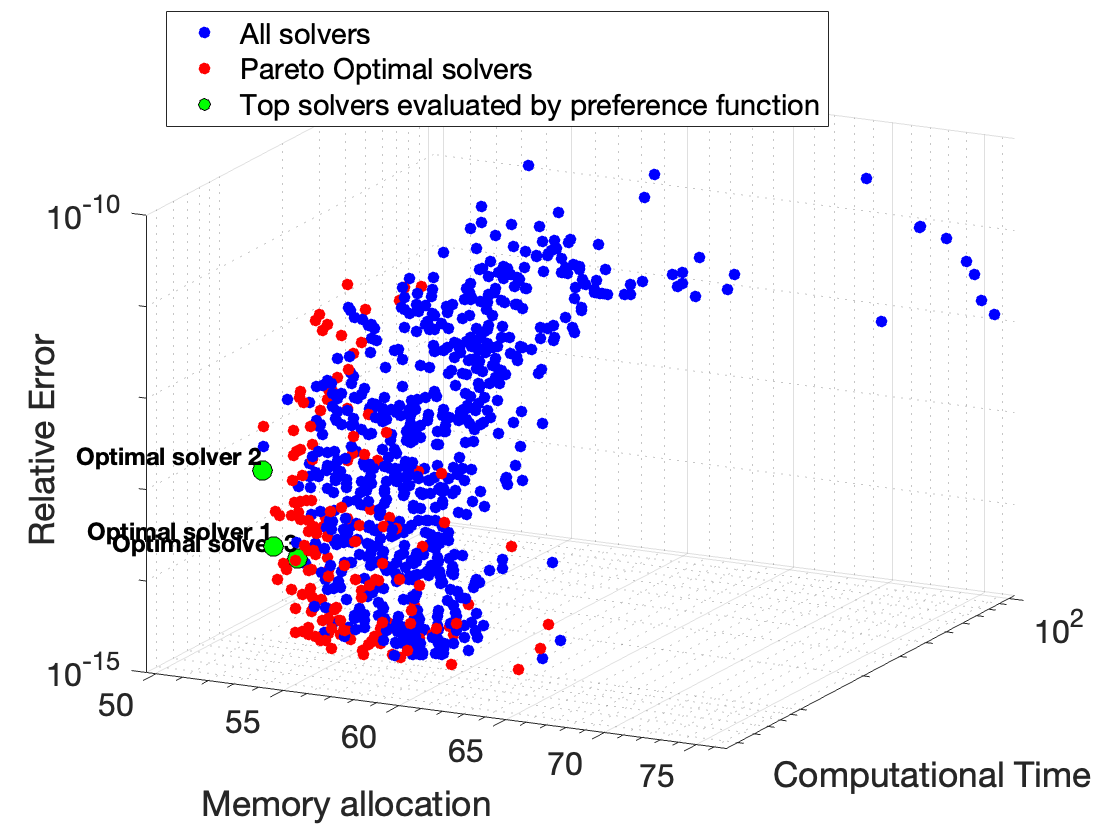}}
\subfigure[The top-3 meta-solvers discovered by preference function $p^2$.]{
\includegraphics[width=0.45\textwidth]{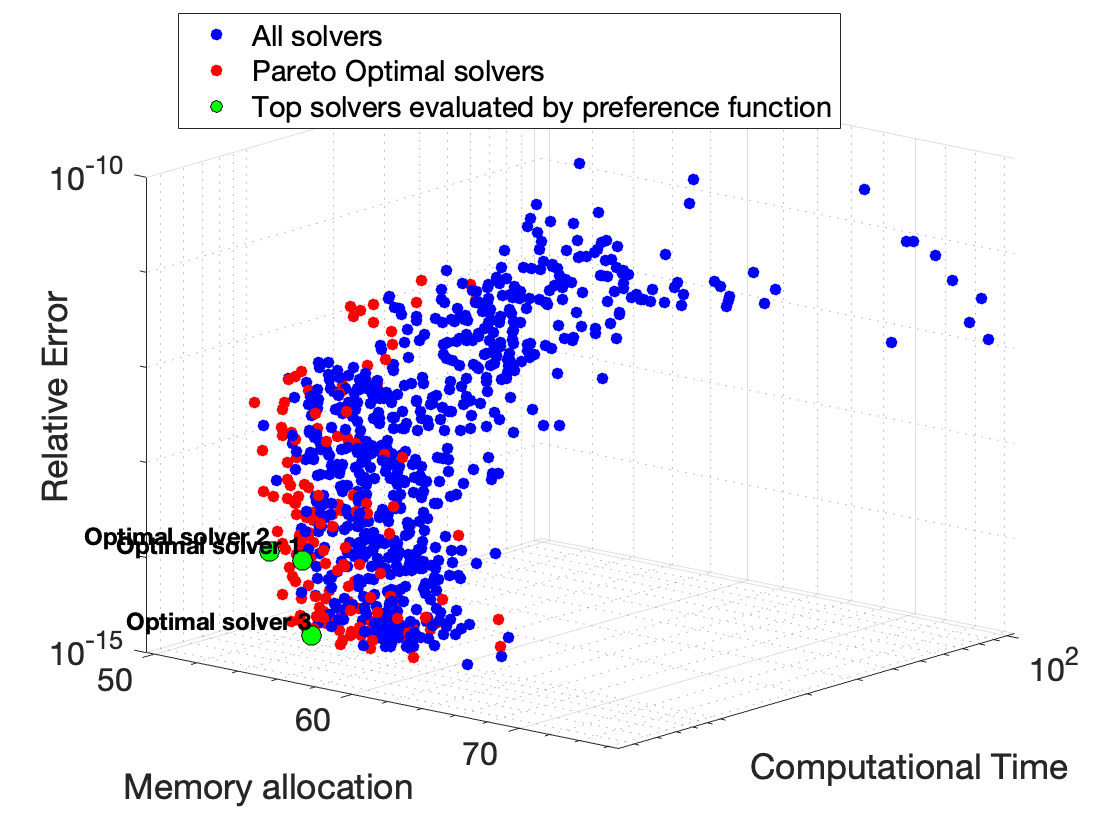}}
\caption{Pareto front: Computational Time -- Relative Error -- Memory allocation, with top 3 solvers evaluated by $p^1$ and $p^2$, for IMEX method.}
\label{discover_im}
\end{figure}

We plot the performance of the optimal meta-solvers discovered by different preference functions, in~\Cref{discover_new} for Newton-Raphspn method and in~\Cref{discover_im} for IMEX method. 

We implement the linear programming (LP) approach to rediscovery of solvers based on a particular type of preference functions, the weighted sum function of the rescaled performance. 

In the following we present the weights in the order s.t. $(\lambda_1,\lambda_2,\dots,\lambda_6)$ for ( $\lambda_1$: relative error, $\lambda_2$: computational time, $\lambda_3$: number of iterations, 
$\lambda_4$: memory allocation, $\lambda_5$: MACs, $\lambda_6$: Training time).

For Newton-Raphson method:

\begin{subequations}
\begin{equation}
\left\{
\begin{aligned}
&\text{Preference function $p(\lambda_a;r) = (0.2641, 0.4428, 0, 0.0938, 0.0718, 0.1274)^Tr$  } \\
&\text{Optimal solver}: x^a = (\text{U-Net}, \text{FGMRES}, \text{2-level}, \text{SSOR}, 5-1-5) \ .
\end{aligned}
\right.
\end{equation}

\begin{equation}
\left\{
\begin{aligned}
&\text{Preference function $p(\lambda_a;r) = (0.2395, 0, 0.0033, 0.0361, 0.3675, 0.3536)^Tr$ } \\
&\text{Optimal solver}: x^a = (\text{JacobiKAN}, \text{BiCGStab}, \text{2-level}, \text{SSOR}, 7-1-7) \ .
\end{aligned}
\right.
\end{equation}

\begin{equation}
\left\{
\begin{aligned}
&\text{Preference function $p(\lambda_a;r) = (0.0430, 0.7707, 0, 0.0301, 0.1561, 0)^Tr$ } \\
&\text{Optimal solver}: x^a = (\text{KAN}, \text{FGMRES}, \text{3-level}, \text{SSOR}, 3-1-3) \ .
\end{aligned}
\right.
\end{equation}

\end{subequations}

For IMEX method:

\begin{subequations}
\begin{equation}
\left\{
\begin{aligned}
&\text{Preference function $p(\lambda_a;r) = (0, 0.1485, 0, 0.6952, 0.1287, 0.0277)^Tr$  } \\
&\text{Optimal solver}: x^a = (\text{U-Net}, \text{BiCGStab}, \text{1-level}, \text{SSOR}, 9-1-9) \ .
\end{aligned}
\right.
\end{equation}

\begin{equation}
\left\{
\begin{aligned}
&\text{Preference function $p(\lambda_a;r) = (0.2510, 0.1022, 0.2329, 0.0979, 0.3161, 0)^Tr$ } \\
&\text{Optimal solver}: x^a = (\text{ChebyKAN}, \text{CG}, \text{2-level}, \text{SSOR}, 5-1-5) \ .
\end{aligned}
\right.
\end{equation}

\begin{equation}
\left\{
\begin{aligned}
&\text{Preference function $p(\lambda_a;r) = (0.2660, 0, 0, 0.1638, 0.4763, 0.0939)^Tr$ } \\
&\text{Optimal solver}: x^a = (\text{DeepONet}, \text{BiCGStab}, \text{2-level}, \text{SSOR}, 9-1-9) \ .
\end{aligned}
\right.
\end{equation}

\end{subequations}

However, not all solvers can be discovered using the weighted sum preference function. In terms of the formulation~\eqref{lp_discover}, it does not necessarily admit a feasible solution. Only the solvers in the boundary of the convex hull of Pareto optimal solvers can be rediscovered by the weighted sum preference function. 
As a preliminary quantification of the nonconvexity of the Pareto front, we count the number of Pareto optimal meta-solvers that lie in its nonconvex region, as described below. 
\begin{table}[H]
\centering
\begin{tabular}{c|c|c}
\hline
   & \# of solvers discovered by~\eqref{lp_discover}  & \# of solvers NOT  \\
  \hline
  Newton-Paphson method & 93   & 45 (38.6 \%) \\ \hline
  IMEX method & 120   & 41 (25.5 \%) \\ \hline
\end{tabular}
\end{table}

\subsection{Time-dependent Reaction Diffusion Equation with Small Correlation Length}

To show the generalizability of our methodology, we consider the time-dependent reaction diffusion equation with a different diffusion coefficient, corresponds to smaller correlation length in the system. 
In particular, we take the new diffusion coefficient to satisfying the following: 
\begin{equation}
k \sim N(10.0, K(\mathbf{x},\mathbf{x})=3.0e^{-\frac{\Vert \mathbf{x}-\mathbf{x}'\Vert^{2}}{2 \times 0.01^{2}}}), \text{ and } k \geq 1.0 \ . 
\end{equation}
The other implementation parameters remain the same as the previous equation.

We implement the same meta-solvers for both Newton-Raphson based method and IMEX based method. 
For the Newton-Raphson method, among 900 solvers, there are 167 Pareto optimal solvers.  
For the IMEX method, among 900 solvers, there are 107 Pareto optimal solvers. 
In the following, we summarize the composition of the set of Pareto optimal solvers, by counting the number of different components in the construction of meta-solvers. 

\begin{table}[H]
    \centering
    \caption{The composition of the set of Pareto optimal solvers by counting the number of elements in each dimension, for time-dependent reaction diffusion with small correlation length.}
    \begin{minipage}{\linewidth}
        \centering
        \captionof{subtable}{Different neural operators.}
\resizebox{\linewidth}{!}{
        \begin{tabular}{c| c c c  c c|c } 
\hline
Neural Op & DeepONet & U-Net  & KAN & JacobiKAN& ChebyKAN&Total \\
\hline
\# in Pareto opt for Newton  & 26 & 41 & 36 & 45 & 19 & 167\\
\hline
\# in Pareto opt for IMEX  & 15 & 22 & 23 & 32 & 15 & 107\\
\hline
\end{tabular}}

    \end{minipage}

    \vspace{1em} 

    \begin{minipage}{\linewidth}
        \centering
        \captionof{subtable}{Different Krylov solvers.}
           \begin{tabular}{c| c c c c    } 
\hline
Classical solvers & FGMRES & CG & BiCGStab \\
\hline
\# in Pareto opt for Newton  & 47 & 13 & 107  \\
\hline
\# in Pareto opt for IMEX  & 13 & 12 & 82  \\
\hline
\end{tabular}
        
    \end{minipage}
    
    \vspace{1em}

    \begin{minipage}{\linewidth}
        \centering
        \captionof{subtable}{Different smoothers.}
           \begin{tabular}{c| c c c c    } 
\hline
Smoother & GS & Jacobi & SOR & SSOR  \\
\hline
\# in Pareto opt for Newton & 40 & 23 & 52 & 52  \\
\hline
\# in Pareto opt for IMEX & 27 & 7 & 33 & 40  \\
\hline
\end{tabular}
        
    \end{minipage}

\vspace{1em}

\begin{minipage}{\linewidth}
        \centering
        \captionof{subtable}{Different strategies of applying smoothers.}
          \begin{tabular}{c| c c c c c    } 
\hline
Strategies for smoother & 1-1-1 & 3-1-3 & 5-1-5 & 7-1-7 & 9-1-9  \\
\hline
\# in Pareto opt for Newton & 16 & 46 & 38 & 39 & 28 \\
\hline
\# in Pareto opt for IMEX & 17 & 28 & 29 & 17 & 16 \\
\hline
\end{tabular}
        
    \end{minipage}
    
\vspace{1em}

\begin{minipage}{\linewidth}
        \centering
        \captionof{subtable}{Different levels of multigrid.}
    \begin{tabular}{c| c c c } 
\hline
Levels in multi-grid & 1-level & 2-level &  3-level   \\
\hline
\# in Pareto opt for Newton & 3 & 96 & 68  \\
\hline
\# in Pareto opt for IMEX & 7 & 64 & 36  \\
\hline
\end{tabular}

       \end{minipage}
    \label{rd_com2}
\end{table}

Moreover, we plot the Pareto fronts for both Newton-Paphson method, in~\Cref{front_2drd_new}, and IMEX method, in~\Cref{front_2drd_im}, for the computational time, relative error and number of iteration, for the relative error, MACs and Memory allocation, and for the  Relative error, convergence rate and MACs.  
One can find that the shapes of the Pareto fronts preserve high similarity compared to the previous equation, showing the generalizability of the Pareto optimility analysis proposed in our methodology.  

\begin{figure}[H]
\centering
\subfigure[Computational time -- Relative error -- \# of iterations]{
\includegraphics[width=0.45\textwidth]{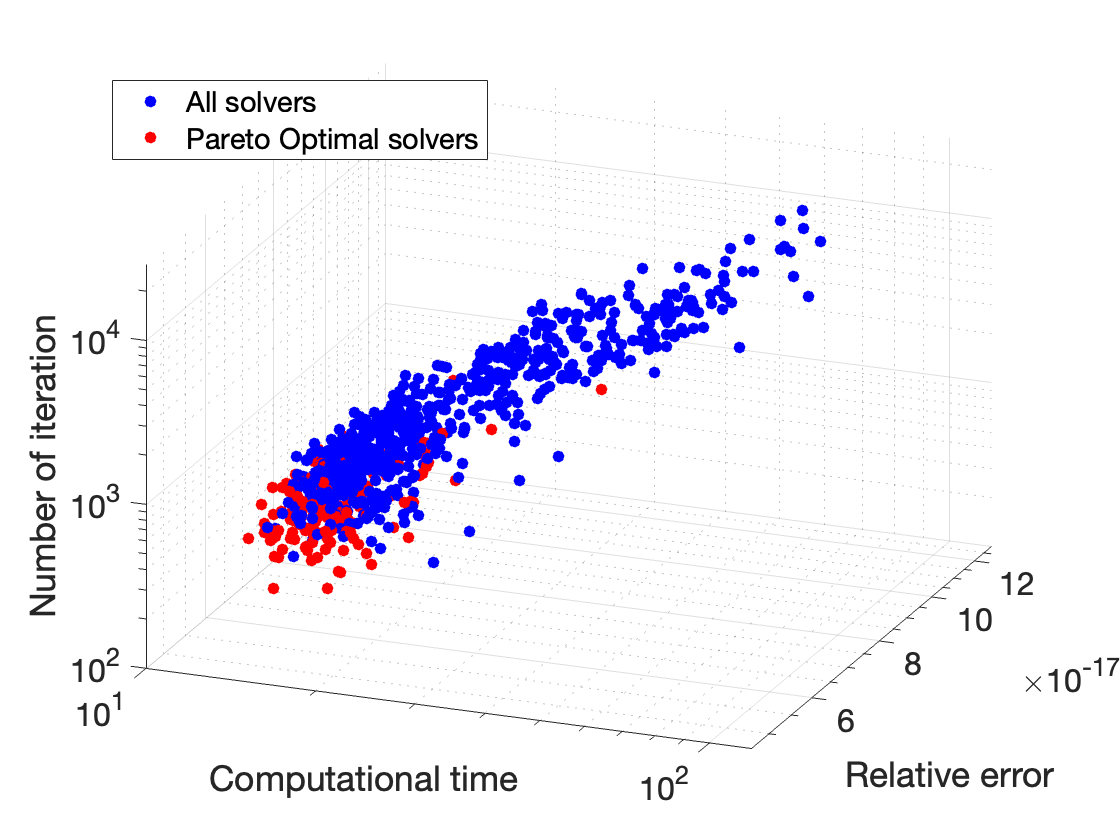}}
\subfigure[Relative error -- MACs -- Memory allocation]{
\includegraphics[width=0.45\textwidth]{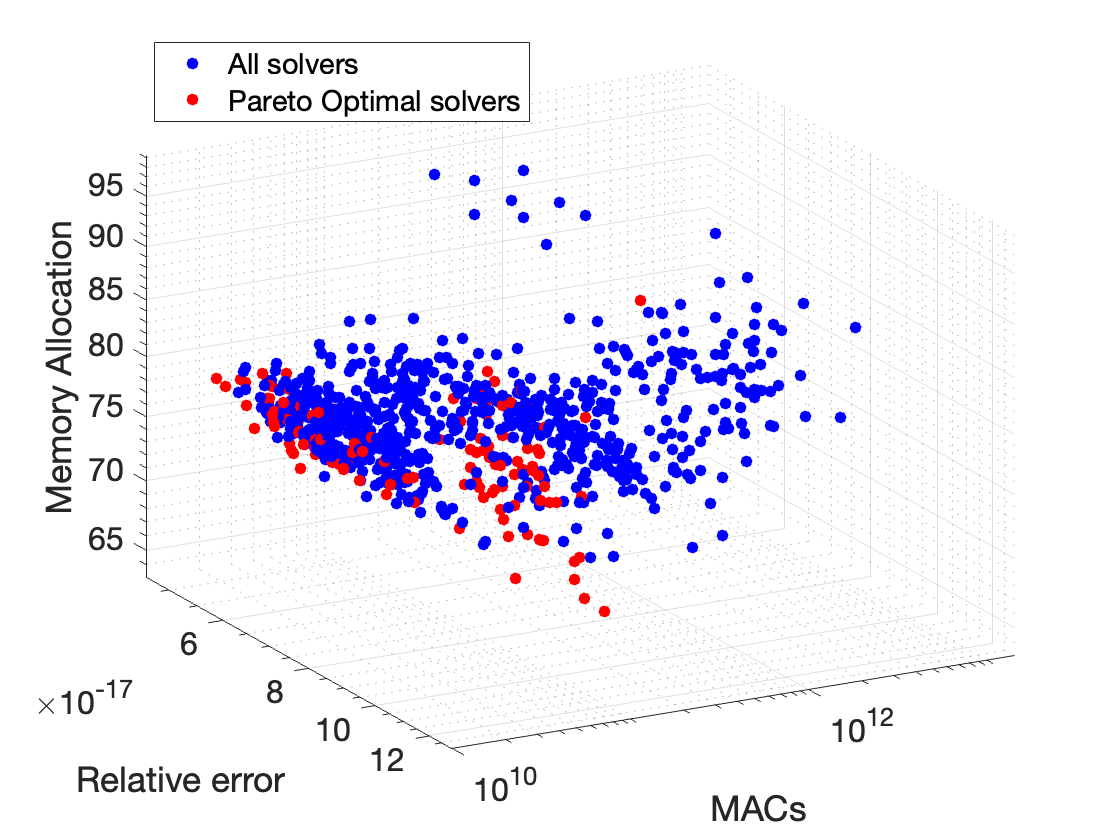}}
\subfigure[Relative error -- MACs -- \# of iterations]{
\includegraphics[width=0.45\textwidth]{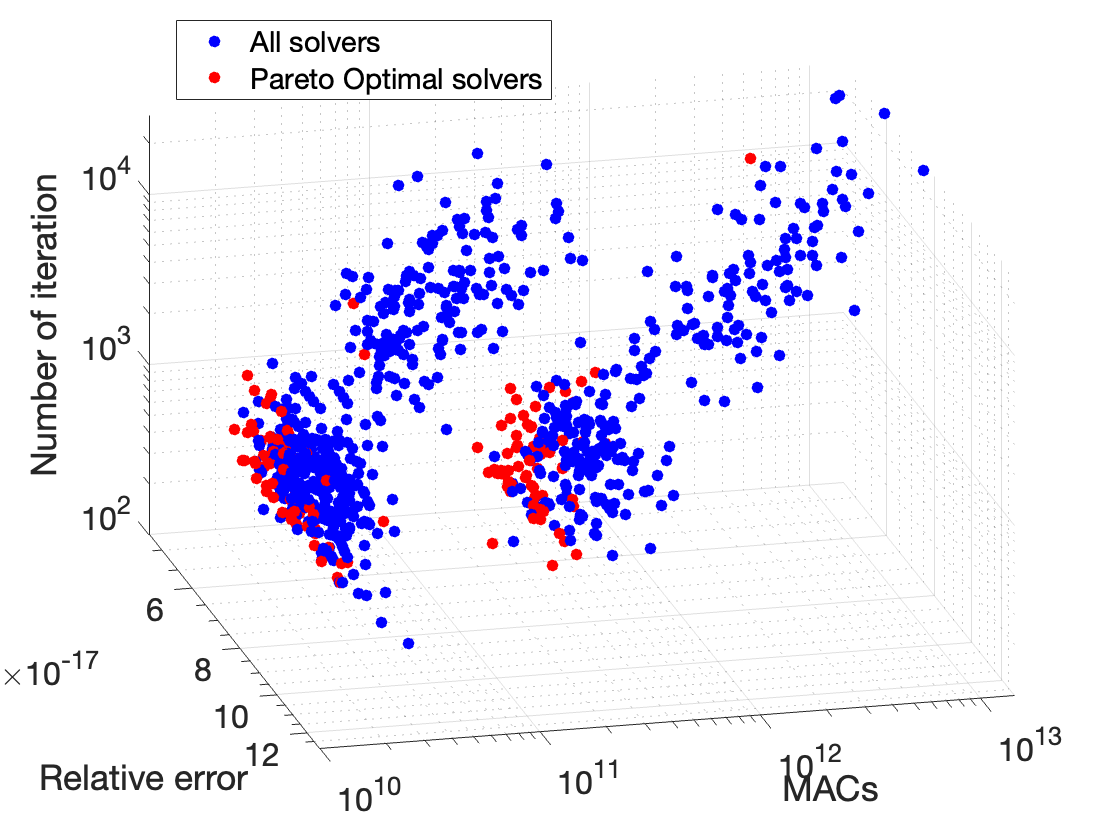}}
\caption{Pareto fronts for solving 2-D time-dependent reaction-diffusion equation, with small correlation length, using Newton-Raphson method. All solvers are depicted in blue while Pareto optimal solvers are highlighted in red. The``gap'' due to the adaption thus the improvement of performance by multi-grid techniques.}
\label{front_2drd_new}
\end{figure}

\begin{figure}[H]
\centering
\subfigure[Computational time -- Relative error -- \# of iterations]{
\includegraphics[width=0.45\textwidth]{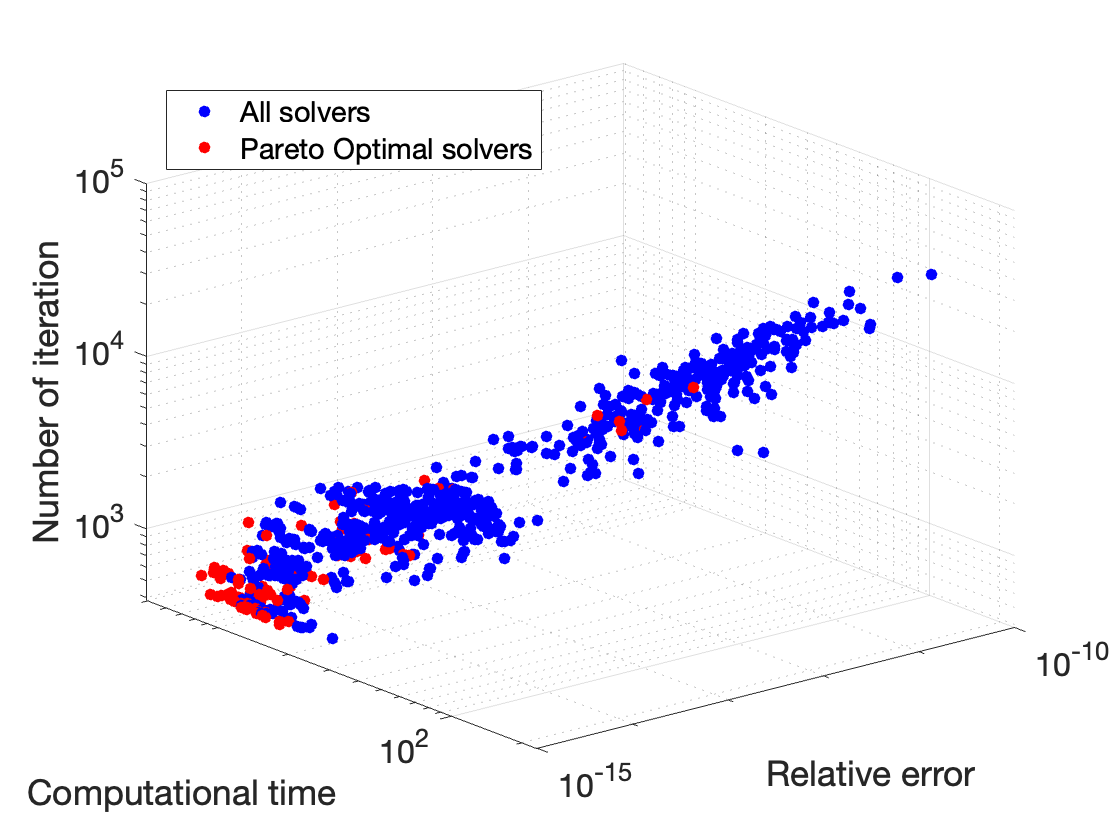}}
\subfigure[Relative error -- MACs -- Memory allocation]{
\includegraphics[width=0.45\textwidth]{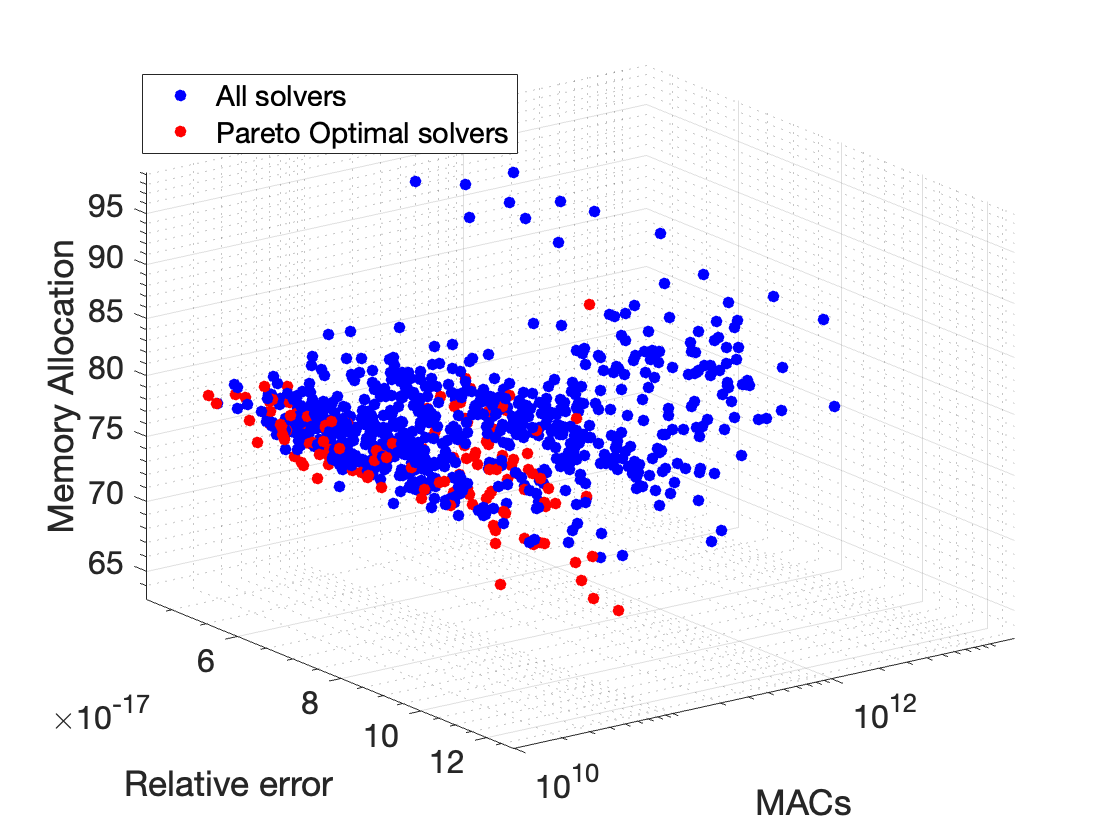}}
\subfigure[Relative error -- MACs -- \# of iterations]{
\includegraphics[width=0.45\textwidth]{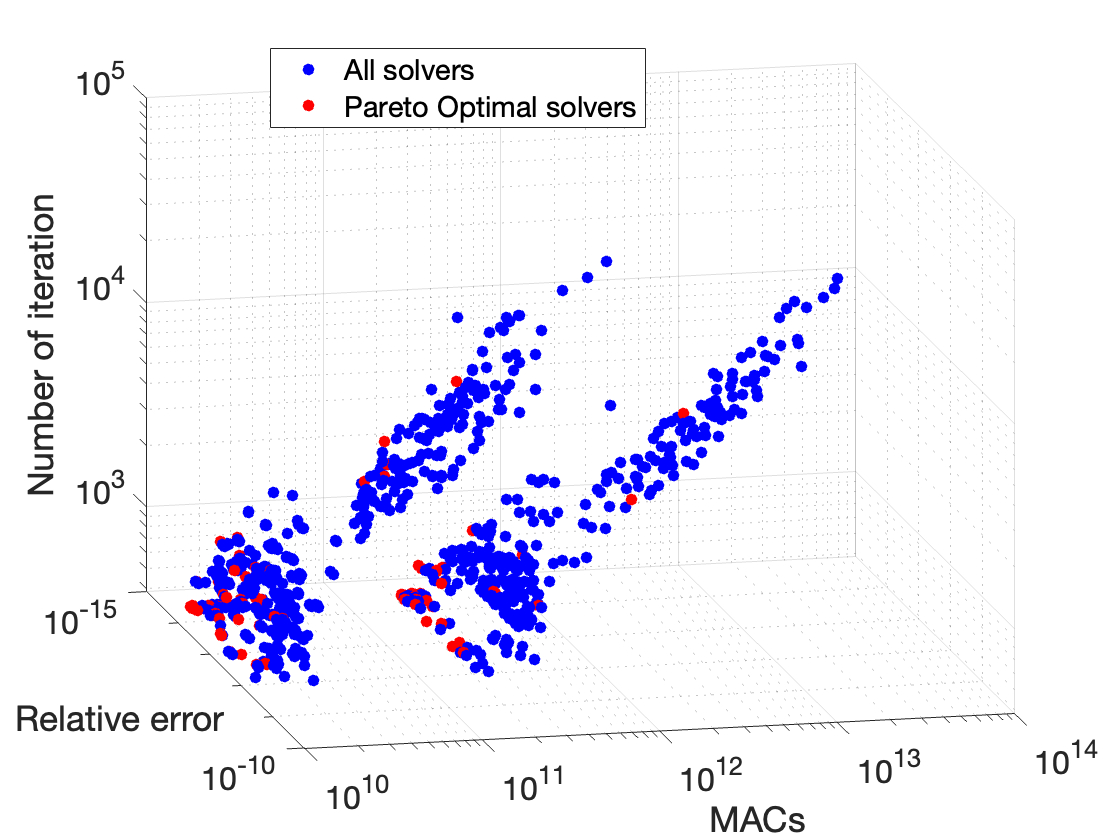}}
\caption{Pareto fronts for solving 2-D time-dependent reaction-diffusion equation, with small correlation length, using IMEX method. All solvers are depicted in blue while Pareto optimal solvers are highlighted in red. The``gap'' due to the adaption thus the improvement of performance by multi-grid techniques.}
\label{front_2drd_im}
\end{figure}

We apply the re-scale function~\eqref{rescale} on the performance data. Then, we use the same preference functions to discover the optimal meta-solvers. In particular, for preference 
\begin{equation}
p^1(r) = \frac{1}{6}(\sum_{i=1}^6 r_i)\ ,
\end{equation}
the average of all the ranks, the results are presented in following. 

\begin{table}[H]
    \centering
    \caption{Top-3 solvers evaluated by preference function $p^1$ for Newton-Raphson method, solving 2-d reaction diffusion equation with small correlation length. }
     \begin{minipage}{\linewidth}
        \centering
        \captionof{subtable}{The top-3 solvers}
        \resizebox{\linewidth}{!}{
        \begin{tabular}{c|c|c|c|c|c}
\hline \hline 
 &Neural operator & Classical solver & Multi-grid & Relaxation & Strategies \\ 
\hline \hline
    Top 1 solver & U-Net & BiCGStab & 3-level & Jacobi & 3-1-3 \\ \hline
Top 2 solver & U-Net & BiCGStab & 2-level & Jacobi &  5-1-5 \\ \hline
Top 3 solver & U-Net & BiCGStab & 3-level & Gauss-Seidel & 7-1-7 \\ \hline
\end{tabular}
}
    \end{minipage}

    \vspace{1em} 

    \begin{minipage}{\linewidth}
        \centering
        \captionof{subtable}{Performance and rank of performance of the top-3 solvers}
        \resizebox{\linewidth}{!}{ 
          \begin{tabular}{c|c|c|c|c|c|c}
\hline \hline 
 &Error & Com. time  & \# of ite & Memory & MACs & Training time \\ 
\hline \hline
Top 1 solver & $6.16\times 10 ^{-17}$ & 15.7099 & 520 & 63.0859 & $3.95\times 10^{11} $   & 9010.38 \\ \hline
Top 2 solver & 6.21 $\times 10^{-17}$ & 17.1465 & 512 & 66.1445 & $3.5 \times 10^{11}$  & 9010.38  \\ \hline
Top 3 solver &  $6.38\times 10^{-17}$ & 15.9111 & 384 & 67.873  & $2.96 \times 10^{11} $& 9010.38 \\\hline
\end{tabular}
        }
    \end{minipage}
    \label{top3_p1_new_rd2}
\end{table}

\begin{table}[H]
    \centering
    \caption{Top-3 solvers evaluated by preference function $p^1$ for IMEX method, solving 2-d reaction diffusion equation with small correlation length.}
     \begin{minipage}{\linewidth}
        \centering
        \captionof{subtable}{The top-3 solvers}
        \resizebox{\linewidth}{!}{
        \begin{tabular}{c|c|c|c|c|c}
\hline \hline 
 &Neural operator & Classical solver & Multi-grid & Relaxation & Strategies \\ 
\hline \hline
Top 1 solver  & U-Net & BiCGStab & 2-level & SSOR & 1-1-1 \\ \hline
Top 2 solver  & U-Net & BiCGStab & 2-level & SOR &  3-1-3 \\ \hline
Top 3 solver  & U-Net & BiCGStab & 2-level & Jacobi & 5-1-5 \\ \hline
\end{tabular}
}
    \end{minipage}

    \vspace{1em} 

    \begin{minipage}{\linewidth}
        \centering
        \captionof{subtable}{Performance and rank of performance of the top-3 solvers}
        \resizebox{\linewidth}{!}{ 
          \begin{tabular}{c|c|c|c|c|c|c}
\hline \hline 
 &Error & Com. time  & \# of ite & Memory & MACs & Training time \\ 
\hline \hline
Top 1 solver & $4.02 \times 10^{-15}$ & 6.28384  & 512 & 51.6621 & $3.38\times 10^{11}$ & 10973 \\ \hline
Top 2 solver & 1.02 $\times 10^{-13}$ & 5.58047  & 384 & 52.249  & $2.55\times 10^{11}$ & 10973 \\ \hline
Top 3 solver & $3.66 \times 10^{-15}$ & 9.5998   & 512 & 52.3398 & $3.44\times 10^{11}$ & 10973 \\ \hline
\end{tabular}
        }
    \end{minipage}
    \label{top3_p1__im_rds}
\end{table}

For the preference 
\begin{equation}
    p^2(r) =  \frac{1}{24} \big( 10 (r_1 + r_2) +  (r_3 + r_4 + r_5 + r_6 ) \big) \ ,
\end{equation}
the results are presented in the following. 

\begin{table}[H]
    \centering
    \caption{Top-3 solvers evaluated by preference function $p^2$ for Newton-Raphson method, solving 2-d reaction diffusion equation with small correlation length.}
     \begin{minipage}{\linewidth}
        \centering
        \captionof{subtable}{The top-3 solvers}
        \resizebox{\linewidth}{!}{
        \begin{tabular}{c|c|c|c|c|c}
\hline \hline 
 &Neural operator & Classical solver & Multi-grid & Relaxation & Strategies \\ 
\hline \hline
Top 1 solver & U-Net & FGMRES & 2-level & Gauss-Seidel & 9-1-9 \\ \hline
Top 2 solver & DeepONet & BiCGStab & 2-level & Gauss-Seidel &  1-1-1 \\ \hline
Top 3 solver & U-Net & FGMRES & 2-level & SSOR & 1-1-1 \\ \hline
\end{tabular}
}
    \end{minipage}

    \vspace{1em} 

    \begin{minipage}{\linewidth}
        \centering
        \captionof{subtable}{Performance and rank of performance of the top-3 solvers}
        \resizebox{\linewidth}{!}{ 
          \begin{tabular}{c|c|c|c|c|c|c|c}
\hline \hline 
 &Error & Com. time &  \# of ite. & Memory & MACs & Training time \\ 
\hline \hline
Top 1 solver & $5.15 \times 10^{-17}$ & 15.3767  & 582 & 72.2656 & $3.88\times 10^{11}$ & 9010.38 \\ \hline
Top 2 solver & $5.07\times 10^{-17}$ & 15.011   & 653 & 73.8701 & $3.78\times 10^{11}$ & 17270.6 \\  \hline
Top 3 solver & $5.40\times 10^{-17}$ & 13.8969  & 768 & 75.8438 & $5.05\times 10^{11}$ & 9010.38 \\ \hline
\end{tabular}
        }
    \end{minipage}
    \label{top3_p2_new_rds}
\end{table}

\begin{table}[H]
    \centering
    \caption{Top-3 solvers evaluated by preference function $p^2$ for IMEX method, solving 2-d reaction diffusion equation with small correlation length.}
     \begin{minipage}{\linewidth}
        \centering
        \captionof{subtable}{The top-3 solvers}
        \resizebox{\linewidth}{!}{
        \begin{tabular}{c|c|c|c|c|c}
\hline \hline 
 &Neural operator & Classical solver & Multi-grid & Relaxation & Strategies \\ 
\hline \hline
Top 1 solver & U-Net & BiCGStab & 2-level & SSOR & 1-1-1 \\ \hline
Top 2 solver & U-Net & BiCGStab & 2-level & SOR &  3-1-3 \\ \hline
Top 3 solver & DeepONet & BiCGStab & 2-level & Gauss-Seidel & 5-1-5 \\ \hline
\end{tabular}
}
    \end{minipage}

    \vspace{1em} 

    \begin{minipage}{\linewidth}
        \centering
        \captionof{subtable}{Performance and rank of performance of the top-3 solvers}
        \resizebox{\linewidth}{!}{ 
          \begin{tabular}{c|c|c|c|c|c|c|c}
\hline \hline 
 &Error & Com. time &  \# of ite. & Memory & MACs & Training time \\ 
\hline \hline
Top 1 solver & $4.02\times 10^{-15}$ & 6.28384  & 512  & 51.6621 & $3.38\times 10^{11}$  & 10973  \\ \hline
Top 2 solver & $1.02\times 10^{-13}$ & 5.58047  & 384  & 52.249  & $2.55\times 10^{11}$  & 10973   \\ \hline
Top 3 solver & $2.09\times 10^{-14}$ & 6.16931  & 384  & 53.8867 & $2.22\times 10^{11}$  & 17235.9 \\ \hline
\end{tabular}
        }
    \end{minipage}
    \label{top3_p2_im_rds}
\end{table}

\begin{figure}[H]
\centering
\subfigure[The top-3 meta-solvers discovered by preference function $p^1$.]{
\includegraphics[width=0.45\textwidth]{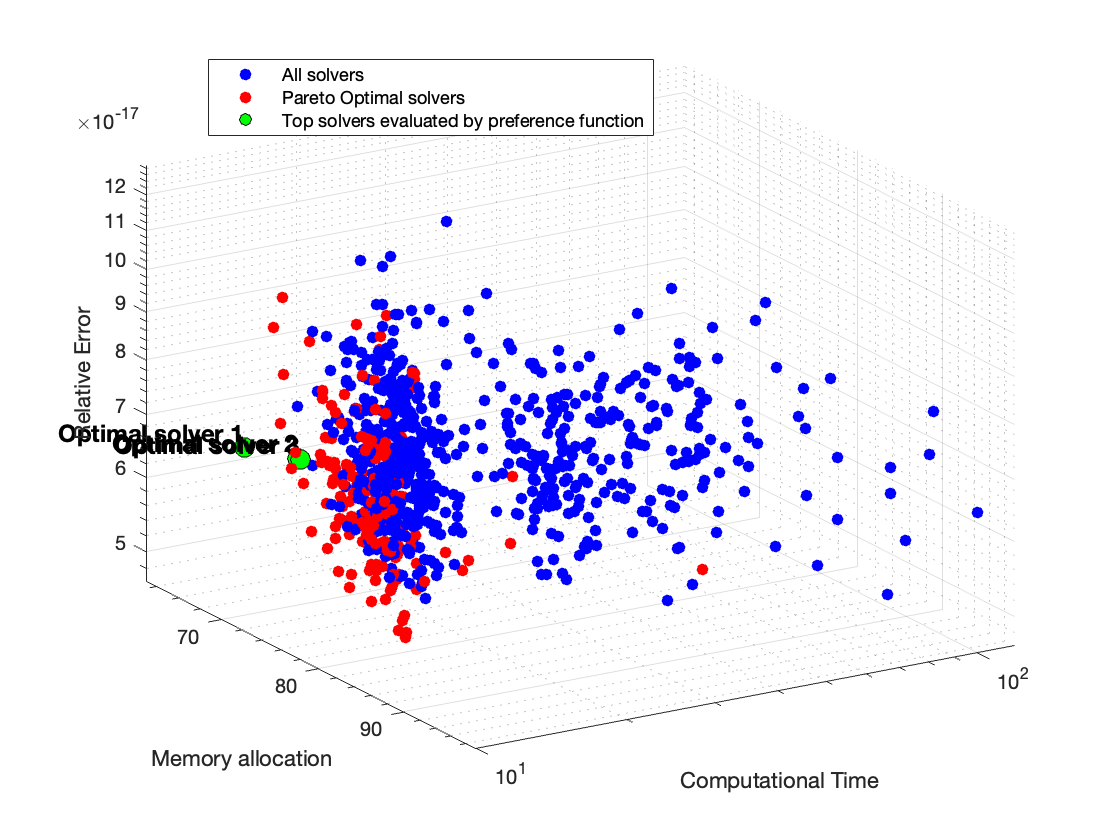}}
\subfigure[The top-3 meta-solvers discovered by preference function $p^2$.]{
\includegraphics[width=0.45\textwidth]{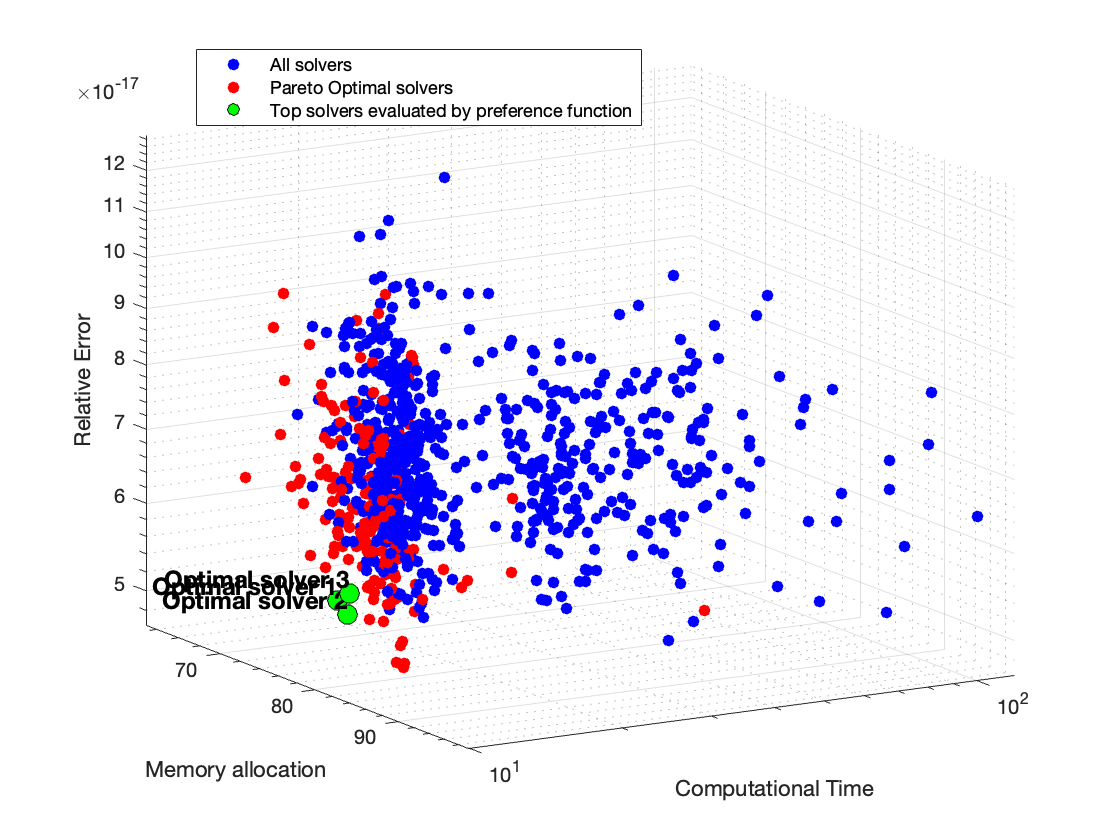}}
\caption{Pareto front: Computational Time -- Relative Error -- Memory allocation, with top 3 solvers evaluated by $p^1$ and $p^2$, for Newton-Raphson method, solving 2-d reaction diffusion equation with small correlation length.}
\label{discover_new_s}
\end{figure}

\begin{figure}[H]
\centering
\subfigure[The top-3 meta-solvers discovered by preference function $p^1$.]{
\includegraphics[width=0.45\textwidth]{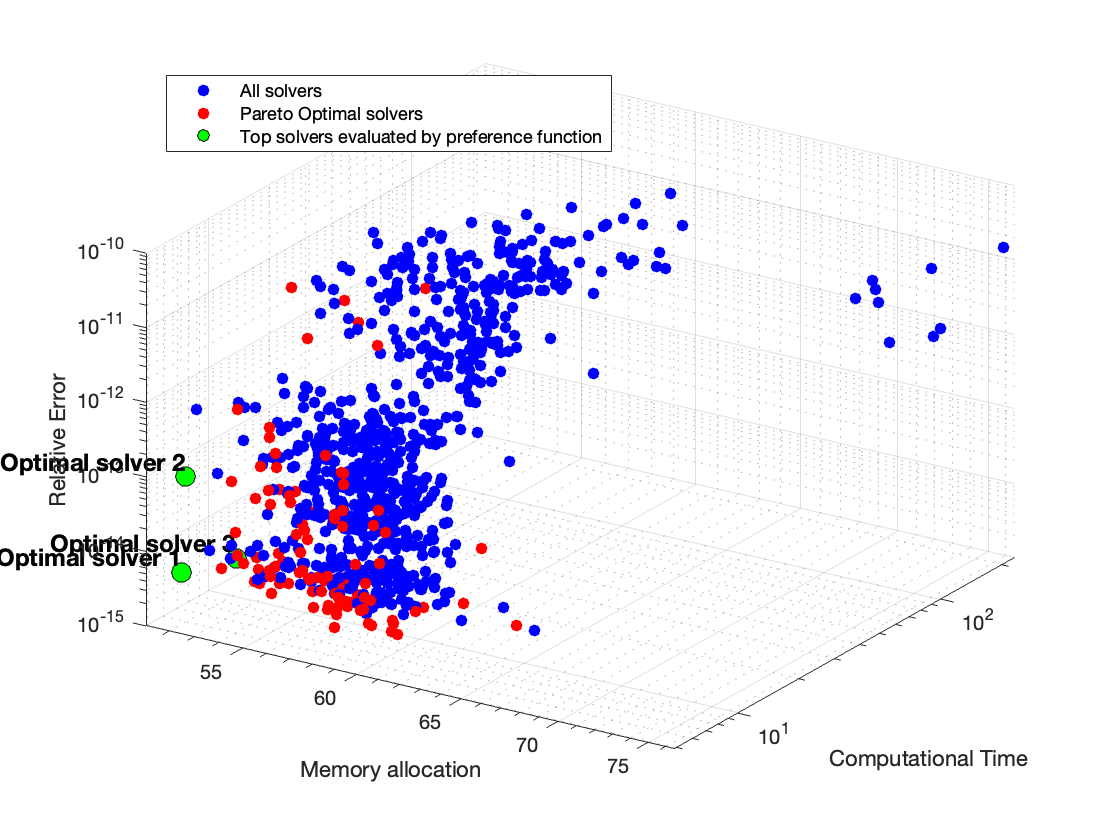}}
\subfigure[The top-3 meta-solvers discovered by preference function $p^2$.]{
\includegraphics[width=0.45\textwidth]{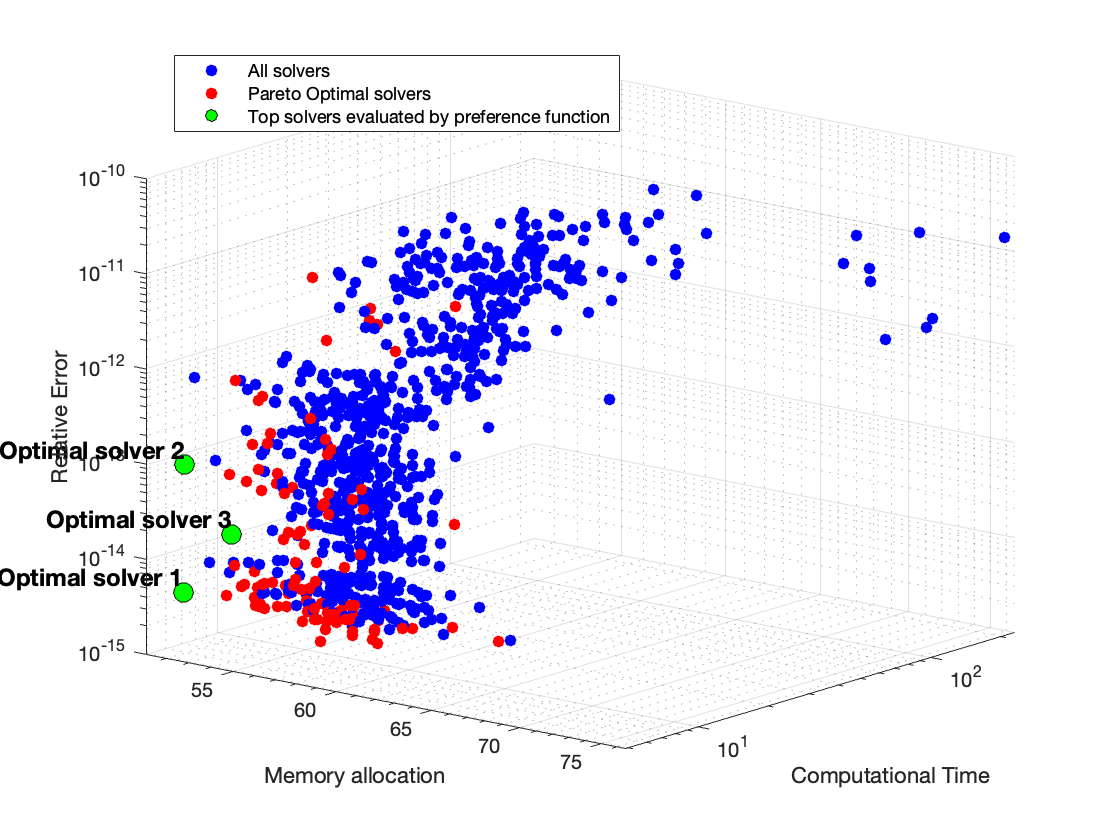}}
\caption{Pareto front: Computational Time -- Relative Error -- Memory allocation, with top 3 solvers evaluated by $p^1$ and $p^2$, for IMEX method, solving 2-d reaction diffusion equation with small correlation length.}
\label{discover_im_s}
\end{figure}

We plot the performance of the optimal meta-solvers discovered by different preference functions, in~\Cref{discover_new_s} for Newton-Raphspn method and in~\Cref{discover_im_s} for IMEX method. 

\section{Implementation details and further numerical results for solving incompressible Navier-Stokes equations}\label{app_ns}

In this section, we consider a vector-valued nonlinear equation, the incompressible Navier Stokes equation. Given domain $\Omega = [0,1]^3 \subseteq \R^{3}$, the incompressible Naiver Stokes equation is given by
\begin{equation}
\left\{\begin{split}
\rho \frac{\partial \mathbf{u}}{\partial t} + \rho (\mathbf{u} \cdot \nabla)\mathbf{u}&= -\nabla p + \mu \Delta \mathbf{u}, \,\,\, \text{ in } \Omega \times (0,1], \\
\nabla \cdot \mathbf{u} &= 0, \,\,\, \text{ in } \Omega \times (0,1], \\
\mathbf{u} &= \mathbf{u}_{0}, \text{ at } t=0, \\
\mathbf{u} &= \mathbf{u}_{D}, \text{ on } \partial\Omega \times (0,1],
\end{split}
\right.
\label{app:eqn:imcompressible-navier-stokes}
\end{equation}
where $\rho = \mu = 1$. 
In this case, the exact velocity $\mathbf{u}=(u_1, u_2, u_3)$ and the pressure $p$ of~\eqref{app:eqn:imcompressible-navier-stokes} are given by
\begin{align*}
    u_{1}(x, y, z, t) &= -a [e^{ax}\sin(ay+dz)+e^{az}\cos(ax+dy)]e^{-d^{2}t}, \\
    u_{2}(x, y, z, t) &= -a [e^{ay}\sin(az+dx)+e^{ax}\cos(ay+dz)]e^{-d^{2}t}, \\
    u_{3}(x, y, z, t) &= -a [e^{az}\sin(ax+dy)+e^{ay}\cos(az+dx)]e^{-d^{2}t}, \\
    p(x, y, z, t) &= -\frac{a^{2}}{2} [e^{2ax} + e^{2ay} + e^{2az} +  2 \sin(ax+dy)\cos(az+dx)e^{a(y+z)} \\
    &+ 2 \sin(ay+dz)\cos(ax+dy)]e^{a(z+x)} + 2 \sin(az+dx)\cos(ay+dz)]e^{a(x+y)}]e^{-2d^{2}t},
\end{align*}
where $a$ and $d$ are given. 

To solve \eqref{app:eqn:imcompressible-navier-stokes}, we employ the \textit{splitting method} instead of solving \textit{saddle point problem}.
We briefly describe the incremental pressure correction scheme~(IPCS).
\begin{subequations}
    
\begin{enumerate}
    \item Find $\mathbf{u}^{\ast}$ using IMEX:
    \begin{equation}
        \langle \rho\frac{\mathbf{u}^{\ast}-\mathbf{u}^{(n)}}{\Delta t}, \mathbf{v} \rangle  + AB_{3}(F, \mathbf{u}^{(n)}, \mathbf{u}^{(n-1)}, \mathbf{u}^{(n-2)}) + \langle \nabla p^{(n)} , \mathbf{v} \rangle + \langle \frac{\mu}{2} \nabla (\mathbf{u}^{(n)} + \mathbf{u}^{\ast}), \nabla \mathbf{v} \rangle = 0,
    \end{equation}
    where $F(\mathbf{u}) := \langle \rho (\mathbf{u} \cdot \nabla) \mathbf{u}, \mathbf{v} \rangle$.
    \item Find $p^{(n+1)}$ satisfying
    \begin{align}
        \langle \nabla p^{(n+1)}, \nabla q\rangle &= \langle \nabla p^{(n)}, \nabla q\rangle - \frac{\rho}{\Delta t}\langle \nabla \cdot \mathbf{u}^{\ast}, q\rangle, \\
        \frac{\partial p^{(n+1)}}{\partial n} &= 0.
    \end{align}
    This relation is derived from
    \begin{equation}
        \rho \frac{\mathbf{u}^{(n+1)}-\mathbf{u}^{\ast}}{\Delta t} + \nabla p^{(n+1)} - \nabla p^{(n)} = 0,
    \end{equation}
    where we use the condition $\nabla \cdot \mathbf{u}^{(n+1)}$.
    \item Find $\mathbf{u}^{(n+1)}$ satisfying
    \begin{equation}
        \rho \langle \mathbf{u}^{(n+1)}-\mathbf{u}^{\ast}, \mathbf{v} \rangle = -\Delta t \langle \nabla(p^{(n+1)}-p^{(n)}), \mathbf{v} \rangle.
    \end{equation}
\end{enumerate}

\end{subequations}
Moreover, we considered two different values of d. 
The implementation details are given in following: For training:
\begin{itemize}
    \item Case 1: $a_{i} \sim \mathcal{U}[\frac{\pi}{2}-0.1, \frac{\pi}{2}+0.1]$ and $d_{i} \sim \mathcal{U}[\frac{\pi}{4}-0.1, \frac{\pi}{4}+0.1]$.
    \item Case 2: $a_{i} \sim \mathcal{U}[{\pi}-0.1, {\pi}+0.1]$ and $d_{i} \sim \mathcal{U}[\frac{\pi}{8}-0.1, \frac{\pi}{8}+0.1]$.
    \item We sampled $3,000$ $(a_{i}, d_{i})$ pairs.
    \item Spatial mesh size: $\Delta x = \frac{1}{14}$.
    \item Time mesh size: $\Delta t = \frac{C  \Delta x}{\lVert \mathbf{u} \rVert}$, $C=0.7$, $\Delta t = \frac{1}{252}$.
    \item We randomly sampled $10$ time snapshots for each $(a_{i},d_{i})$ pair.
    \item Total number of training samples: $30,000$.
    \item Polynomial order for $\mathbf{u}$ and $p$ is $1$. 
\end{itemize}

For the test case:

\begin{itemize}
    \item Case 1: $a = \frac{\pi}{2}$ and $d=\frac{\pi}{4}$.
    \item Case 2: $a = {\pi}$ and $d=\frac{\pi}{8}$.
    \item Spatial mesh size: $\Delta x = \frac{1}{32}$.
    \item Time mesh size: $\Delta t = \frac{C  \Delta x}{\lVert \mathbf{u} \rVert}$, $C=0.7$, $\Delta t = \frac{1}{1215}$.
    \item We solved the equation in $\Omega \times [0, 1]$.
    \item Polynomial order for $\mathbf{u}$ and $p$ is $1$. 
\end{itemize}

\subsection{Identify all optimal meta-solvers in Pareto sense}

We implement the parametrized IMEX based meta-solvers. In particular, we consider $\NO= \{\text{DeepONet, U-DeepONet,}$ \ $ \text{KAN, JacobiKAN,  CheyKAN}\}$, in the choose of neural operators. 
For the classical Krylov methods, we consider $\kr = \{ \text{FGMRES, CG, BiCGStab} \}$. Moreover, for relaxation method we consider $\rela = \{\text{Jacobi, Gauss-Seidel, SOR, SSOR}\}$, and for the strategies of applying the relaxation method, we consider $S = \{\text{1-1-1, 3-1-3, 5-1-5, 7-1-7, 9-1-9}\}$. Moreover, we also implement the multi-grid techniques, where the levels of multi-grid is choose from $\{1,2,3\}$. 
As a results, we generate 900 different meta-solvers. 

For the performnace matrices, the evaluation is conducted in a eight-dimensional performance function $f = (f_1,\dots, f_8)$, where each dimension represents a particular performance criterion. In particular, $f_1$ is the computational time, $f_2$ is the relative error of velocity $\textbf{u}$,  $f_3$ is the relative error of pressure $\textbf{p}$, $f_4$ is the number of iterations, $f_5$ is the memory allocation, $f_6$ is the MACs, $f_7$ is the average MACs and $f_8$ is the training time, respectively. 

For $d=\frac{\pi}{4}$, there are 258 Pareto optimal solvers among 900. For $d=\frac{\pi}{8}$, there are 247 Pareto optimal solvers among 900. 
In the following, we first summarize the composition of the set of Pareto optimal solvers, by counting the number of different components in the construction of meta-solvers. 

\begin{table}[H]
    \centering
    \caption{The composition of the set of Pareto optimal solvers by counting the number of elements in each dimension, for 3d incompressible Navier-Stokes equation.}
    \begin{minipage}{\linewidth}
        \centering
        \captionof{subtable}{Different neural operators.}
        \begin{tabular}{c| c c c  c c|c } 
        \hline
Neural Op & DeepONet & U-Net  & KAN & JacobiKAN& ChebyKAN & Total\\
\hline
\# in Pareto opt, $d = \frac{\pi}{4}$ & 62 & 60 & 39 & 54 & 43 & 258\\
\hline
\# in Pareto opt, $d = \frac{\pi}{8}$ & 55 & 54 & 44 & 51 & 43 & 247  \\
\hline
\end{tabular}

    \end{minipage}

    \vspace{1em} 

    \begin{minipage}{\linewidth}
        \centering
        \captionof{subtable}{Different Krylov solvers.}
           \begin{tabular}{c| c c c c    } 
\hline
Classical solvers & FGMRES & CG & BiCGStab \\
\hline
\# in Pareto opt, $d = \frac{\pi}{4}$ & 43 & 129 & 86  \\
\hline
\# in Pareto opt, $d = \frac{\pi}{8}$ & 48 & 112 & 87  \\
\hline
\end{tabular}
        
    \end{minipage}
    
    \vspace{1em}

    \begin{minipage}{\linewidth}
        \centering
        \captionof{subtable}{Different smoothers.}
           \begin{tabular}{c| c c c c    } 
\hline
Smoother & GS & Jacobi & SOR & SSOR  \\
\hline
\# in Pareto opt, $d = \frac{\pi}{4}$  & 85 & 1 & 51 & 121  \\
\hline
\# in Pareto opt, $d = \frac{\pi}{8}$ & 74 & 3 & 54 & 116  \\
\hline
\end{tabular}
        
    \end{minipage}

\vspace{1em}

\begin{minipage}{\linewidth}
        \centering
        \captionof{subtable}{Different strategies of applying smoothers.}
          \begin{tabular}{c| c c c c c    } 
\hline
Strategies for smoother & 1-1-1 & 3-1-3 & 5-1-5 & 7-1-7 & 9-1-9  \\
\hline
\# in Pareto opt, $d = \frac{\pi}{4}$& 49 & 60 & 50 & 51 & 48 \\
\hline
\# in Pareto opt, $d = \frac{\pi}{8}$& 37 & 59 & 55 & 43 & 53 \\
\hline
\end{tabular}
        
    \end{minipage}
    
\vspace{1em}

\begin{minipage}{\linewidth}
        \centering
        \captionof{subtable}{Different levels of multigrid.}
    \begin{tabular}{c| c c c } 
\hline
Levels in multi-grid & 1-level & 2-level &  3-level   \\
\hline
\# in Pareto opt, $d = \frac{\pi}{4}$ & 47 & 125 & 86  \\
\hline
\# in Pareto opt, $d = \frac{\pi}{8}$ & 58 & 128 & 61  \\
\hline
\end{tabular}

       \end{minipage}
    \label{ns_3d_com}
\end{table}

We plot the Pareto fronts, for MACs, average MACs, for computational time, relative error (p), number of iterations, for relative error (p), memory allocation, MACs and for computational time, relative error (p), memory allocation. 

\begin{figure}[H]
\centering
\subfigure[Computational time -- relative error (p) -- number of iterations]{
\includegraphics[width=0.45\textwidth]{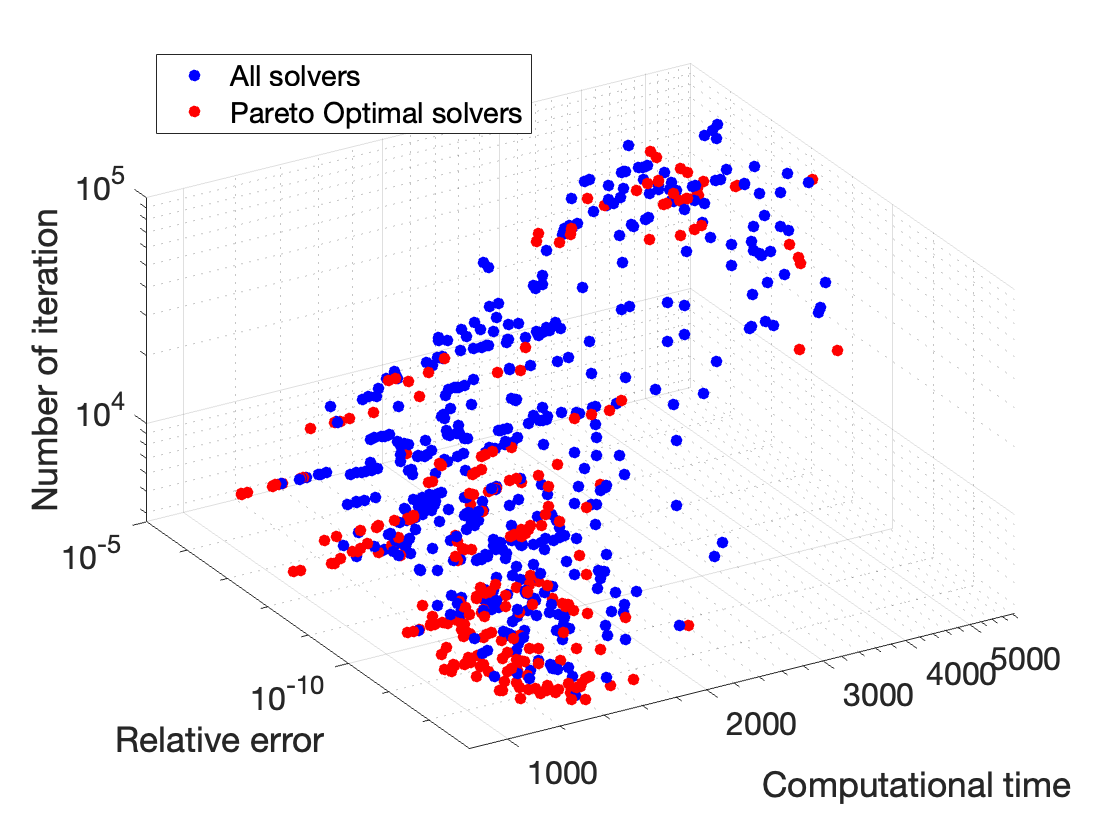}}
\subfigure[Relative error (p) -- memory allocation -- MACs]{
\includegraphics[width=0.45\textwidth]{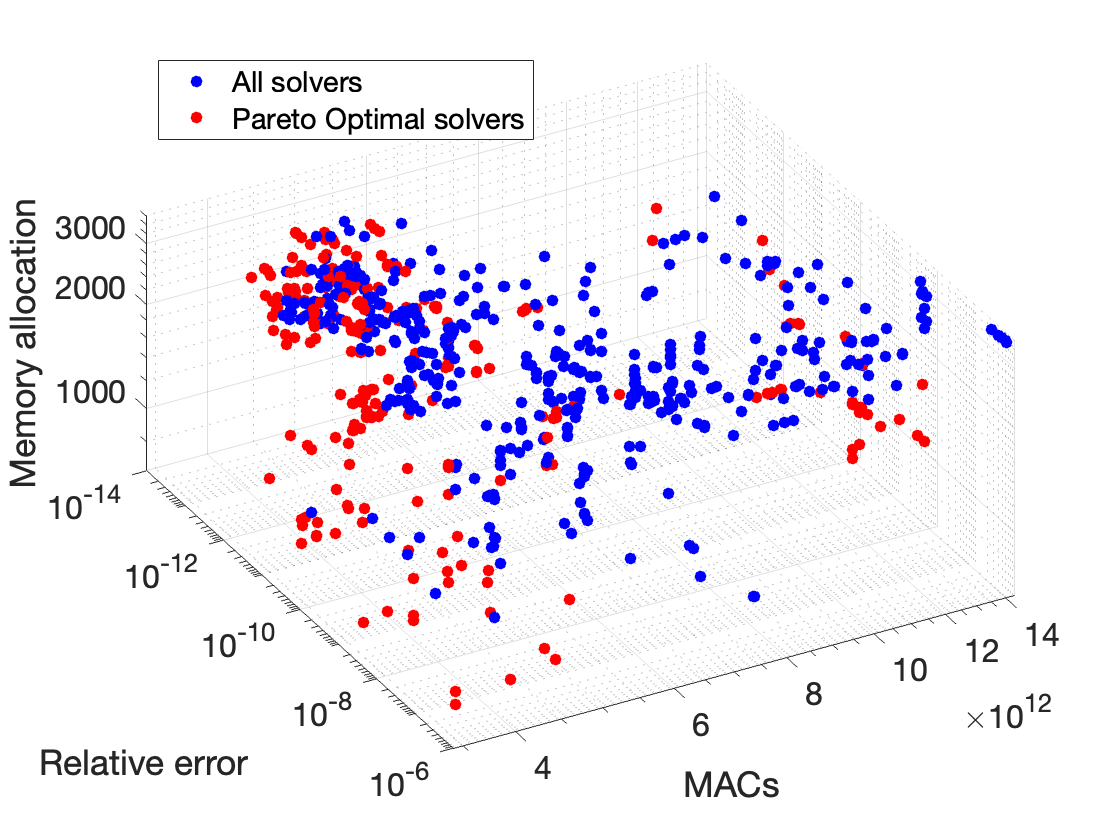}}
\subfigure[computational time -- relative error (p) -- memory allocation]{
\includegraphics[width=0.45\textwidth]{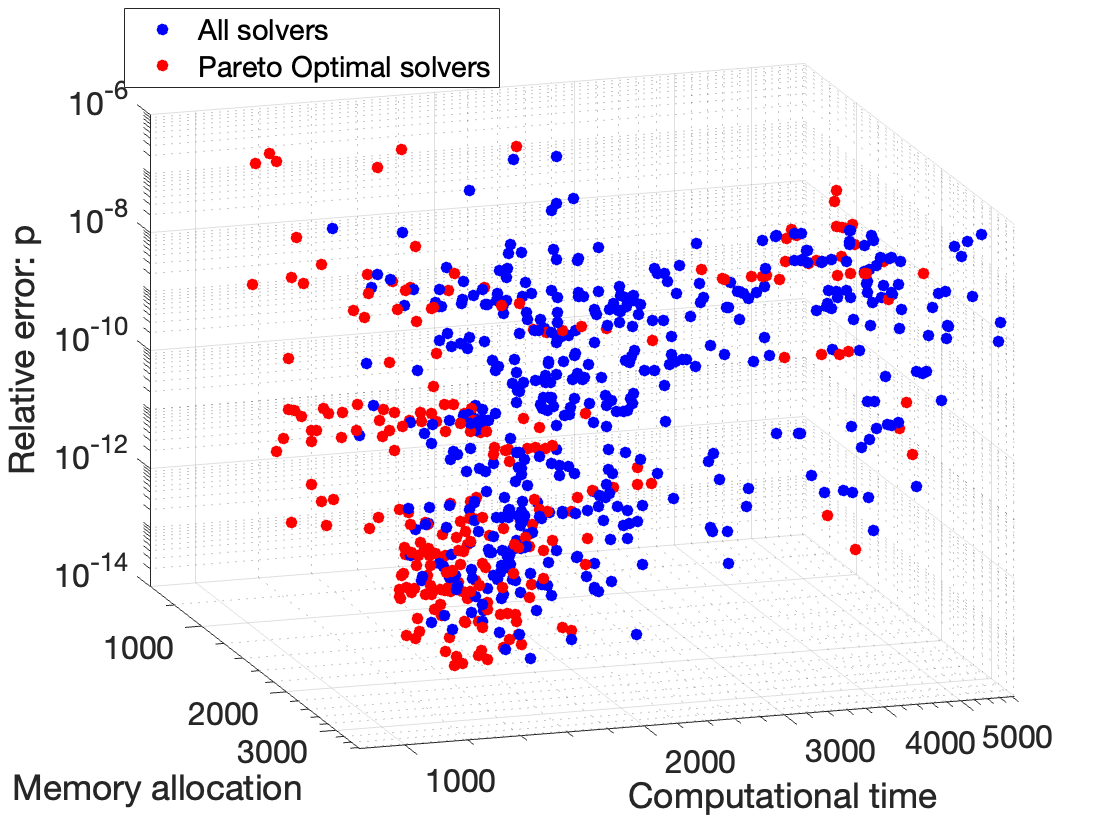}}
\subfigure[computational time -- MACs -- Average MACs]{
\includegraphics[width=0.45\textwidth]{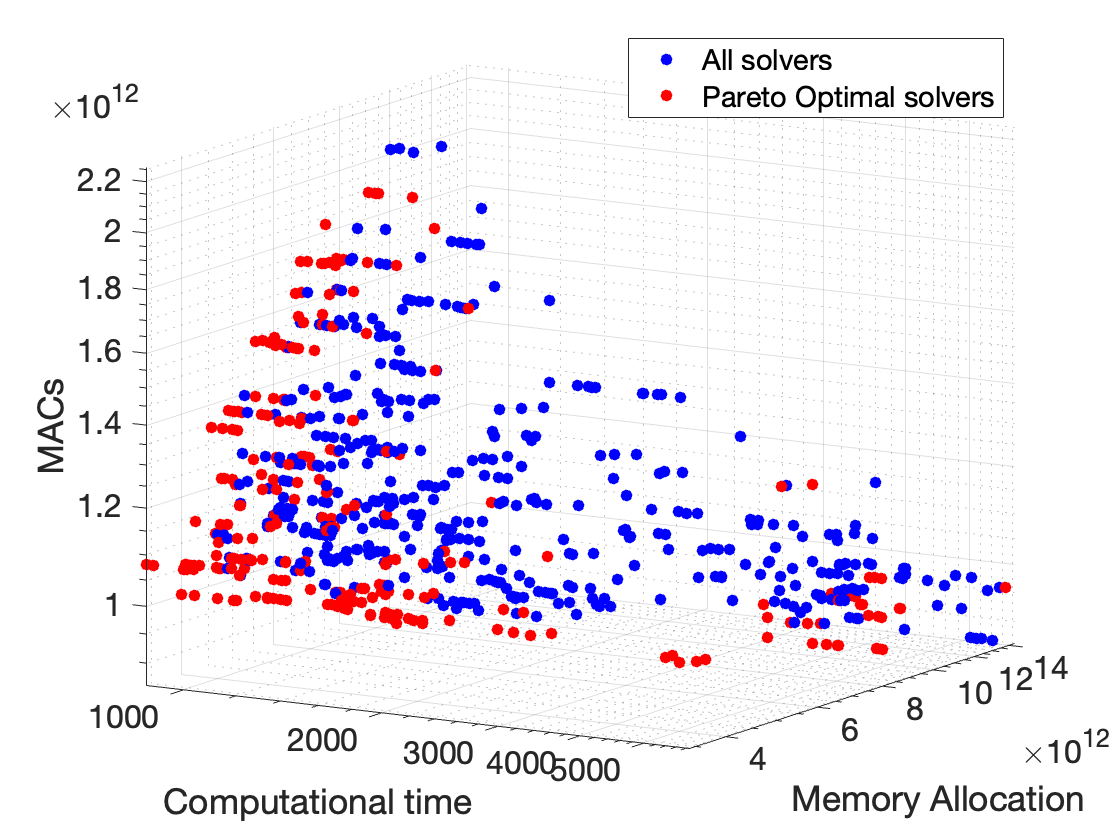}}
\caption{Pareto fronts for solving 3-D incompressive Navier-Stokes equation, using Newton-Raphson method, when $d = \frac{\pi}{4}$. All solvers are depicted in blue while Pareto optimal solvers are highlighted in red.}
\label{front_3ndsd1}
\end{figure}

\newpage

\begin{figure}[H]
\centering
\subfigure[Computational time -- relative error (p) -- number of iterations]{
\includegraphics[width=0.45\textwidth]{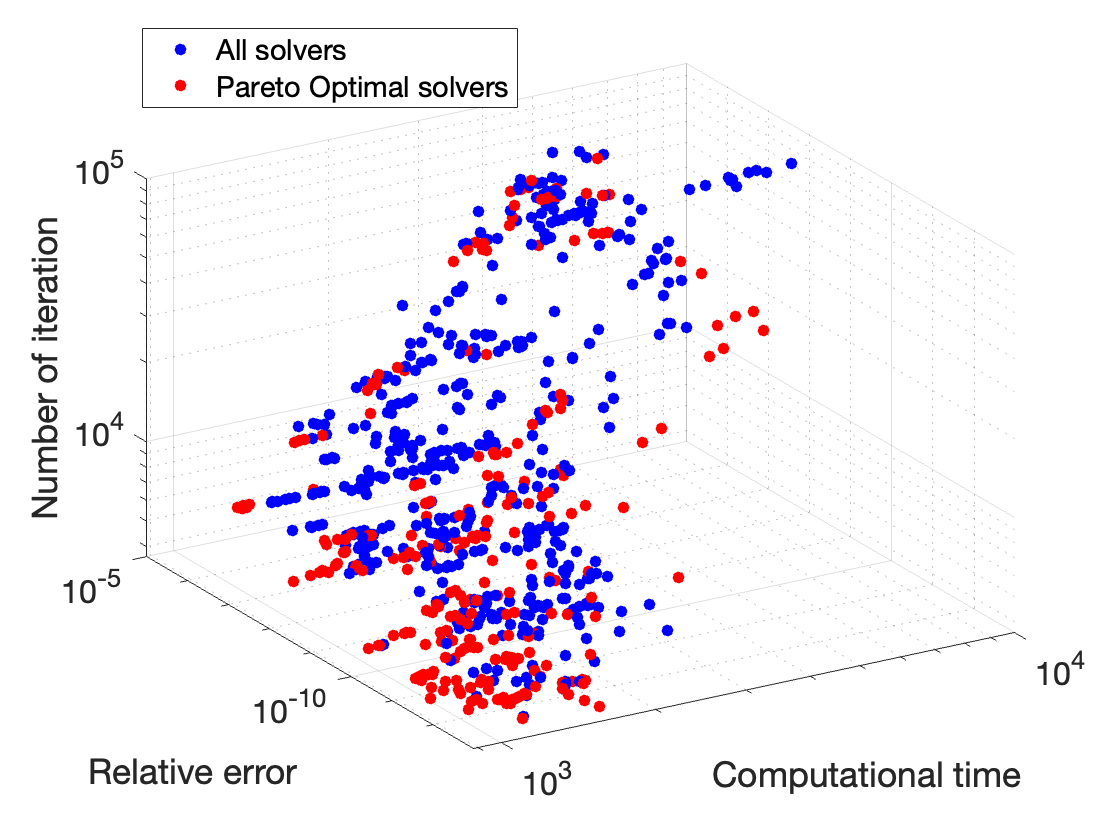}}
\subfigure[Relative error (p) -- memory allocation -- MACs]{
\includegraphics[width=0.45\textwidth]{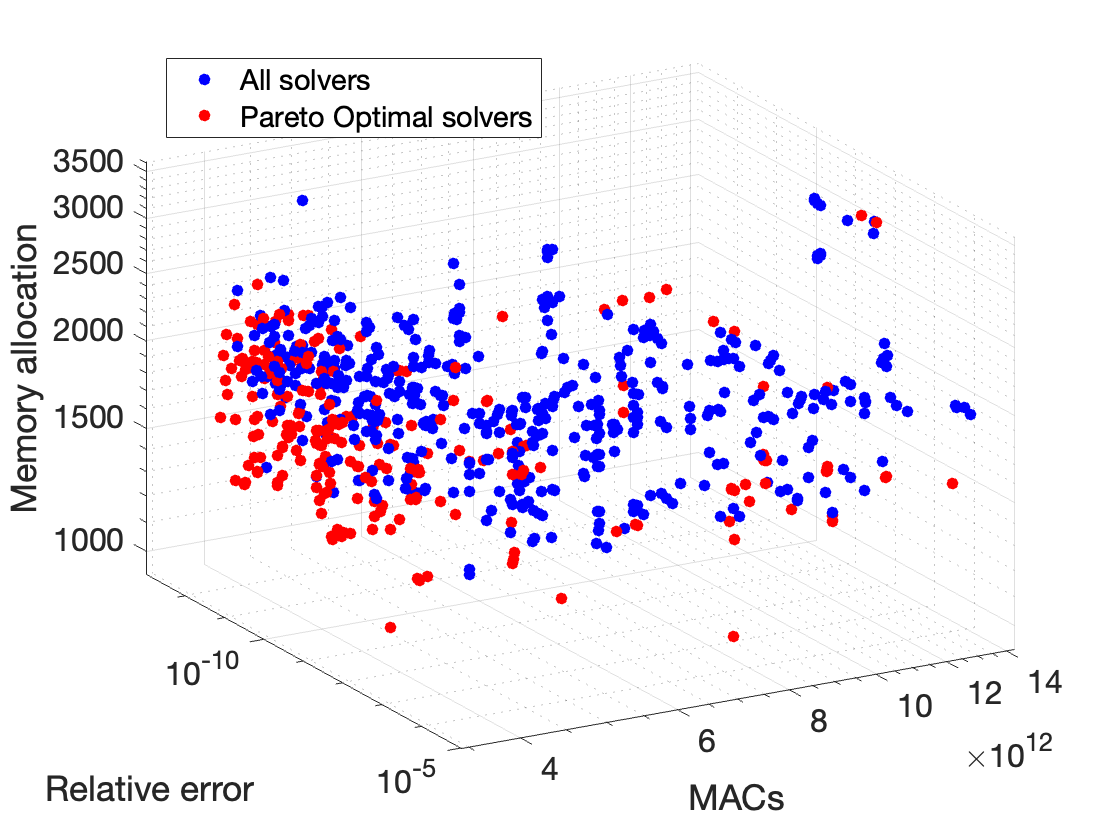}}
\subfigure[computational time -- relative error (p) -- memory allocation]{
\includegraphics[width=0.45\textwidth]{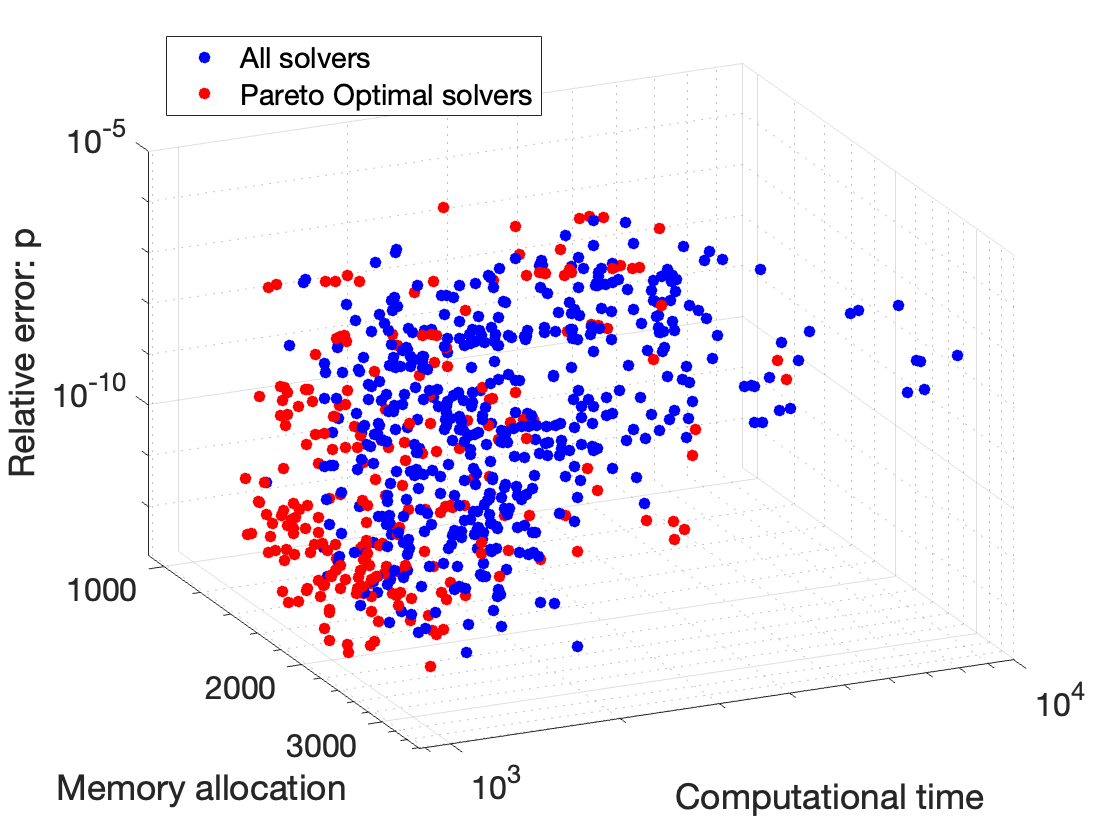}}
\subfigure[computational time -- MACs -- Average MACs]{
\includegraphics[width=0.45\textwidth]{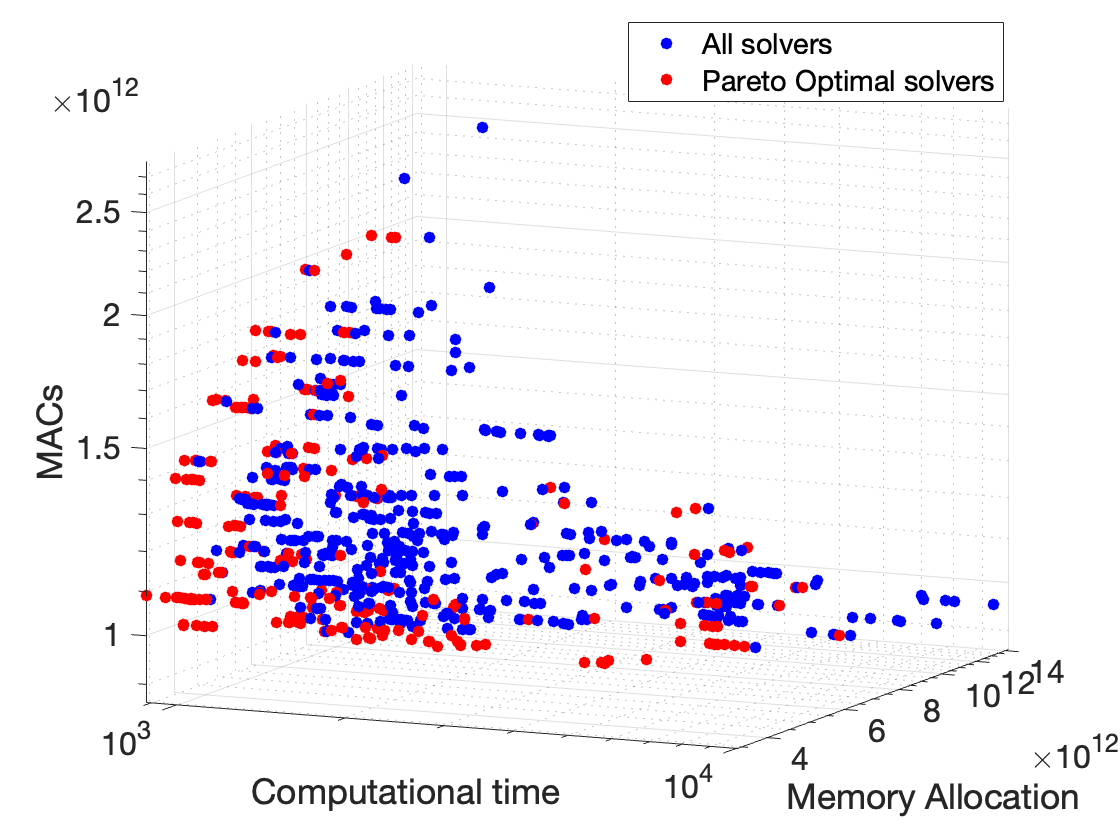}}
\caption{Pareto fronts for solving 3-D incompressive Navier-Stokes equation, using IMEX method, when $d = \frac{\pi}{8}$. All solvers are depicted in blue while Pareto optimal solvers are highlighted in red.}
\label{front_3dnsd2}
\end{figure}

\subsection{Discovery and re-discovery of optimal meta-solvers by preference functions}
We implement the preference function based discover methodology, for both cases. In particular, the first preference function is
\begin{pref}
$p^1(f') = \frac{1}{8}(\sum_{i=1}^8 f'_i)$. The average of all rescaled data.
\end{pref}

We present the results of top-3 solvers, for $d = \frac{\pi}{4}$:
\begin{table}[H]
    \centering
    \caption{Top-3 solvers evaluated by preference function $p^1$, for 3-D incompressible Navier-Stokes equation with  $d = \frac{\pi}{4}$, using IMEX based meta-solvers. }
     \begin{minipage}{\linewidth}
        \centering
        \captionof{subtable}{The top-3 solvers}
        \resizebox{\linewidth}{!}{
        \begin{tabular}{c|c|c|c|c|c}
\hline \hline 
 &Neural operator & Classical solver & Multi-grid & Relaxation & Strategies \\ 
\hline \hline
Top 1 solver & KAN & CG & 3-level & SSOR & 5-1-5 \\ \hline
Top 2 solver & DeepONet & CG & 3-level & Gauss-Seidel &  7-1-7 \\ \hline
Top 3 solver & U-DeepONet & CG & 3-level & SSOR & 3-1-3 \\ \hline
\end{tabular}
}
    \end{minipage}

    \vspace{1em} 

    \begin{minipage}{\linewidth}
        \centering
        \captionof{subtable}{Performance and rank of performance of the top-3 solvers}
        \resizebox{\linewidth}{!}{ 
          \begin{tabular}{c|c|c|c|c|c|c|c|c}
\hline \hline 
  &relative error:u& relative error:p& Com. time  & \# of ite & Memory & MACs&Ave. MACs & Training time \\ 
\hline 
Top 1 solver &$3.36 \times 10^{-11}$ & $2.28 \times 10^{-9}$ & 1095.7 & 6543 & 809.414 & $3.5248 \times 10^{12}$ & $1.11013 \times 10^{12}$ & 1086.76 \\
\hline
Top 2 solver & $9.55 \times 10^{-14}$ & $1.34 \times 10^{-11}$ & 1229.37 & 7959 & 784.054 & $3.74227 \times 10^{12}$ & $1.10191 \times 10^{12}$ & 26.4055 \\
\hline
Top 3 solver & $1.60 \times 10^{-13}$ & $5.33 \times 10^{-12}$ & 1317.67 & 8505 & 777.641 & $4.04756 \times 10^{12}$ & $1.06778 \times 10^{12}$ & 49.6242 \\
\hline
\end{tabular}
        }
    \end{minipage}
    \label{top3_p1_ns3dd1}
\end{table}

For $d = \frac{\pi}{8}$:

\begin{table}[H]
    \centering
    \caption{Top-3 solvers evaluated by preference function $p^1$, for 3-D incompressible Navier-Stokes equation with  $d = \frac{\pi}{8}$, using IMEX based meta-solvers.}
     \begin{minipage}{\linewidth}
        \centering
        \captionof{subtable}{The top-3 solvers}
        \resizebox{\linewidth}{!}{
        \begin{tabular}{c|c|c|c|c|c}
\hline \hline 
 &Neural operator & Classical solver & Multi-grid & Relaxation & Strategies \\ 
\hline \hline
Top 1 solver & JacobiKAN & CG & 2-level & Gauss-Seidel & 5-1-5 \\ \hline
Top 2 solver & DeepONet & CG & 2-level & Gauss-Seidel &  7-1-7 \\ \hline
Top 3 solver & U-DeepONet & CG & 2-level & SSOR & 5-1-5 \\ \hline
\end{tabular}
}
    \end{minipage}

    \vspace{1em} 

    \begin{minipage}{\linewidth}
        \centering
        \captionof{subtable}{Performance and rank of performance of the top-3 solvers}
        \resizebox{\linewidth}{!}{ 
          \begin{tabular}{c|c|c|c|c|c|c|c|c}
\hline \hline 
  &relative error:u& relative error:p& Com. time  & \# of ite & Memory & MACs&Ave. MACs & Training time \\ 
\hline 
Top 1 solver &$9.36 \times 10^{-11}$ & $1.72 \times 10^{-9}$ & 991.732 & 8519 & 1403.84 & $3.61028 \times 10^{12}$ & $1.02183 \times 10^{12}$ & 115.963 \\
\hline
Top 2 solver & $6.57 \times 10^{-14}$ & $8.21 \times 10^{-12}$ & 1016.13 & 7730 & 1408.26 & $3.65938 \times 10^{12}$ & $1.08166 \times 10^{12}$ & 26.4055 \\
\hline
Top 3 solver & $3.34 \times 10^{-11}$ & $2.60 \times 10^{-9}$ & 1008.1 & 6487 & 1546.77 & $3.46252 \times 10^{12}$ & $1.08682 \times 10^{12}$ & 49.6242 \\
\hline
\end{tabular}
        }
    \end{minipage}
    \label{top3_p1_ns3dd2}
\end{table}

For the second preference function, we consider 
\begin{pref}
$p^2(f') = \frac{1}{25}* \Big( 5* (f'_1+f'_2)+ 
10*f'_3+(\sum_{i=4}^8 f'_i) \Big)$. 
\end{pref}

The top-3 solvers for $d = \frac{\pi}{4}$:
\begin{table}[H]
    \centering
    \caption{Top-3 solvers evaluated by preference function $p^2$, for 3-D incompressible Navier-Stokes equation with  $d = \frac{\pi}{4}$, using IMEX based meta-solvers.  }
     \begin{minipage}{\linewidth}
        \centering
        \captionof{subtable}{The top-3 solvers}
        \resizebox{\linewidth}{!}{
        \begin{tabular}{c|c|c|c|c|c}
\hline \hline 
 &Neural operator & Classical solver & Multi-grid & Relaxation & Strategies \\ 
\hline \hline
Top 1 solver & KAN & CG & 3-level & SSOR & 5-1-5 \\ \hline
Top 2 solver & KAN & CG & 2-level & Gauss-Seidel &  7-1-7 \\ \hline
Top 3 solver & DeepONet & CG & 3-level & Gauss-Seidel & 7-1-7 \\ \hline
\end{tabular}
}
    \end{minipage}

    \vspace{1em} 

    \begin{minipage}{\linewidth}
        \centering
        \captionof{subtable}{Performance and rank of performance of the top-3 solvers}
        \resizebox{\linewidth}{!}{ 
          \begin{tabular}{c|c|c|c|c|c|c|c|c}
\hline \hline 
  &relative error:u& relative error:p& Com. time  & \# of ite & Memory & MACs&Ave. MACs & Training time \\ 
\hline 
Top 1 solver &$3.36 \times 10^{-11}$ & $2.28 \times 10^{-9}$ & 1095.7 & 6543 & 809.414 & $3.5248 \times 10^{12}$ & $1.11013 \times 10^{12}$ & 1086.76 \\
\hline
Top 2 solver & $8.03 \times 10^{-14}$ & $1.41 \times 10^{-11}$ & 1008.61 & 7580 & 1652.64 & $3.6277 \times 10^{12}$ & $1.07766 \times 10^{12}$ & 1086.76 \\
\hline
Top 3 solver & $9.55 \times 10^{-14}$ & $1.34 \times 10^{-11}$ & 1229.37 & 7959 & 784.054 & $3.74227 \times 10^{12}$ & $1.10191 \times 10^{12}$ & 26.4055 \\
\hline
\end{tabular}
        }
    \end{minipage}
    \label{top3_p2_ns3dd1}
\end{table}

For $d = \frac{\pi}{8}$:
\begin{table}[H]
    \centering
    \caption{Top-3 solvers evaluated by preference function $p^2$, for 3-D incompressible Navier-Stokes equation with  $d = \frac{\pi}{8}$, using IMEX based meta-solvers.  }
     \begin{minipage}{\linewidth}
        \centering
        \captionof{subtable}{The top-3 solvers}
        \resizebox{\linewidth}{!}{
        \begin{tabular}{c|c|c|c|c|c}
\hline \hline 
 &Neural operator & Classical solver & Multi-grid & Relaxation & Strategies \\ 
\hline \hline
Top 1 solver & DeepONet & CG & 2-level & SSOR & 5-1-5 \\ \hline
Top 2 solver & JacobiKAN & CG & 2-level & Gauss-Seidel &  5-1-5 \\ \hline
Top 3 solver & DeepONet & CG & 2-level & Gauss-Seidel & 7-1-7 \\ \hline
\end{tabular}
}
    \end{minipage}

    \vspace{1em} 

    \begin{minipage}{\linewidth}
        \centering
        \captionof{subtable}{Performance and rank of performance of the top-3 solvers}
        \resizebox{\linewidth}{!}{ 
          \begin{tabular}{c|c|c|c|c|c|c|c|c}
\hline \hline 
  &relative error:u& relative error:p& Com. time  & \# of ite & Memory & MACs&Ave. MACs & Training time \\ 
\hline 
Top 1 solver &$3.34 \times 10^{-11}$ & $2.60 \times 10^{-9}$ & 883.912 & 6249 & 2042.37 & $3.43909 \times 10^{12}$ & $1.08671 \times 10^{12}$ & 26.4055 \\
\hline
Top 2 solver & $9.36 \times 10^{-11}$ & $1.72 \times 10^{-9}$ & 991.732 & 8519 & 1403.84 & $3.61028 \times 10^{12}$ & $1.02183 \times 10^{12}$ & 115.963 \\
\hline
Top 3 solver & $6.57 \times 10^{-14}$ & $8.21 \times 10^{-12}$ & 1016.13 & 7730 & 1408.26 & $3.65938 \times 10^{12}$ & $1.08166 \times 10^{12}$ & 26.4055 \\
\hline
\end{tabular}
        }
    \end{minipage}
    \label{top3_p2_ns3dd2}
\end{table}

\begin{figure}[H]
\centering
\subfigure[Discovery of optimal solvers by $p^1$]{
\includegraphics[width=0.45\textwidth]{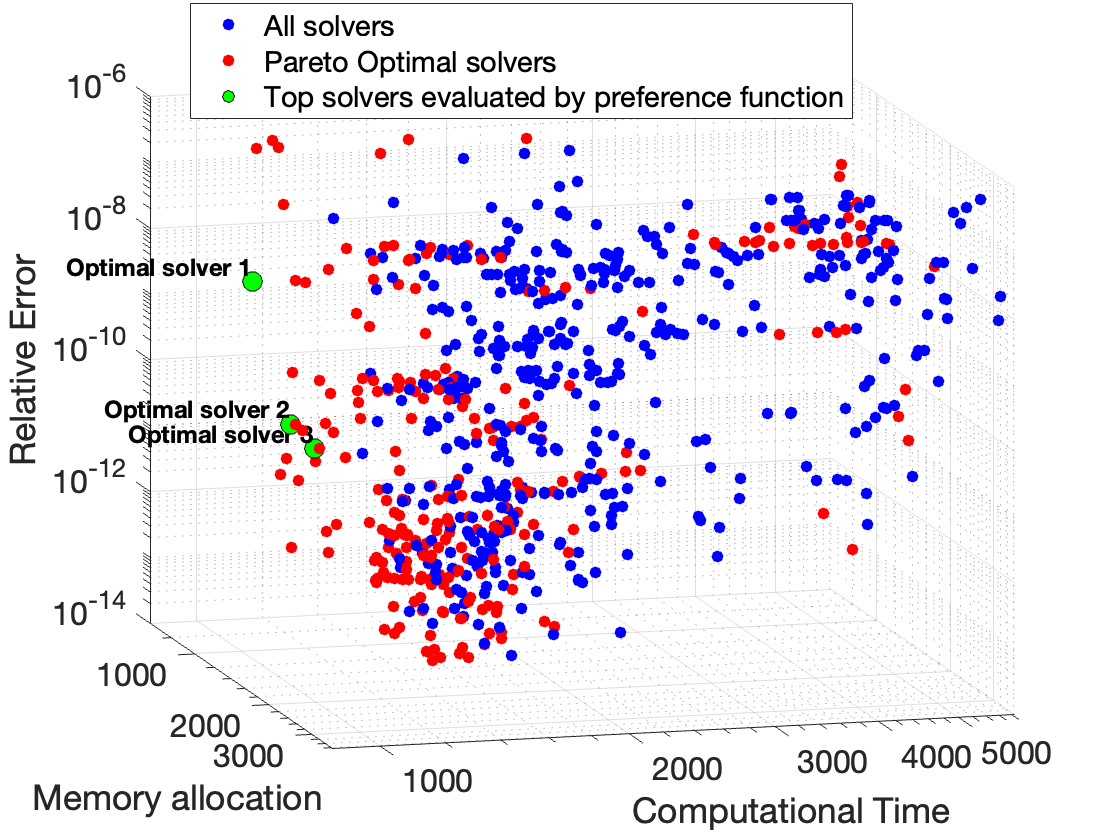}}
\subfigure[Discovery of optimal solvers by $p^2$]{
\includegraphics[width=0.45\textwidth]{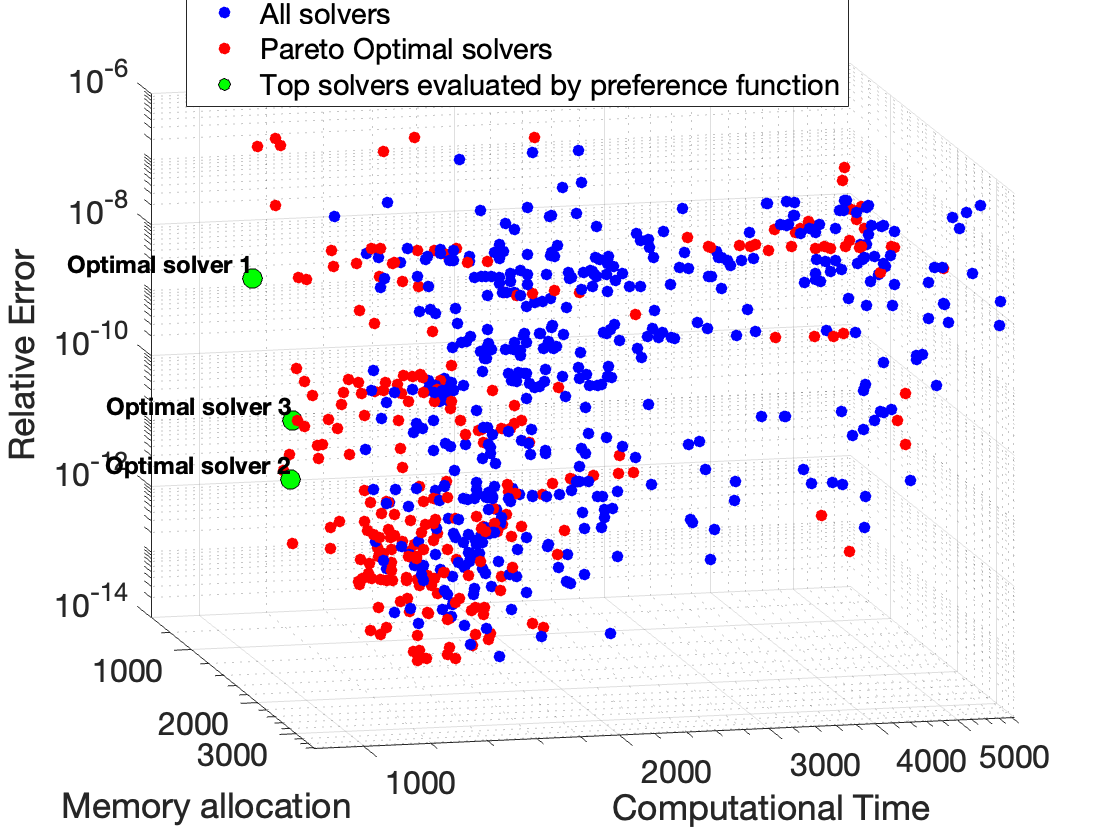}}
\caption{Pareto front: Computational Time -- Relative Error -- Memory allocation, with top 3 solvers evaluated by $p^1$ and $p^2$, for solving 3d Incompressible Navier Stokes with $d=\frac{\pi}{4}$ equation using IMEX based meta-solvers.}
\label{discover_ns3d_1}
\end{figure}

\begin{figure}[H]
\centering
\subfigure[Discovery of optimal solvers by $p^1$]{
\includegraphics[width=0.45\textwidth]{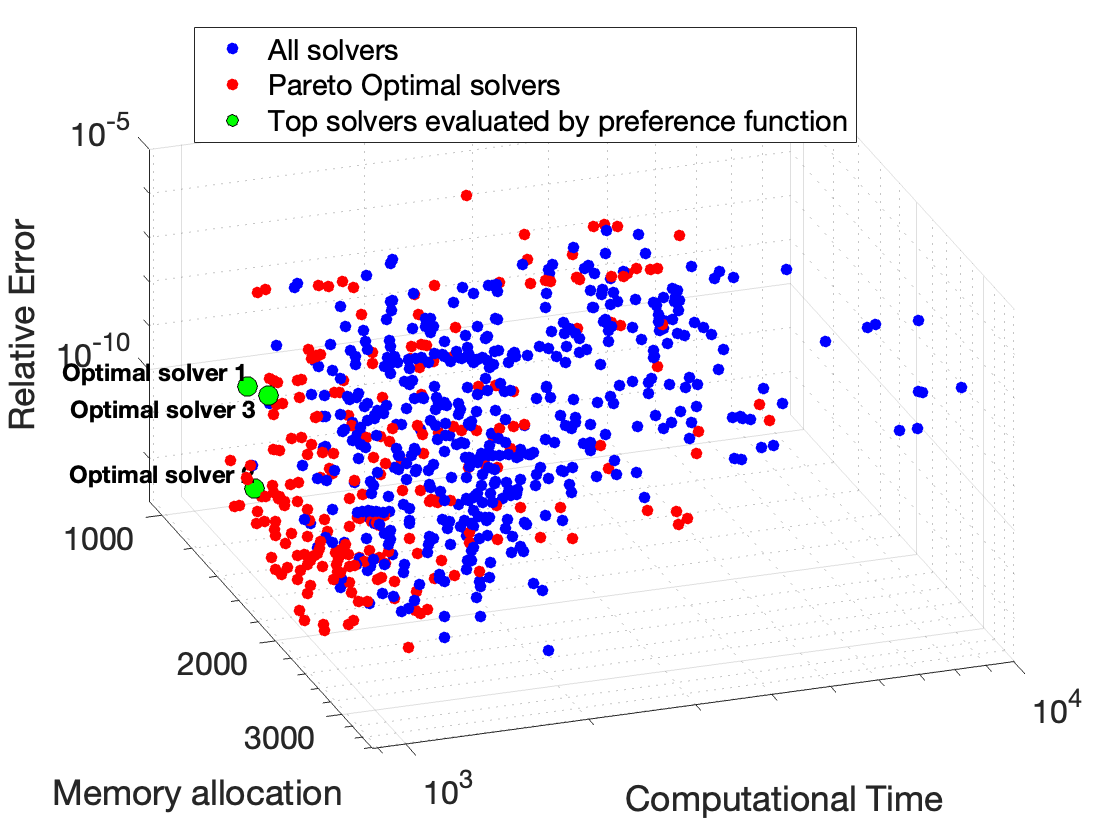}}
\subfigure[Discovery of optimal solvers by $p^2$]{
\includegraphics[width=0.45\textwidth]{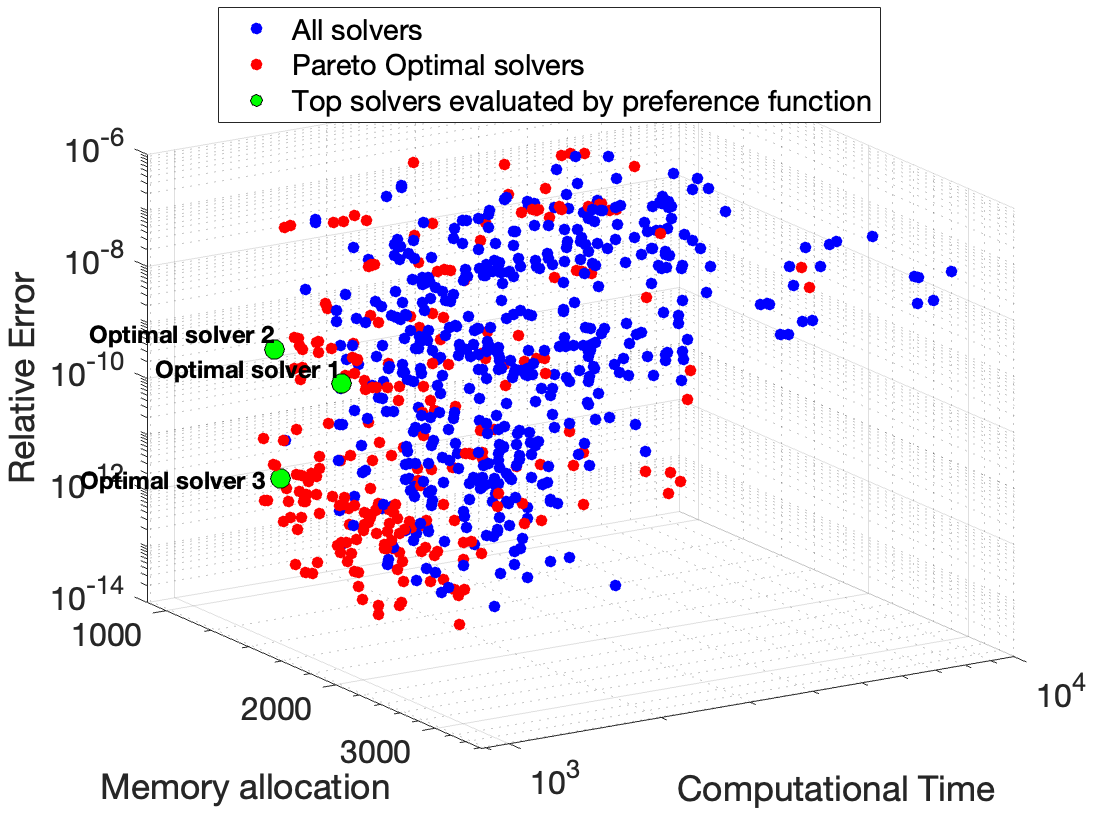}}
\caption{Pareto front: Computational Time -- Relative Error -- Memory allocation, with top 3 solvers evaluated by $p^1$ and $p^2$, for solving 3d Incompressible Navier Stokes with $d=\frac{\pi}{8}$ equation using IMEX based meta-solvers.}
\label{discover_ns3d_2}
\end{figure}

We plot the performance of the optimal meta-solvers discovered by different preference functions, in~\Cref{discover_ns3d_1} for case $d=\frac{\pi}{4}$ and in~\Cref{discover_ns3d_2} for case $d=\frac{\pi}{8}$

Moreover, we present the results of using the Linear Programming approach to re-discover the particular, weighted sum, preference functions, for given Pareto meta-solvers. 
In the following we present the weights in the order s.t. $(\lambda_1,\lambda_2,\dots,\lambda_8)$ for ( $\lambda_1$: relative error of $\textbf{u}$, $\lambda_2$: relative error of $\textbf{p}$,  $\lambda_3$: computational time, $\lambda_4$: number of iterations, 
$\lambda_5$: memory allocation, $\lambda_6$: MACs, $\lambda_7$: average MACs, $\lambda_8$: Training time), that is for a given Pareto optimal meta-solver.

For $d=\frac{\pi}{4}$, we have the following rediscovery results:

\begin{subequations}
\begin{equation}
\left\{
\begin{aligned}
&\text{Preference function $p(\lambda_a;r) = (0,0.0394,0.1674,0,0.0178,0,0.7755,0)^Tr$  } \\
&\text{Optimal solver}: x^a = (\text{KAN}, \text{CG}, \text{1-level}, \text{SSOR}, 1-1-1) \ .
\end{aligned}
\right.
\end{equation}

\begin{equation}
\left\{
\begin{aligned}
&\text{Preference function $p(\lambda_a;r) = (0,0,0,0.2656,0.0344,0.2228,0.4772,0)^Tr$ } \\
&\text{Optimal solver}: x^a = (\text{U-DeepONet}, \text{CG}, \text{3-level}, \text{Gauss-Seidel}, 5-1-5) \ .
\end{aligned}
\right.
\end{equation}

\begin{equation}
\left\{
\begin{aligned}
&\text{Preference function $p(\lambda_a;r) = (0.0194,0,0.0318,0.3270,0.5327,0,0.0841,0)^Tr$ } \\
&\text{Optimal solver}: x^a = (\text{CheyKAN}, \text{BiCGStab}, \text{3-level}, \text{Gauss-Seidel}, 5-1-5) \ .
\end{aligned}
\right.
\end{equation}

\end{subequations}

For $d=\frac{\pi}{8}$, we have the following rediscovery results:

\begin{subequations}
\begin{equation}
\left\{
\begin{aligned}
&\text{Preference function $p(\lambda_a;r) = (0,0,0,0.5510,0.2619,0,0.1871,0)^Tr$  } \\
&\text{Optimal solver}: x^a = (\text{JacobiKAN}, \text{CG}, \text{3-level}, \text{SSOR}, 5-1-5) \ .
\end{aligned}
\right.
\end{equation}

\begin{equation}
\left\{
\begin{aligned}
&\text{Preference function $p(\lambda_a;r) = (0.5109,0,0,0.0228,0.0066,0.1049,0.3536,0.0011)^Tr$ } \\
&\text{Optimal solver}: x^a = (\text{JacobiKAN}, \text{CG}, \text{2-level}, \text{Gauss-Seidel}, 3-1-3) \ .
\end{aligned}
\right.
\end{equation}

\begin{equation}
\left\{
\begin{aligned}
&\text{Preference function $p(\lambda_a;r) = (0,0.5158,0,0.4347,0.0014,0,0.0478,0.0003)^Tr$ } \\
&\text{Optimal solver}: x^a = (\text{JacobiKAN}, \text{CG}, \text{2-level}, \text{SSOR}, 7-1-7) \ .
\end{aligned}
\right.
\end{equation}

\end{subequations}

\section{Implementation details and further numerical results for solving brittle fracture problem~(phase-field modeling)}\label{app-bf} 

\subsection{The model and implementation details}
Assume a domain $\Omega \subset \R^{2}$ is given.
The external boundary $\Gamma$ is decomposed into the Dirichlet boundary $\Gamma_{D}$ and Neumann boundary $\Gamma_{N}$, i.e., $\Gamma = \Gamma_{D} \cup \Gamma_{N}$.
Moreover, the Dirichlet boundary $\Gamma_{D}$ consists of a homogeneous boundary $\Gamma_{D,0}$ and non-homogeneous~(loading) boundary $\Gamma_{D,1}$. 
Let $\Omega_{c} \subset \Omega$ be a crack set.
Then, we consider the $n$-th loading step of the phase-field modeling of Brittle fracture.
The pair of solution spaces $\{V,Q\}$ is given by
\begin{align}
    V &= \{\mathbf{u} \in [H^{1}(\Omega)]^{d}\colon \mathbf{u}=0 \text{ on } \Gamma_{D,0}, \mathbf{u}=\mathbf{u}_{n} \text{ on } \Gamma_{D,1}\}, \\
    Q &= \{\alpha \in H^{1}(\Omega)\colon 0 \leq \alpha \leq 1, \alpha \geq \alpha_{n-1} \text{ in } \Omega\},
\end{align}
Then, the pair of solution $\{\mathbf{u}, \alpha\}$ is obtained by solving the following minimization problem:
\begin{equation}
\min_{\mathbf{u} \in V, \alpha \in Q} \mathcal{E}_{n}(\mathbf{u}, \alpha) := \mathcal{E}_{el} + \mathcal{E}_{d} - \mathcal{E}_{w} + \mathcal{E}_{irr},
\label{app:eqn:phase-field}
\end{equation}
where
\begin{equation}
\left\{\begin{split}
\mathcal{E}_{el}&=\int_{\Omega} \left(a(\alpha)\Psi^{+}(\varepsilon(\mathbf{u}))+\Psi^{-}(\varepsilon(\mathbf{u})\right)dx, \\
\mathcal{E}_{d}&=\frac{G_{c}}{c_{w}}\int_{\Omega} (\frac{w(\alpha)}{\ell}+\ell \vert \nabla \alpha Bro\vert^{2})dx, \\
\mathcal{E}_{w} &= \int_{\Omega} \mathbf{f}_{n} \cdot \mathbf{u} dx + \int_{\Gamma_{N}}\mathbf{t}_{n} \cdot \mathbf{u} ds, \\
\mathcal{E}_{irr} &= \frac{\gamma}{2}\int_{\Omega} \langle \alpha - \alpha_{n-1} \rangle_{-}^{2} dx.
\end{split}
\right.
\label{eqn:phase-field-def}
\end{equation}
Note that 
\begin{equation}
    \left\{ \begin{split}
        \varepsilon(\mathbf{u}) &= \frac{1}{2} \left[ \nabla \mathbf{u} + (\nabla \mathbf{u})^{T} \right], \\
        \lambda &= \frac{E \nu}{(1+\nu)(1-2\nu)}, \\
        \mu &= \frac{E}{2(1+\nu)}, \\
        \Psi^{+}(\varepsilon) &= \frac{1}{2} \kappa \langle \tr(\varepsilon)\rangle_{+}^{2} + 2\mu \mathbf{e}(\varepsilon) \cdot \mathbf{e}(\varepsilon), \\
        \kappa &= \lambda + \frac{2}{d}\mu,\\
        \mathbf{e}(\varepsilon) &= \varepsilon - \frac{\tr(\varepsilon)}{d}I, \\
        \Psi^{-}(\varepsilon) &= \frac{1}{2} \kappa \langle \tr(\varepsilon)\rangle_{-}^{2}, \\
        a(\alpha) &= (1 - k)(1 - \alpha)^{2} +k, \\
        w(\alpha) &= \begin{cases}
            \alpha, &\text{ AT1 model},\\
            \alpha^{2}, &\text{ AT2 model},\\
        \end{cases} \\
        c_{w} &= \begin{cases}
            \frac{8}{3}, &\text{ AT-1 model},\\
            2, &\text{ AT-2 model},\\
        \end{cases} \\
        \gamma &= \begin{cases}
            \frac{G_{c}}{\ell}\left(\frac{27}{64\tau^{2}}\right), &\text{ AT-1 model},\\
            \frac{G_{c}}{\ell}\left(\frac{1}{\tau^{2}}-1\right), &\text{ AT-2 model}.\\
        \end{cases} \\
    \end{split} \right.
\label{eqn:phase-field-parameter}
\end{equation}

In this experiment, we use the AT-2 model, $E=210.0$, $\nu =0.3$, $\ell = 0.01$, $G_{c}=2.7 \times 10^{-3}$, $k=10^{-8}$, $\tau=10^{-2}$, and $\gamma=2.7 \times 10^{3}$.

We also use the alternating minimization technique.
Let $\alpha^{(0)} = 0$.
For $n = 1, 2, \ldots$,
\begin{enumerate}
    \item Find $\mathbf{u}^{(n)} \in V$ satisfying
\begin{equation*}
    \frac{\partial \mathcal{E}_{n}}{\partial \mathbf{u}}(\mathbf{u}^{(n)}, \alpha^{(n-1)}) = 0,
\end{equation*}
where the symbol $\alpha^{(n-1)}$ denotes the previous damage.
\item After that, find $\alpha^{(n)} \in Q$ satisfying
\begin{equation*}
    \frac{\partial \mathcal{E}_{n}}{\partial \alpha}(\mathbf{u}^{(n)}, \alpha^{(n)}) = 0.
\end{equation*}
\item If $\Vert \alpha^{(n)} - \alpha^{(n-1)}\Vert_{\infty} < \varepsilon$, then $\mathbf{u}=\mathbf{u}^{(n)}$ and $\alpha=\alpha^{(n)}$. If not, $n \gets n+1$.
\end{enumerate}

For all local problems, we utilize the Newton-Raphson method. The implementation details are given in the following:

\subsubsection*{Common details}
\begin{itemize}
    \item The domain $\Omega$ is descretized by $6,252$ nodes with mesh size varying from $h_{min}=\frac{1}{500}$ to $h_{max}=\frac{1}{25}$.
    \item $\mathbf{f}_{n}=\mathbf{t}_{n}=0$.
    \item $\varepsilon=10^{-4}$ (Stop criteria for alternating minimization)
    \item Maximum number of iterations for the alternating minimization: $1000$.
    \item Maximum number of iterations for the Newton-Rapshon method: $1000$.
    \item Maximum number of iterations for the inner linearized problem: $2000$.
\end{itemize}
\subsubsection*{Training dataset}
\begin{itemize}
    \item We sampled $7$ pairs.
    \item Spatial mesh size for network: $\Delta x = \frac{1}{30}$
    \item Time step~(Quasi-static): $\Delta t = \frac{1}{50}$.
    \item $u_{T} \sim \mathcal{U}(0.001, 0.01)$, $\mathbf{u}_{n}=[0.0, u_{T}*n*\Delta t]$
    \item Total number of training/validation samples: $335/15$.
    \item Polynomial order for $\mathbf{u}$ and $\alpha$ is $1$. 
\end{itemize}

\subsubsection*{Test case 1}
\begin{itemize}
    \item The domain $\Omega$ is discretized by $6,252$ nodes with mesh size varying from $h_{min}=\frac{1}{500}$ to $h_{max}=\frac{1}{25}$.
    \item $u_{T}=0.006$, $\Delta t = \frac{1}{30}$ $\Rightarrow$ $\delta u_{T}= \frac{u_{T}}{\Delta T}=0.0002$.
\end{itemize}

\subsubsection*{Test case 2~(much more fine grid)}
\begin{itemize}
    \item The domain $\Omega$ is discretized by $37,949$ nodes with mesh size varying from $h_{min}=\frac{3}{2000}$ to $h_{max}=\frac{3}{100}$.
    \item $u_{T}=0.006$, $\Delta t = \frac{1}{30}$ $\Rightarrow$ $\delta u_{T}= \frac{u_{T}}{\Delta T}=0.0002$.
\end{itemize}

\subsubsection*{Modified approach}
\begin{itemize}
    \item Instead of using the penalty $\mathcal{E}_{irr}=\frac{\gamma}{2}\int_{\Omega}\langle \alpha - \alpha_{n-1} \rangle_{-}^{2}dx$, we use simple projection: we apply
    \begin{equation*}
        \displaystyle (\alpha_{n})_{j} \gets \max((\alpha_{n})_{j}, (\alpha_{n-1})_{j})
    \end{equation*}
    \item $u_{T}=0.007$, $\Delta t = \frac{1}{35}$ $\Rightarrow$ $\delta u_{T}= \frac{u_{T}}{\Delta T}=0.0002$.
\end{itemize}

\begin{figure}
    \centering
    \includegraphics[width=\linewidth, trim={5cm, 0cm, 5cm, 0cm}, clip]{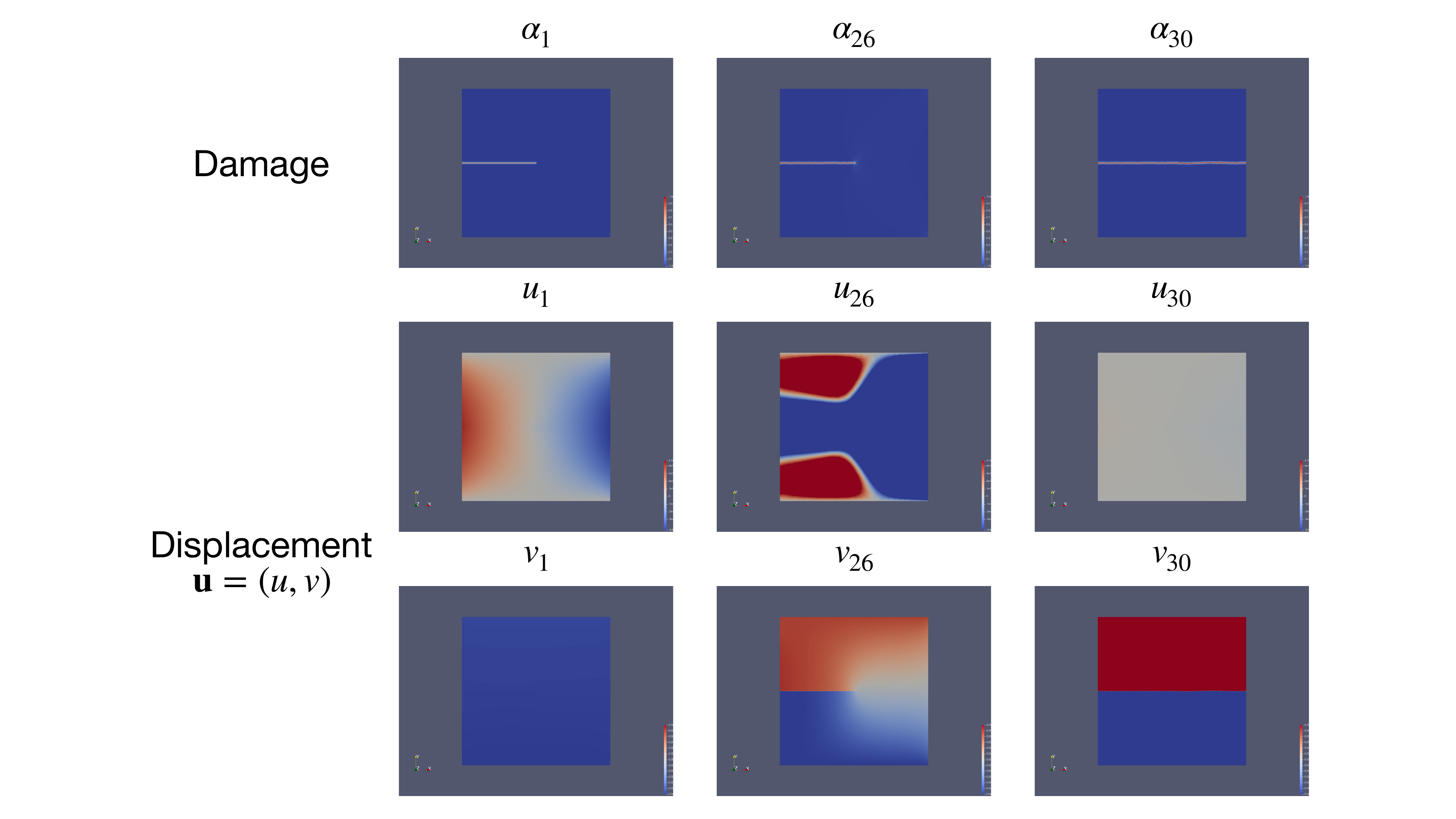}
    \caption{The snapshots of solution for SEN specimen when $u_{T}=0.006$ and $n=1,26,30$.}
    \label{fig:sen_sol}
\end{figure}

Here, we also briefly describe the Chebyshev semi-iterative method which can be used for smoother. Given a linear system $\mathbf{A}\mathbf{x} = \mathbf{b}$, the classical preconditioned Richardson iteration is given by
\begin{align}
\begin{cases}
    \mathbf{r}^{(i)} &= \mathbf{b} - \mathbf{A}\mathbf{x}^{(i)},\\
    \mathbf{x}^{(i+1)} &= \mathbf{x}^{(i)} + \mathbf{M}^{-1}(\mathbf{r}^{(i)}),
\end{cases}
\end{align}
where $\mathbf{M}^{-1}$ denotes a preconditioner.
Note that the classical relaxation methods can be represented by 
\begin{align}
    \mathbf{M} = \begin{cases}
        \mathbf{D}, &\text{Jacobi,}\\
        \mathbf{D}+\mathbf{L}, &\text{Gauss-Seidel,}\\
        \frac{1}{\omega} \mathbf{D}+\mathbf{L}, &\text{SOR,}\\
        \frac{1}{\omega (2- \omega)} (\mathbf{D} + \omega \mathbf{L})\mathbf{D}^{-1}(\mathbf{D} + \omega \mathbf{U}), &\text{SSOR,}\\
    \end{cases}
\end{align}
where $\mathbf{A} = \mathbf{D} + \mathbf{L} + \mathbf{U}$.

On the other hand, let us assume that $\lambda_{min}, \lambda_{max}$ are the smallest and largest eigenvalues of $A$, respectively.
In addition, a symbol $T_{n}(\xi)$ denotes the $n$-th order Chebyshev polynomial defined in the interval $[-1,1]$.
The translated and scaled residual polynomial $p_{n}(\xi)$ is defined as
\begin{equation}
    p_{n}(\xi) := \frac{T_{n}(\frac{\xi - a}{c})}{T_{n}(-\frac{a}{c})},
\end{equation}
where $a = \frac{\lambda_{max} + \lambda_{min}}{2}$ and $c = \frac{\lambda_{max} - \lambda_{min}}{2}$.
Then, the $i$-th approximation and residual of the Chebyshev semi-iterative method satisfy
\begin{equation}
    \mathbf{r}^{(i)} = \mathbf{b} - \mathbf{A}\mathbf{x}^{(i)} = p_{i}(\mathbf{A})\mathbf{r}^{(0)}.
\end{equation}

\begin{algorithm}
\caption{Chebyshev semi-iterative method}\label{alg:cheb}
\begin{algorithmic}
\Require $\lambda_{min}$, $\lambda_{max}$, $\mathbf{x}^{(0)}$, 
$\mathbf{r}^{(0)} = \mathbf{b}-\mathbf{A}\mathbf{x}^{(0)}$,
$\mathbf{x}^{(-1)}=\mathbf{r}^{(-1)}=\mathbf{0}$
\State $a \gets \frac{\lambda_{max}+\lambda_{min}}{2}$, $c \gets \frac{\lambda_{max}-\lambda_{min}}{2}$, $\eta \gets -\frac{a}{c}$
\State $\beta_{-1} \gets 0$, $\beta_{0} \gets \frac{c}{2}\frac{1}{\eta}=-\frac{c^{2}}{a}$, $\gamma_{0} \gets -a$
\State $i \gets 0$
\While{until convergence}
\State $\mathbf{x}^{(i+1)} \gets -(\mathbf{r}^{(i)} +a\mathbf{x}^{(i)}+\beta_{i-1}\mathbf{x}^{(i-1)})/\gamma_{i}$
\State $\mathbf{r}^{(i+1)} \gets -(\mathbf{A}\mathbf{r}^{(i)} -a\mathbf{r}^{(i)}-\beta_{i-1}\mathbf{r}^{(i-1)})/\gamma_{i}$
\State $i \gets i + 1$
\State $\beta_{i-1} \gets \frac{c}{2}\frac{T_{i-1}(\eta)}{T_{i}(\eta)}=\left(\frac{c}{2}\right)^{2}\frac{1}{\gamma_{i-1}}$ if $i \geq 2$
\State $\gamma_{i} \gets \frac{c}{2}\frac{T_{i+1}(\eta)}{T_{i}(\eta)} = -(a+\beta_{i-1})$
\EndWhile
\end{algorithmic}
\end{algorithm}

\subsection{Identify all optimal meta-solvers in Pareto sense}

For the coarse mesh-step, there are 144 Pareto optimal solvers. For the fine mesh-step, there are 54 Pareto optimal solvers. In the following, we first summarize the composition of the set of Pareto optimal solvers, by counting the number of different components in the construction of meta-solvers. 

\begin{table}[H]
    \centering
    \caption{The composition of the set of Pareto optimal solvers by counting the number of elements in each dimension, for 2d brittle fracture problem.}
    \begin{minipage}{\linewidth}
        \centering
        \captionof{subtable}{Different neural operators.}
        \resizebox{\linewidth}{!}{
        \begin{tabular}{c| c c c  c c |c} 
\hline
Neural Op & DeepONet & U-Net  & KAN & JacobiKAN& ChebyKAN & Total\\
\hline
\# in Pareto opt for coarse mesh & 38 & 36 & 24 & 22 & 24 & 144 \\
\hline
\# in Pareto opt for fine mesh & 18 & 17 & 7 & 6 & 6 & 54 \\
\hline

\end{tabular}}

    \end{minipage}

    \vspace{1em} 

    \begin{minipage}{\linewidth}
        \centering
        \captionof{subtable}{Different Krylov solvers.}
           \begin{tabular}{c| c c c c    } 
\hline
Classical solvers & GMRES & CG & BiCGStab \\
\hline
\# in Pareto opt for coarse mesh  & 25 & 48 & 71  \\
\hline
\# in Pareto opt for fine mesh  & 4 & 19 & 31  \\
\hline
\end{tabular}
        
    \end{minipage}
    
    \vspace{1em}

    \begin{minipage}{\linewidth}
        \centering
        \captionof{subtable}{Different smoothers.}
           \begin{tabular}{c| c c c c   } 
\hline
Smoother & GS  & SOR & SSOR &Cheyshev  \\
\hline
\# in Pareto opt for coarse mesh  & 17 & 19 & 13 & 95  \\
\hline
\# in Pareto opt for fine mesh & 11 & 0 & 2 & 41  \\
\hline
\end{tabular}
        
    \end{minipage}

\vspace{1em}

\begin{minipage}{\linewidth}
        \centering
        \captionof{subtable}{Different strategies of applying smoothers.}
          \begin{tabular}{c| c c c c c    } 
\hline
Strategies for smoother & 1-1-1 & 3-1-3 & 5-1-5 & 7-1-7 & 9-1-9  \\
\hline
\# in Pareto opt for coarse mesh & 22 & 32 & 33 & 28 & 29 \\
\hline
\# in Pareto opt & 13 & 15 & 11 & 7 & 8 \\
\hline

\end{tabular}
        
    \end{minipage}
    
\vspace{1em}

\begin{minipage}{\linewidth}
        \centering
        \captionof{subtable}{Adaption of algebraic multigrid method.}
   \begin{tabular}{c| c c     } 
\hline
Adaption of AMG & Yes & No    \\
\hline
\# in Pareto opt for coarse mesh & 104 & 40   \\
\hline
\# in Pareto opt for fine mesh & 48 & 6  \\
\hline

\end{tabular}

       \end{minipage}
    \label{bf_2d_com}
\end{table}

We plot the Pareto fronts, for MACs, average MACs, for computational time, relative error (p), number of iterations, for relative error (p), memory allocation, MACs and for computational time, relative error (p), memory allocation

\begin{figure}[H]
\centering
\subfigure[Computational time -- relative error -- number of iterations]{
\includegraphics[width=0.45\textwidth]{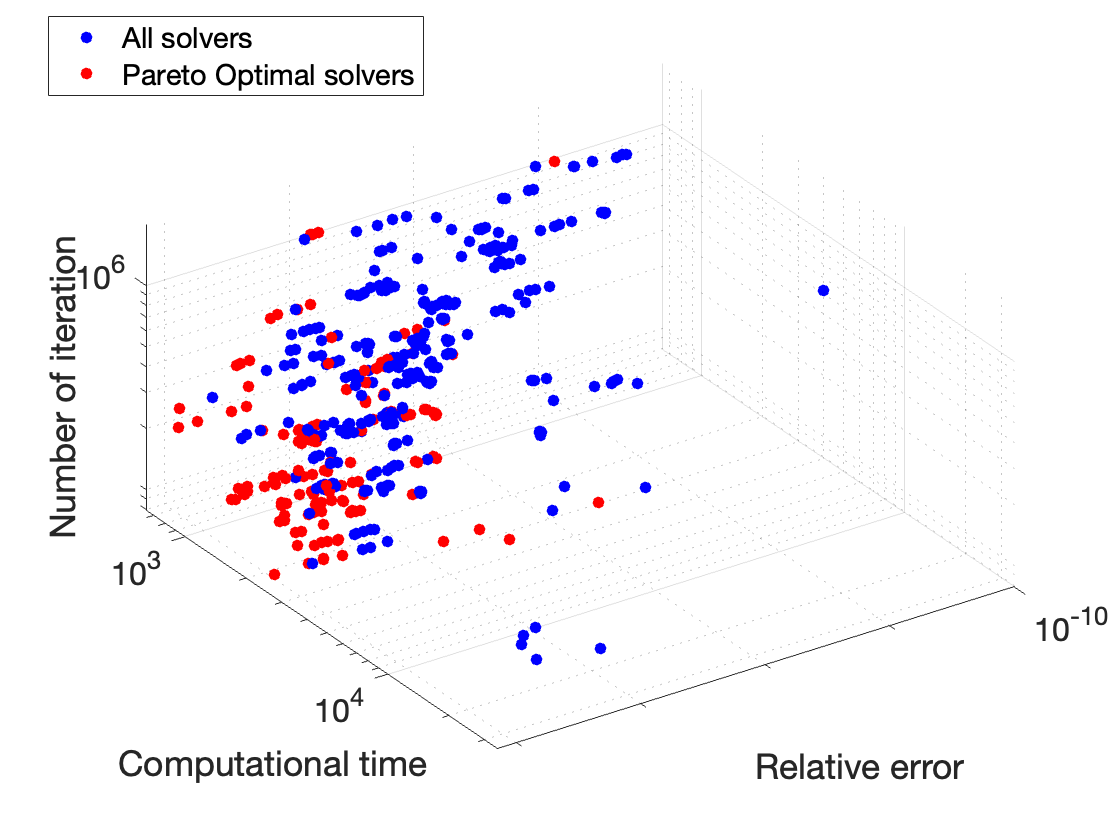}}
\subfigure[Relative error (p) -- memory allocation -- MACs]{
\includegraphics[width=0.45\textwidth]{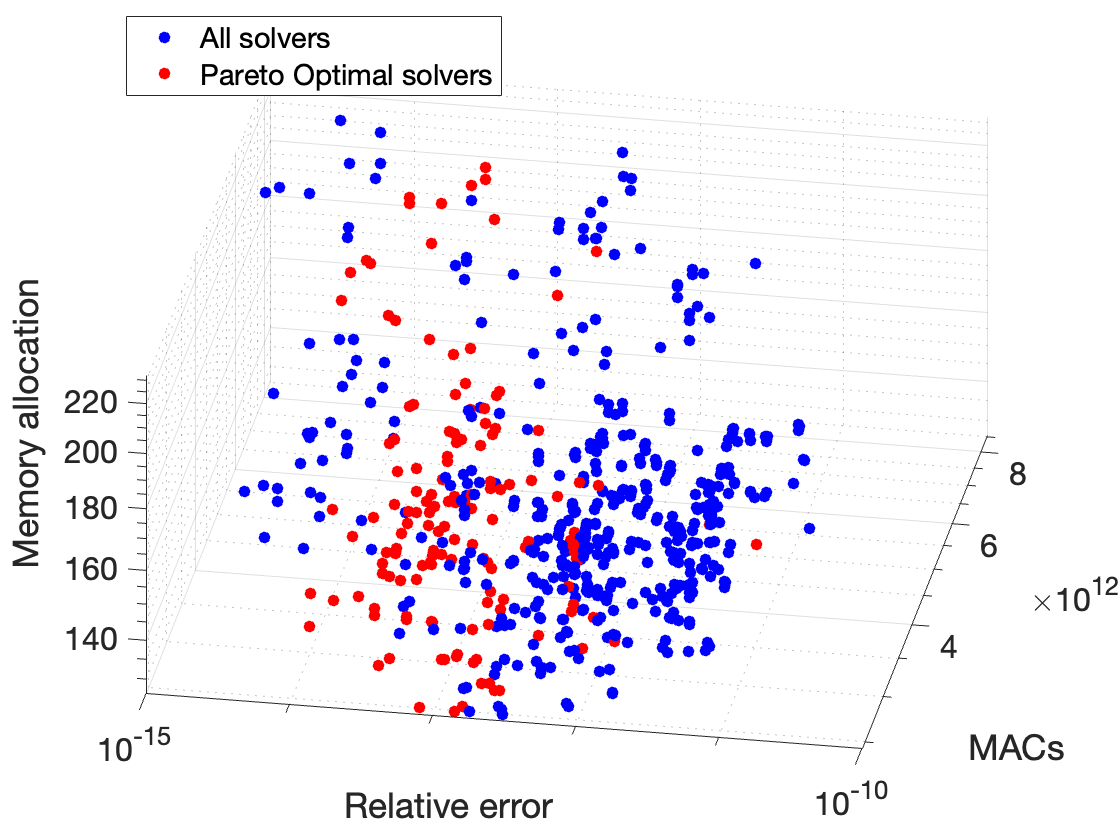}}
\subfigure[computational time -- relative error -- memory allocation]{
\includegraphics[width=0.45\textwidth]{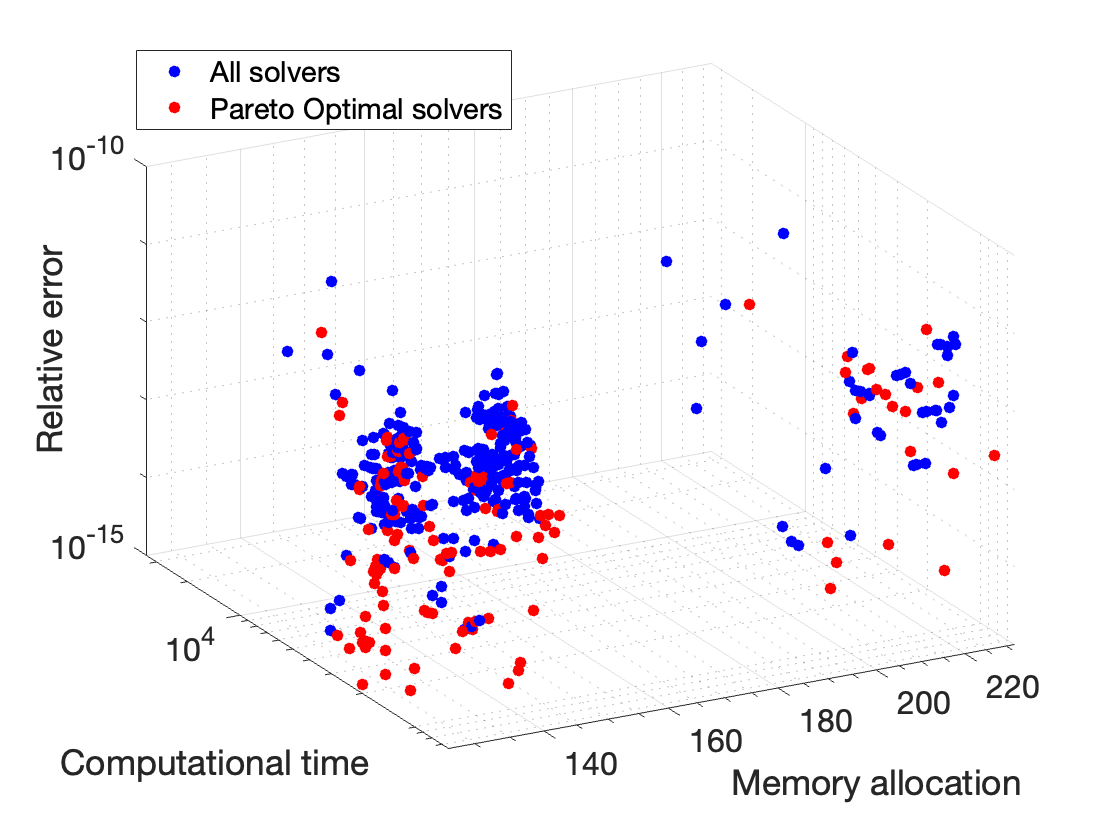}}
\subfigure[computational time -- MACs -- Average MACs]{
\includegraphics[width=0.45\textwidth]{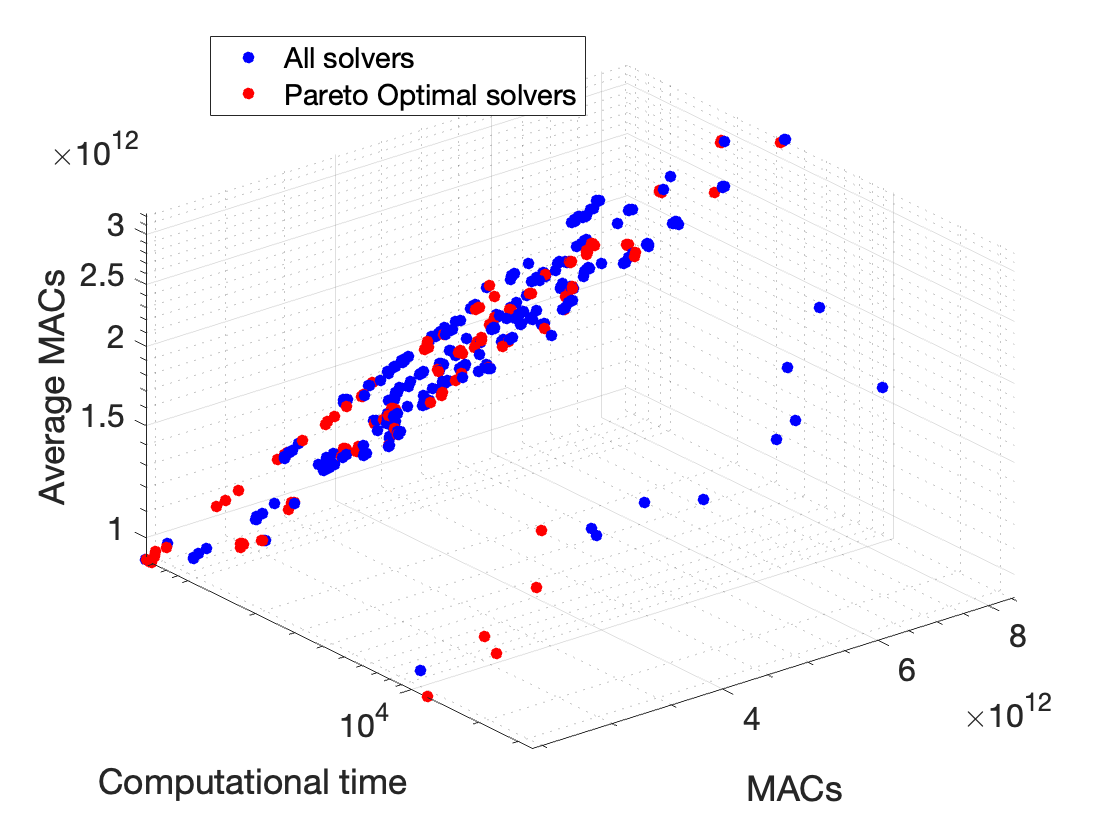}}
\caption{Pareto fronts for solving 2-D brittle fracture problem with coarse mesh, using Newton-Raphson method. All solvers are depicted in blue while
Pareto optimal solvers are highlighted in red.}
\label{front_bf_kry}
\end{figure}

\begin{figure}[H]
\centering
\subfigure[Computational time -- relative error -- number of iterations]{
\includegraphics[width=0.45\textwidth]{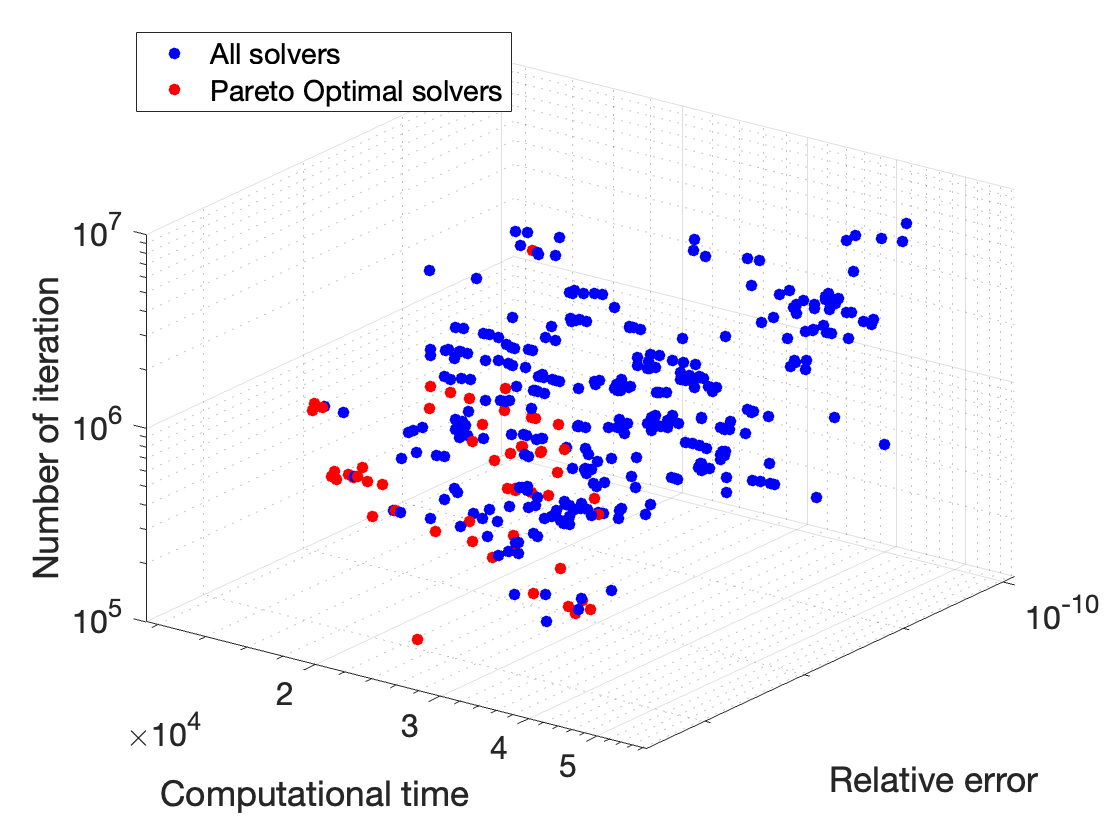}}
\subfigure[Relative error (p) -- memory allocation -- MACs]{
\includegraphics[width=0.45\textwidth]{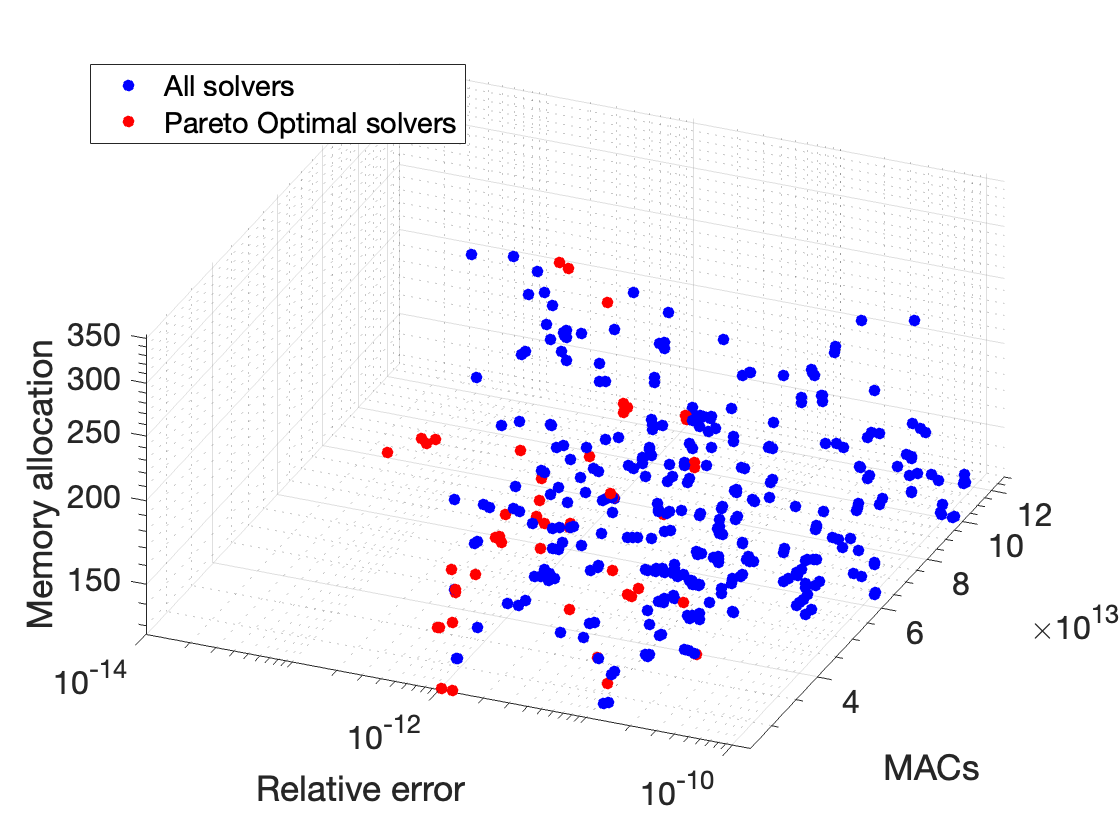}}
\subfigure[computational time -- relative error -- memory allocation]{
\includegraphics[width=0.45\textwidth]{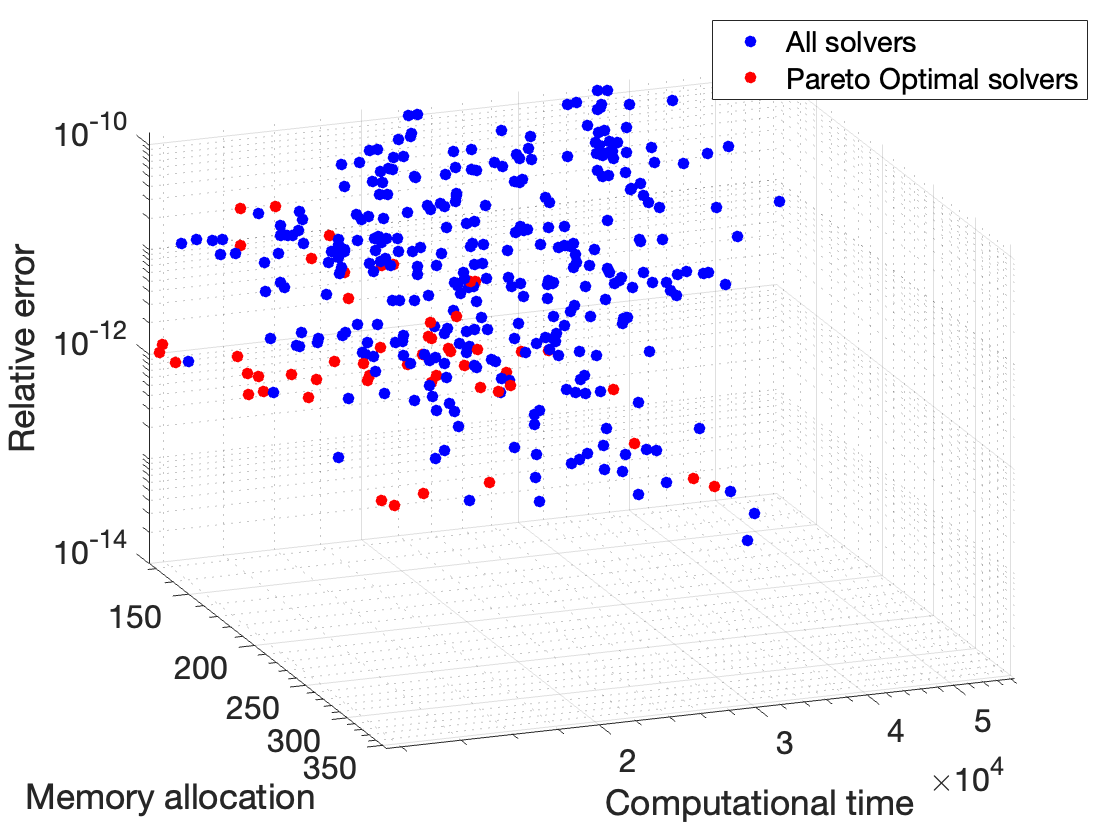}}
\subfigure[computational time -- MACs -- Average MACs]{
\includegraphics[width=0.45\textwidth]{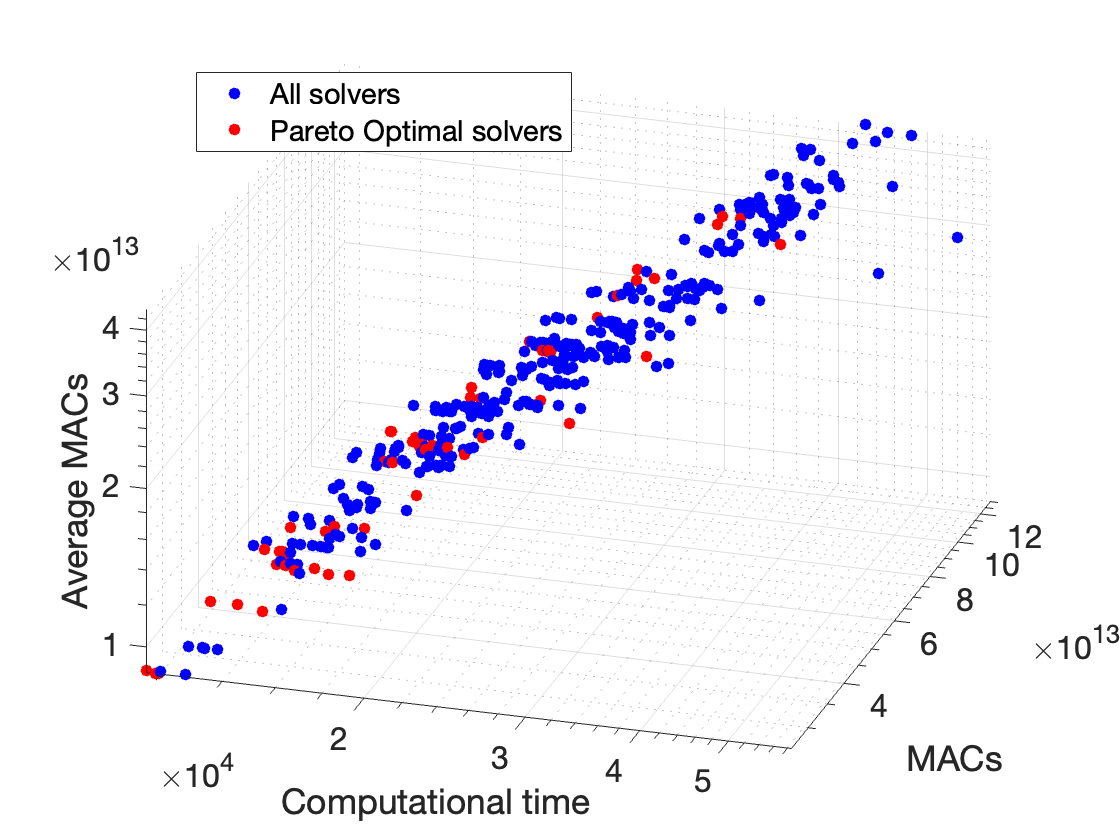}}
\caption{Pareto fronts for solving 2-D brittle fracture problem with fine mesh, using Newton-Raphson method. All solvers are depicted in blue while
Pareto optimal solvers are highlighted in red.}
\label{front_bf_kry_2}
\end{figure}

\subsection{Discovery of optimal meta-solvers by preference functions}

We implement the preference function that represents the average of all performance, for both coarse mesh and fine mesh. The top-3 meta-solvers are presented in the following: 
\begin{table}[H]
    \centering
    \caption{Top-3 solvers evaluated by the average preference function  for 2-d brittle fracture problem with coarse mesh  }
     \begin{minipage}{\linewidth}
        \centering
        \captionof{subtable}{The top-3 solvers}
        \resizebox{\linewidth}{!}{
        \begin{tabular}{c|c|c|c|c|c}
\hline \hline 
 &Neural operator & Classical solver  & Relaxation & Strategies & AMG \\ 
\hline \hline
Top 1 solver & U-DeepONet & CG  & Chebyshev & 3-1-3& NO \\ \hline
Top 2 solver & U-DeepONet & CG  & Chebyshev &  1-1-1& YES \\ \hline
Top 3 solver & U-DeepONet & CG  & Chebyshev & 5-1-5& NO \\ \hline
\end{tabular}
}
    \end{minipage}

    \vspace{1em} 

    \begin{minipage}{\linewidth}
        \centering
        \captionof{subtable}{Performance and rank of performance of the top-3 solvers}
        \resizebox{\linewidth}{!}{ 
          \begin{tabular}{c|c|c|c|c|c|c|c|c}
\hline \hline 
  &relative error:u& relative error:a& Com. time  & \# of ite & Memory & MACs&Ave. MACs & Training time \\ 
\hline 
Top 1 solver &$7.84 \times 10^{-14}$ & $1.77 \times 10^{-13}$ & 660.96 & 460328 & 126.856 & $2.49 \times 10^{12}$ & $9.45 \times 10^{11}$ & 2110.97 \\
\hline
Top 2 solver & $2.29 \times 10^{-13}$ & $1.41 \times 10^{-13}$ & 1147.43 & 305761 & 127.271 & $2.7 \times 10^{12}$ & $1.02 \times 10^{12}$ & 2110.97 \\
\hline
Top 2 solver & $3.21 \times 10^{-14}$ & $2.08 \times 10^{-13}$ & 746.842 & 341859 & 127.066 & $2.89 \times 10^{12}$ & $1.09 \times 10^{12}$ & 2110.97 \\
\hline
\end{tabular}
        }
    \end{minipage}
    \label{top3_p1_1d}
\end{table}

Compare with vanilla method, around 3.19 times faster in computational time and 69.1 times faster in number of iterations. 

\begin{table}[ht]
\centering
\begin{tabular}{c|c|c|c|c|c}
\hline \hline
 Vanilla Krylov method & Com. time & \# of iterations & Memory & MACs & Ave. MACs\\
\hline \hline
BiCGStab & 3930 & 21647901 & 102.958 & $2.09\times 10^{13}$ & $7.98 \times10^{12}$ \\ \hline
CG & 2110 & 21114235 & 103.419 & $1.09 \times 10^{13}$ & $4.13 \times10^{12}$ \\ \hline
\end{tabular}
\caption{Performance of vanilla method in coarse mesh}
\end{table}

\begin{table}[H]
    \centering
    \caption{Top-3 solvers evaluated by the average preference function  for 2-d brittle fracture problem with fine mesh  }
     \begin{minipage}{\linewidth}
        \centering
        \captionof{subtable}{The top-3 solvers}
        \resizebox{\linewidth}{!}{
        \begin{tabular}{c|c|c|c|c|c}
\hline \hline 
 &Neural operator & Classical solver  & Relaxation & Strategies & AMG \\ 
\hline \hline
Top 1 solver & U-DeepONet & CG  & Chebyshev & 1-1-1& YES \\ \hline
Top 2 solver & DeepONet & CG  & Chebyshev &  1-1-1& YES \\ \hline
Top 3 solver & U-DeepONet & BiCGStab  & Chebyshev & 1-1-1& YES \\ \hline
\end{tabular}
}
    \end{minipage}

    \vspace{1em} 

    \begin{minipage}{\linewidth}
        \centering
        \captionof{subtable}{Performance and rank of performance of the top-3 solvers}
        \resizebox{\linewidth}{!}{ 
          \begin{tabular}{c|c|c|c|c|c|c|c|c}
\hline \hline 
  &relative error:u& relative error:a& Com. time  & \# of ite & Memory & MACs&Ave. MACs & Training time \\ 
\hline 
Top 1 solver &$4.10 \times 10^{-6}$  & $1.23 \times 10^{-12}$ & 11937.4 & 521754 & 126.896 & $2.63 \times 10^{13}$ & $8.96 \times 10^{12}$ & 2110.97 \\
\hline
Top 2 solver & $4.10 \times 10^{-6}$  & $1.04 \times 10^{-12}$ & 11844.2 & 521426 & 126.744 & $2.63 \times 10^{13}$ & $8.95 \times 10^{12}$ & 2388.54 \\
\hline
Top 3 solver & $4.10 \times 10^{-6}$  & $4.77 \times 10^{-13}$ & 15300.9 & 363977 & 127.058 & $3.55 \times 10^{13}$ & $1.21 \times 10^{13}$ & 2110.97 \\
\hline
\end{tabular}
        }
    \end{minipage}
    \label{top3_p1_1d}
\end{table}

Compare with vanilla method, around 3.63 times faster in computational time and 158.72 times faster in number of iterations. 

\begin{table}[ht]
\centering
\begin{tabular}{c|c|c|c|c|c}
\hline \hline
 Vanilla Krylov method & Com. time & \# of iterations & Memory & MACs & Ave. MACs\\
\hline \hline
CG & 43000 & 57771176 & 89.1924 & $1.84\times 10^{14}$ & $6.22 \times10^{13}$ \\ \hline
\end{tabular}
\caption{Performance of vanilla method in fine mesh}
\end{table}

\end{document}